%% file: main.tex
\documentclass{siamart190516}
\usepackage{amsmath,amsfonts,amssymb,amscd,xcolor}
\usepackage{standalone}
\usepackage{tikz}
\usetikzlibrary{spy}
\usepackage{placeins}
\usepackage{adjustbox}
\usepackage{pgfplots}
\usepackage{pgfplotstable}
\usepgfplotslibrary{patchplots}
\pgfplotsset{compat=newest}
\usetikzlibrary{calc}
\usetikzlibrary{shapes.misc}
\usetikzlibrary{patterns.meta}
\usetikzlibrary{shadows.blur}
\usetikzlibrary{positioning}
\usepgfplotslibrary{colorbrewer}

\usepackage{pbox}
\usepackage{bm}
\DeclareMathSymbol{\shortminus}{\mathbin}{AMSa}{"39}

\usepackage{pdfpages}
\AtEndDocument{\includepdf[pages={-}]{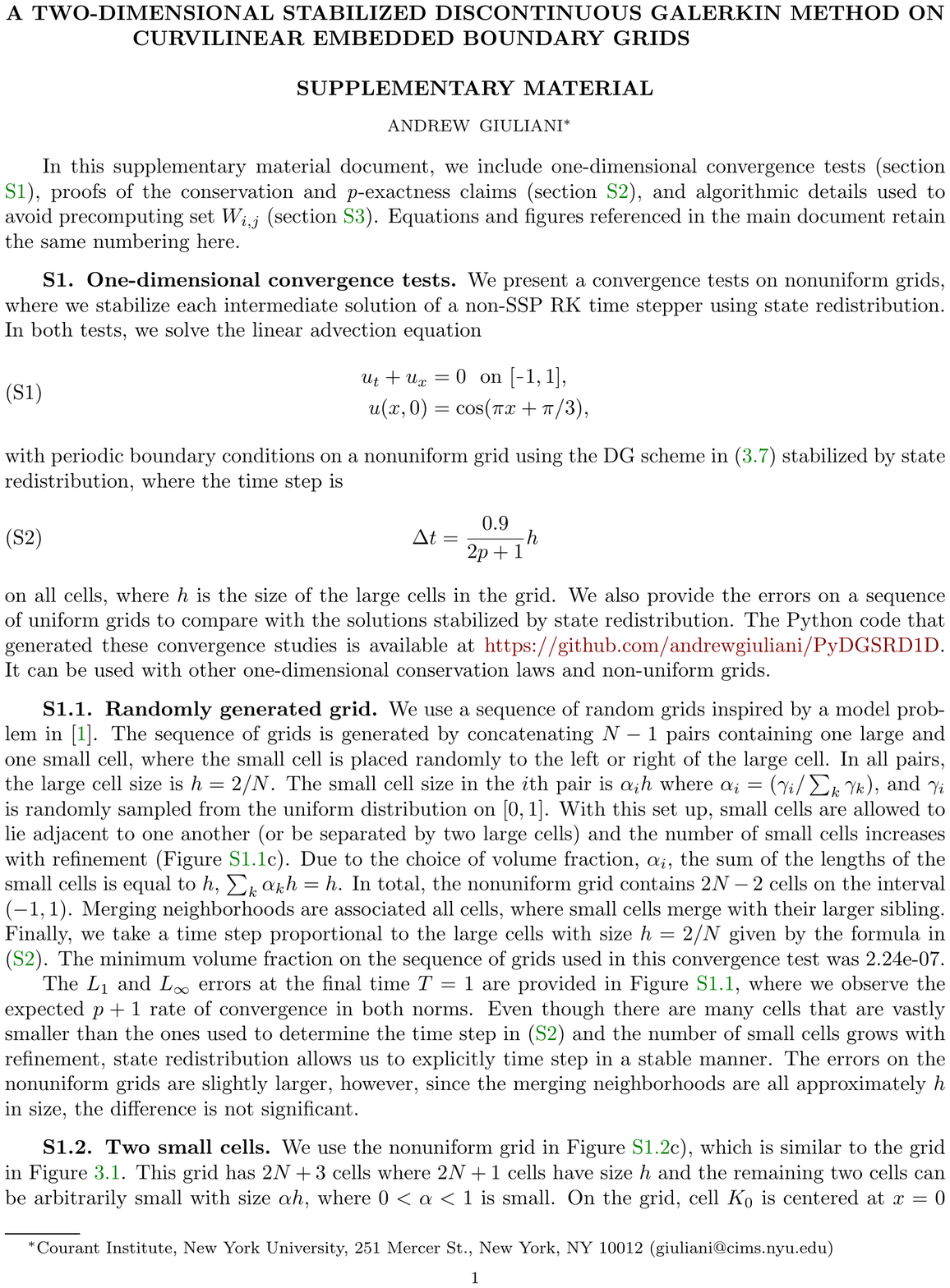}}

\title{A two-dimensional stabilized discontinuous Galerkin method on curvilinear embedded boundary grids}
\author{Andrew Giuliani\thanks{Courant Institute, New York University, 251 Mercer St., New York, NY 10012 (giuliani@cims.nyu.edu)}   }

\begin{document}
\maketitle
\begin{abstract}
    We propose a state redistribution method for high order discontinuous Galerkin methods on curvilinear embedded boundary grids.  
    State redistribution relaxes the overly restrictive CFL condition that results from arbitrarily small cut cells and explicit time stepping. 
    Thus, the scheme can take time steps that are proportional to the size of cells in the background grid. 
    The discontinuous Galerkin scheme is stabilized by postprocessing the numerical solution after each stage or step of an explicit time stepping method.
    This is done by temporarily merging the small cells into larger, possibly overlapping neighborhoods using a special weighted inner product.
    Then, the numerical solution on the neighborhoods is returned to the base grid in a conservative fashion.
    The advantage of this approach is that it uses only basic mesh information that is already available in many cut cell codes and does not require complex geometric manipulations.
    Finally, we present a number of test problems that demonstrate the encouraging potential of this technique for applications on curvilinear embedded geometries.
    Numerical experiments reveal that our scheme converges with order $p+1$ in $L_1$ and between $p$ and $p+1$ in $L_\infty$ for problems with smooth solutions.
    We also demonstrate that state redistribution is capable of capturing shocks.  
\end{abstract}

\begin{keywords}
state redistribution, discontinuous Galerkin methods, cut cell grids, explicit time stepping, curvilinear embedded boundaries.
\end{keywords}

\begin{AMS}
	65M60, 65M20, 35L02, 35L65
\end{AMS}

\section{Introduction}
Practical problems in computational fluid dynamics require the solution of hyperbolic conservation laws on complex domains.  
Typically, the domain is first discretized into a mesh of elements, then algorithms such as the finite volume (FV) or discontinuous Galerkin (DG) methods are used to approximate the exact solution on that mesh.  
When the boundary of the domain is complex, the most time-consuming aspect of the workflow can be the mesh generation phase \cite{osti_1149042}.
One technique that solves this problem is the use of embedded-boundary grids whereby the complex boundary is superimposed on a Cartesian grid.  
The computational mesh is then composed of Cartesian cells on the domain interior and irregular ``cut" cells on the domain boundary.
These mesh generation algorithms are robust and automated.  

A problem that arises when using explicit time stepping methods on embedded boundary grids is the small cell problem.  
Explicit methods require a time step that is proportional to the cell size in the mesh.  
Since irregular cut cells on the domain boundary can be arbitrarily small, this results in an overly stiff system of ordinary differential equations that would require impractically small time steps.  
A special treatment of these cut cells is therefore necessary for explicit time integrators to be used.

Many solutions to the small cell problem in finite volume methods have been proposed, such as flux redistribution \cite{colella2006cartesian}, $h$-box methods \cite{Berger2012ASH, mjb-hel-rjl:hbox2}, implicit/explicit time stepping \cite{JEBENS20121380,May-Berger:JSC}, dimensionally-split flux stabilization \cite{klein2009well,gokhale2018dimensionally}, cell merging \cite{quirk1994alternative}, and more recently state redistribution (SRD) \cite{berger2020state}.
Flux redistribution is an approach where each cell is updated with a locally stable, but possibly nonconservative, time step \cite{colella2006cartesian}.  Excess flux is then locally redistributed in order to maintain conservation.  This is straightforward to implement in three dimensions, but is only first order accurate at the embedded boundary.
\textit{H}-box methods are second order accurate, stable, and have desirable theoretical properties \cite{Berger2012ASH, mjb-hel-rjl:hbox2}.  This is achieved by increasing the domain of dependence of fluxes on cut cells, however, doing so in three dimensions seems difficult to implement.
Implicit/explicit time integration schemes have also been developed, where the small cells are integrated implicitly and the large cells are integrated explicitly \cite{May-Berger:JSC}.

Approaches that have been used for cut cell DG discretizations include subcycling in time \cite{doi:10.1002/nme.6343}, and penalization terms \cite{engwer2020stabilized}.
When strong anisotropic refinement is required, triangular cut cell approaches for the Navier-Stokes equations have also been examined in \cite{FIDKOWSKI20071653, fidkowski2007adaptive}.
A different technique is used in discontinuous Galerkin difference methods, where basis functions are no longer local to an element and are defined on a neighborhood of cells \cite{dgd_kaur}.
Additional penalty terms can also be designed to allow the neighbors of cut cells to communicate, as proposed in with the Domain of Dependence stabilization \cite{dod1, dod2, dod3}, where it can be shown that this results in a monotone scheme for piecewise constant solutions.
Alternatively, ghost penalty stabilization terms can be used \cite{fu2021high, ghost1}.
Cell merging is popular in two dimensions \cite{doi:10.1002/nme.6343, QIN201324,  krause2017incompressible,muller2017high, SAYE2017647, gulizzi2021coupled,gulizzi2021modelling} and is gaining more traction in three dimensions \cite{SAYE2017647, gulizzi2021coupled,gulizzi2021modelling}. However, this technique carries with it some design choices that could be intricate to implement on complex engineering geometries in three dimensions.

Our approach to the small cell problem is an extension of state redistribution method, developed for second order finite volume methods in two dimensions \cite{berger2020state}.  
State redistribution temporarily merges cut cells into larger, possibly overlapping neighborhoods in a postprocessing step applied after each explicit time step.  
Then, these neighborhood polynomials are recombined in a conservative manner back onto the base grid.
Applying state redistribution to DG numerical solutions follows a similar strategy, however, since the DG scheme does not rely on reconstructions, there are notable differences.
Specifically, the DG numerical solution is composed of discontinuous, locally defined polynomials that are written as linear combinations of polynomial basis functions.
Thus, we adapt our state redistribution approach to DG methods by associating to each possibly overlapping, merging neighborhood a tailored polynomial basis.
The DG solution is projected onto these overlapping merging neighborhoods using a particular weighted inner product.
Then, the numerical solution is projected back onto basis functions defined on elements of the embedded boundary grid. 

The structure of this article is as follows. In section \ref{sec:dgmethod} we briefly describe the spatial discretization used to approximate solutions to hyperbolic conservation laws. In section \ref{sec:srd1d}, we recall the state redistribution method for finite volume schemes and outline how to extend it to arbitrarily high order discontinuous Galerkin methods.  This will involve the projections that we discussed previously.  
For ease of understanding and to reduce the necessary notation, this section uses a model problem in one space dimension.  
In sections \ref{sec:preprocessing} and \ref{sec:srd2d}, we go into more detail for two-dimensional computations.  
In section \ref{sec:numex}, we present a number of numerical experiments. In particular, we examine the error at cut cells, confirming the results in \cite{QIN201324,engwer2020stabilized} that the cut cell methods suffer some loss of accuracy at the boundary.

\section{The discontinuous Galerkin method on cut cell meshes} \label{sec:dgmethod}
We approximate solutions to hyperbolic conservation laws
\begin{equation}\label{eq:hyp}
    \textbf{u}_t + \nabla \cdot \mathbf{F}(\textbf{u}) = 0 ~\text{ on } \Omega \times (0,T),
\end{equation}
using the discontinuous Galerkin method coupled with explicit time steppers on cut cell meshes.  
In the above, $\Omega \subset \mathbb{R}^2$ is the spatial domain, $T$ is the final time, $\textbf{u}: \Omega \times (0,T) \rightarrow \mathbb{R}^m$, $m \in \mathbb{N}$ is a vector of conserved quantities and $\mathbf{F}$ is a flux function.

\subsection{Mesh generation} \label{sec:genmesh}
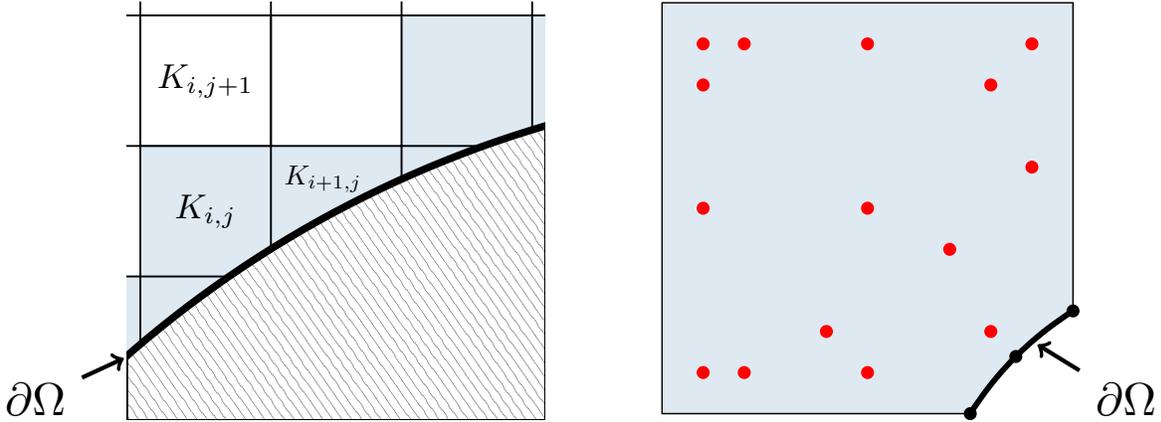
\begin{figure}
\centering
        \begin{tikzpicture}
        \node (11) at (-8,0) 
        {\includestandalone[width = 0.35\textwidth]{1-dg/example_mesh}};
        
        \node [scale=1.75] (A) at (-12, -2.5) {$\partial \Omega$};
        \node (B) at (-10.7, -1.9) {};
        \draw[->,ultra thick]  (A) -- (B);
        
        \node (12) at (-.9,0){\includestandalone[width = 0.35\textwidth]{1-dg/cc}};
        \node [scale=1.75] (AA) at (2.5, -2.5) {$\partial \Omega$};
        \node (BB) at (1.2, -1.7) {};
        \draw[->,ultra thick]  (AA) -- (BB);
        \end{tikzpicture}
        \caption{
        On the left, we plot a zoom of a curvilinear cut cell grid at the boundary, where whole cells are white and cut cells are light blue.  The boundary of the domain $\partial \Omega$ is the bold black line.
        On the right, we show the quadrature rule for $p=2$ of order accuracy 4 on $K_{i,j}$. The intersection of the cell with the embedded boundary are indicated by black circles ($\bullet$) and the volume quadrature points of the cut cell are red circles ($\textcolor{red}{\bullet}$). 
        } \label{fig:mesh}
\end{figure}

The domain $\Omega$ is discretized with a mesh of elements $K_{i,j}$, by superimposing the boundary of $\partial \Omega$ on a background Cartesian grid and deleting the exterior of the domain from the grid.  
The geometry in two dimensions is defined as a function that takes a point $(x,y)$ and returns whether it is inside or outside the domain.
In the volume, the mesh comprises regular Cartesian cells of dimension $\Delta x$ by $\Delta y$.
On the boundary, there are irregular polygons with both straight and curved edges that are determined by the intersections of the boundary, $\partial \Omega$, with the background Cartesian grid.
The $x$ and $y$ coordinates of the curved edge are polynomials of degree $q$ that interpolate $\partial \Omega$ at $q+1$ points.
For example, the cut cell on the right in Figure \ref{fig:mesh} has one curved edge, where $q = 2$.
We compute interpolation points that are equispaced in arclength on curved boundaries.  For circular boundaries, this is done exactly since closed form expressions are available.  The approach for non-circular boundaries is approximate and relies on a parametric representation of the boundary $\Gamma(s) = (x(s),y(s))$.  Say that a cut cell's curved edge is defined by the set of points $\{ \Gamma(s) : s \in [s_1,s_2] \}$.
In order to find additional interpolation points along the boundary, define the function 
$$
f(s) = \frac{||\Gamma(s) - \Gamma(s_1))||_2}{||\Gamma(s_2)-\Gamma(s_1)||_2} - w
$$ 
for some $s \in [s_1,s_2]$ and $w \in (0,1)$.  This function is zero when $s$ maps to a point, $\Gamma(s)$, that is exactly $w || \Gamma(s_2) - \Gamma(s_1)||_2$ units away from $\Gamma(s_1)$.  The $q+1$ interpolation points along the boundary are computed by finding the roots of $f$ when $w = i/q \text{ for } i = 1 \hdots q-1$.  These $q-1$ additional boundary points, along with the edge endpoints $\Gamma(s_1)$ and $\Gamma(s_2)$, form a set of $q+1$ edge interpolation points. 
Note that this approach does not extend to three dimensions and is only appropriate when the boundary is sufficiently well-resolved.

As reported in \cite{bassi1997high, KRIVODONOVA2006492}, discontinuous Galerkin methods require at least piecewise quadratic approximations to curved boundaries when solving the Euler equations.  
In \cite{bassi1997high}, it was found that a nonphysical boundary layer would develop and not disappear with refinement when using a second order accurate DG scheme ($p=1$) and piecewise linear approximations to the curved boundary ($q=1$).  This numerical difficulty resolved when a piecewise quadratic representation of the boundary was used, pairing $p=1$ and $q=2$.  No problems were found for $p=q$ when $p>1$. 
Thus, we always pair the second order DG method ($p=1$) with a piecewise quadratic approximation of the boundary ($q=2$).  
Otherwise, we let $p=q$ when $p>1$.
In numerical experiments (section \ref{sec:numex}), we sometimes use the notation, e.g. DG-$P5Q5$, to refer to simulations where the number after the $P$ refers to the polynomial degree of accuracy on the cells ($p$), and the number after $Q$ refers to the polynomial degree of the boundary interpolant ($q$).
Each cut cell is assumed to have at most two, possibly curved, irregular edges associated to the embedded boundary ($\partial \Omega$). For example, the cut cells in Figure \ref{fig:mesh} (highlighted in blue), only have one irregular edge.  Cells that have two irregular edges lie on sharp corners of the embedded boundary (Figure \ref{fig:overlaps_ringleb}).  
We also do not allow split or tunneled cells for ease of code development.
 This restriction is not fundamental and we will need to implement these features in three dimensions for complicated geometries.
We have made our mesh generation code written in Python available at \url{https://github.com/andrewgiuliani/PyGrid2D}.  
This code can be used to reproduce all the high order embedded boundary grids used in this work, as well as generate other embedded boundary grids in two dimensions.

\subsection{Semidiscretization}
We use the modal discontinuous Galerkin method as the base scheme on this cut cell grid.  
Multiplying \eqref{eq:hyp} by a test function $v(x,y) \in H^1(K_{i,j})$ and integrating on cell $K_{i,j}$, we obtain
$$
\int_{K_{i,j}} \textbf{u}_t v ~dx ~dy + \int_{K_{i,j}} \nabla \cdot \mathbf{F}(\textbf{u}) v ~dx ~dy = 0 ~~\forall v \in H^1(K_{i,j}).
$$
Integrating the second term by parts and using the divergence theorem, we have
$$
\int_{K_{i,j}} \textbf{u}_t v ~dx ~dy + \int_{\partial K_{i,j}} v \mathbf{F}(\textbf{u})\cdot \mathbf{n} ~dl - \int_{K_{i,j}} \mathbf{F}(\textbf{u}) \cdot \nabla v ~dx ~dy = 0 ~~\forall v \in H^1(K_{i,j}).
$$
The solution $\mathbf{u}$ and test function $v$ are then restricted to $S^p(K_{i,j})$,  the space of polynomials of degree less than or equal to $p$ on $K_{i,j}$.
On each cell, Cartesian and cut, we use a non-tensor product, orthonormal basis for $S^p(K_{i,j})$ that has $N_p + 1$ basis functions, where $N_p = (p+1)(p+2)/2-1$.  
Alternatively, a tensor product basis \cite{QIN201324} can be used.
Tensor and non-tensor product bases have respectively $(p+1)^2$ and $(p+1)(p+2)/2$ basis functions.
The choice of orthonormal basis functions was a design decision that results in a diagonal mass matrix for convenience.
We prefer a non-tensor product basis (section \ref{sec:ccbasis}) due to the computational savings that result from a smaller set of basis functions.
The numerical solution on each cell is written
\begin{equation} \label{eq:sol}
\textbf{U}_{i,j}(x,y,t) = \sum_{k = 0}^{N_p} \textbf{c}_{i,j,k}(t)\varphi_{i,j,k}(x,y),
\end{equation}
where $\textbf{c}_{i,j,k} \in \mathbb{R}^m$ is the $k$th solution coefficient associated to basis function $\varphi_{i,j,k}$ on cell $K_{i,j}$, where the $\varphi_{i,j,k}$ are orthonormal with respect to the inner product
\begin{equation} \label{eq:l2}
\langle f,g \rangle_{K_{i,j}} = \frac{1}{|K_{i,j}|} \int_{K_{i,j}} fg ~dx~dy,
\end{equation}
where $|K_{i,j}|$ is the volume of cell $i,j$. 
Normalizing the inner product by the cell volume results in the constant basis function $\varphi_{i,j,0} = 1$, independent of the cell geometry.
Without this normalization, the constant basis function on the cell is $\varphi_{i,j,0} = 1/\sqrt{|K_{i,j}|}$, which would blow up on the small cut cells.

The modal discontinuous Galerkin method is then given by the system of ordinary differential equations
\begin{equation} 
    \frac{d}{dt}\textbf{c}_{i,j,k} = -\frac{1}{|K_{i,j}|} \left[  \int_{\partial K_{i,j}} \varphi_{i,j,k}\mathbf{F}^*(\textbf{U}^-_{i,j}, \textbf{U}_{i,j}^+)\cdot \mathbf{n} ~dl - \int_{K_{i,j}}\mathbf{F}(\textbf{U}_{i,j})\cdot \nabla \varphi_{i,j,k} ~dx~dy \right],  \label{eq:dg}
\end{equation}
for $k = 0, \hdots, N_p$, where $\textbf{U}_{i,j}^-$ and $\textbf{U}_{i,j}^+$ are the numerical solutions corresponding to $K_{i,j}$, and a neighboring cell, respectively, located along the element boundary $\partial K_{i,j}$.  Finally, $\mathbf{n}$ is an outward facing unit normal, and $\mathbf{F}^*$ is a numerical flux function.

The semidiscretization \eqref{eq:dg} is integrated in time using a $(p+1)$th order accurate Runge Kutta time stepping algorithm.
Since the cut cell grid can contain arbitrarily small cut cells, \eqref{eq:dg} can be stiff, and be subject to an arbitrarily small time step restriction.  
State redistribution will allow the use of explicit time stepping algorithms on cut cell grids with a time step that is proportional to the Cartesian cell size
\begin{equation} \label{eq:cfl}
\Delta t \left(\frac{|a|}{\Delta x} + \frac{|b|}{\Delta y} \right) \leq \frac{1}{2p+1},
\end{equation}
where $|a|,|b|$ are the maximum wave speeds in the $x$ and $y$ directions \cite{COCKBURN1998199, qinkrivo}.

The DG method above requires an orthogonal polynomial basis on each polygonal cut cell.  We numerically compute orthogonal basis functions and their gradients at quadrature points using the QR factorization of a Vandermonde-like matrix (section \ref{sec:ccbasis}).

\subsection{Quadrature on curved polygons} \label{sec:quadrature}
One issue with cut cell DG approaches that arises independently of the small cell problem is the evaluation of inner products on arbitrary curved polygons.  
Quadrature rules are required to approximate volume integrals on elements
$$
\int_{K_{i,j}} f(x,y)g(x,y) ~dx~dy \approx \sum^{N^{\text{vol}}_{i,j}}_{q = 0} w_{i,j,q} f(x_{i,j,q},y_{i,j,q})g(x_{i,j,q},y_{i,j,q}),
$$
where $N^{\text{vol}}_{i,j}+1$ is the number of volume quadrature points on $K_{i,j}$, and $\sum^{N^{\text{vol}}_{i,j}}_{q = 0} w_{i,j,q} = |K_{i,j}|$.
Some approaches are based on the divergence theorem \cite{sudhakar2014accurate}, hierarchical moment fitting \cite{doi:10.1002/nme.4569, doi:10.1002/nme.5288}, ``speckling" of quadrature points \cite{FIDKOWSKI20071653, fidkowski2007adaptive}, and others \cite{mousavi2010generalized,saye2015high}.
We use the software described in \cite{Sommariva2020,quad_codes} for quadrature rules on polygons with curved edges, but other quadrature generation codes can be used as well.  By construction, this algorithm always results in quadrature nodes interior to the polygon with positive weights.
Another well-known approach is to triangulate the polygon and apply standard quadrature rules on each subtriangle \cite{QIN201324}.

Variational crimes due to quadrature errors in under-resolved solutions can lead to instabilities for nonlinear conservation laws \cite{hesthaven2007nodal}.   
Despite this, on cut and Cartesian cells, we use quadrature rules of order accuracy $2p$ and $2p+1$, respectively, as we did not encounter these problems in the presented numerical experiments.  We note that increasing the order of accuracy of the quadrature rules can mitigate these issues should they occur \cite{gassnerwinters}.
Surface integrals are approximated using standard Gauss-Legendre rules of order accuracy $2p+1$, with $p+1$ points.

\section{State redistribution in one dimension}\label{sec:srd1d}
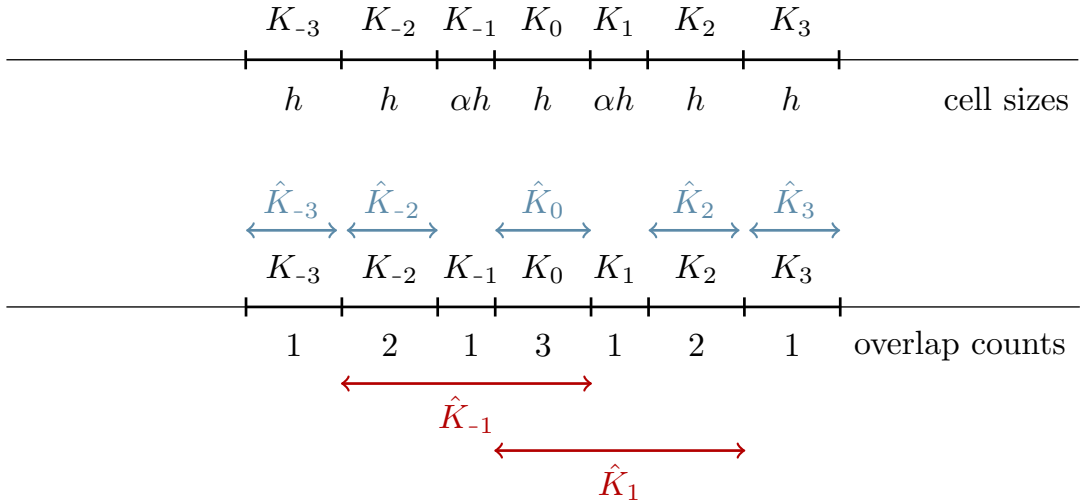
\begin{figure}
    \centering
    \begin{tikzpicture}
    \node (11) at (-4, 0)
    {\includestandalone[width=0.9\textwidth]{2-srd/nonuniform1}};
    \node (12) at (-4, -3.75)
    {\includestandalone[width=0.9\textwidth]{2-srd/nonuniform2}};
    \end{tikzpicture}
    
    \caption{Nonuniform grid for model problem used to describe the state redistribution method for both finite volume and discontinuous Galerkin methods in Section \ref{sec:srd1d}.  
    The top figure indicates the cell sizes.
    The bottom figure shows the merging neighborhoods and the overlap counts on each cell in the base grid.
    $\hat K_i$ refers to the merging neighborhood associated to cell $K_i$ in the base grid.
    Neighborhoods that contain only one cell in the base grid are drawn in blue, while those that contain multiple, temporarily merged cells are drawn in red.
    }
    \label{fig:nonunform}
\end{figure}

In this section, we demonstrate state redistribution with a simple example that solves the linear advection equation
\begin{equation}\label{eq:hpde}
u_t + au_x = 0, \quad a>0
\end{equation}
on the nonuniform grid in Figure \ref{fig:nonunform}, called the base grid.   On full cells, $h_i = h$ and on
the small cells, $h_i = \alpha h$ for $ 0 < \alpha < 1$.
Five cells (indexed by $\shortminus 3$, $\shortminus 2$, $0$, $2$, $3$) are large with size $h$ and the remaining two cells (indexed by $\shortminus 1$ and $1$) 
are small with size $\alpha h$.
For completeness, we first  recall how to stabilize a first order upwind finite volume scheme using the approach described in \cite{berger2020state}.  Note that this scheme can be viewed as a discontinuous Galerkin method when $p=0$.
Then, we will describe how to extend this approach to arbitrarily high order discontinuous Galerkin methods.

\subsubsection*{State redistribution with finite volume methods}
We discretize \eqref{eq:hpde} using the first order, upwind finite volume scheme
\begin{equation}
\begin{aligned}\label{eq:unstable1d_fv}
\hat{U}_{i} &= U^n_{i} - \frac{a \Delta t} {h_i}  (U^n_{i} - U^n_{i-1}).
\end{aligned}
\end{equation}
To apply state redistribution to the above finite volume scheme, we must execute a simple geometric preprocessing stage on the grid before the explicit time stepping portion of the finite volume solver.  We associate a merging neighborhood $\hat K_{i}$ to each cell in the base grid, $K_{i}$.
Cells that are smaller than half the regular grid size ($h_i < h/2$), are merged to the left and right to form the merging neighborhoods $\hat K_{\shortminus 1} = K_{\shortminus2} \cup K_{\shortminus1} \cup K_{0}$ and $\hat K_{1} = K_{0} \cup K_{1} \cup K_{2}$ (Figure \ref{fig:nonunform}, in red).
Cells that are larger than half the regular grid size ($h_i \geq h/2$) do not require merging and we have $\hat K_{i} = K_i$ (Figure \ref{fig:nonunform}, in blue).
Additionally, each cell in the base grid must count the number of neighborhoods that overlap it, $N_i$.  For example, since $K_0$ is a member of three neighborhoods: $\hat K_{\shortminus1}$, $\hat K_{0}$, and $\hat K_{1}$, its overlap count is $N_0 = 3$.  We have provided the overlap counts for the cells on the base grid in Figure \ref{fig:nonunform}.

State redistribution is applied as a postprocessing step after each forward Euler update of \eqref{eq:unstable1d_fv}.  On each neighborhood, we compute a special weighted solution average.  For example, on $\hat K_{\shortminus1}$ the weighted solution average is
\begin{equation}\label{eq:qm1_fv}
    \hat Q_{\shortminus1} = \frac{1}{h/2 + \alpha h + h/3} \left( \frac{h}{2}\hat U_{\shortminus2} + \alpha h\hat U_{\shortminus1} + \frac{h}{3}\hat U_{0} \right),
\end{equation}
where both cell size and solution average in the above sums are divided by their associated overlap counts.  
Similarly, on $\hat K_{1}$ the weighted solution average is
$$
\hat Q_{1} = \frac{1}{h/2 + \alpha h + h/3} \left( \frac{h}{2}\hat U_{2} + \alpha h\hat U_{1} + \frac{h}{3}\hat U_{0} \right).
$$
On neighborhoods that only contain one cell, the weighted solution average is
$$
\hat Q_i = \hat U_i ~\text{ for } i = \shortminus 3, \shortminus 2, 0, 2, 3.
$$
Finally, these weighted solution averages are recombined back onto the base grid.
The final update of a cell in the base grid is given by the average of the neighborhood values that overlap it.
For example, since $K_0$ is overlapped by $\hat K_{\shortminus1}$, $\hat K_{0}$, and $\hat K_{1}$, we have
\begin{equation} \label{eq:final1}
U_0^{n+1} = \frac{1}{3}(\hat Q_{\shortminus1} + \hat Q_{0} + \hat Q_{1}).
\end{equation}
$K_{\shortminus2}$ and $K_2$ are overlapped by two neighborhoods, so their final updates are
\begin{equation} \label{eq:final2}
U^{n+1}_{\shortminus2} = \frac{1}{2}(\hat Q_{\shortminus2} + \hat Q_{\shortminus1}) \text{ and } U^{n+1}_{2} = \frac{1}{2}(\hat Q_{2} + \hat Q_{1}) .
\end{equation}
On the cells overlapped by only one neighborhood, the final update is simply
\begin{equation} \label{eq:final3}
U^{n+1}_{i} = \hat Q_{i} \text{ for } i = \shortminus 3, \shortminus 1, 1, 3.
\end{equation}
We prove in \cite{berger2020state} that the above scheme in \eqref{eq:final1}-\eqref{eq:final3}, is conservative.  Additionally, we show that it can be extended to second order FV methods by reconstructing a slope on each merging neighborhood and using that slope to recombine the neighborhood averages back onto the base grid in a linearity preserving fashion.
In the following section, we describe how to extend this idea to discontinuous Galerkin methods.

\subsubsection*{State redistribution with discontinuous Galerkin methods}

We discretize \eqref{eq:hpde} in space using an upwind RK-DG scheme where the degree of the polynomial approximation is $p$.
Assuming a strong stability preserving Runge Kutta (SSP-RK) scheme, the first forward Euler update can be written
\begin{equation}
\begin{aligned}\label{eq:unstable1d}
\hat{c}_{i,k} = c^n_{i,k} - \frac{a \Delta t} {h_i} \biggl[ U^n_{i}(x_{i+1/2})\varphi_{i,k}(x_{i+1/2}) &- U^n_{i-1}(x_{i-1/2})\varphi_{i,k}(x_{i-1/2}) \\
&- \int_{K_i} U^n_i \varphi_{i,k}' ~dx \biggr], ~\text{ for } k = 0 \hdots p,
\end{aligned}
\end{equation}
where $\hat c_{i,k}$ are the provisionally updated solution coefficients, $U^n_i$ is the DG solution on cell $K_i$ at time $t^n$ and $\hat U_i$ is its provisionally updated (possibly unstable) counterpart
$$
U^n_i(x) = \sum_{k = 0}^{N_p} c^n_{i,k} \varphi_{i,k}(x), \quad \hat U_i(x) = \sum_{k = 0}^{N_p} \hat c_{i,k} \varphi_{i,k}(x).
$$
Additionally, the time step is $\Delta t = h/a(2p+1)$ on all cells, and $\varphi_{i,k}$ are the orthonormal basis functions of $S^p(K_i)$ with respect to the inner product in \eqref{eq:l2}.

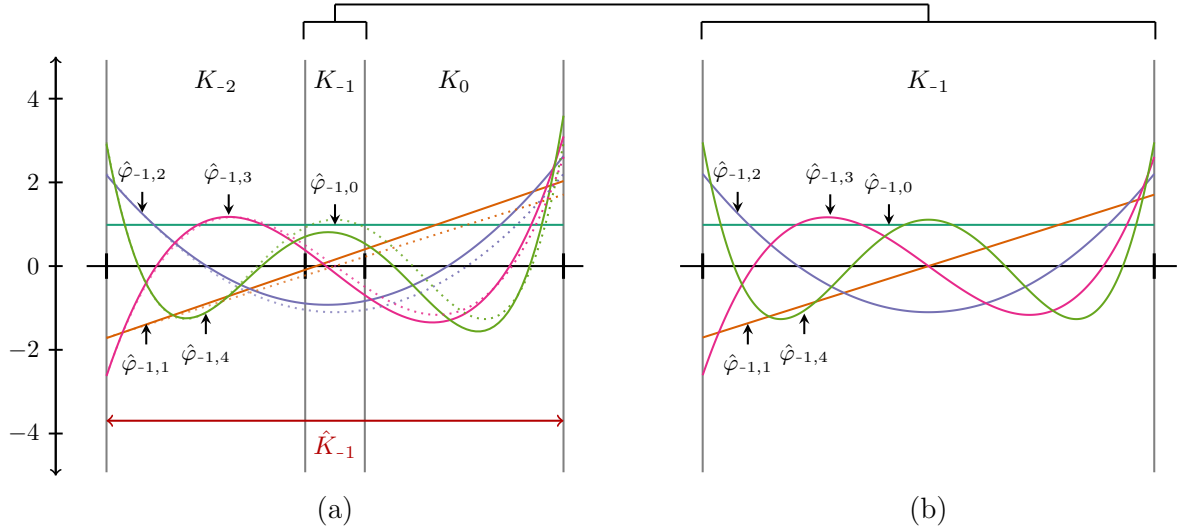
\begin{figure}
\centering
\resizebox{\columnwidth}{!}{%
    \begin{tikzpicture}
    \draw[<->, thick] (-8,-3)--(-8,3) ;
    \draw[-, thick]   (-8,-3)--(-8,3) ;
    
    \draw[-, thick]   (-8.1,-2.4)--(-7.9,-2.4)  node[left,xshift=-.2cm]{$-4$};
    \draw[-, thick]   (-8.1,-1.2)--(-7.9,-1.2)  node[left,xshift=-.2cm]{$-2$};
    \draw[-, thick]   (-8.1, 0.0)--(-7.9, 0.0)  node[left,xshift=-.2cm]{$0$};
    \draw[-, thick]   (-8.1, 1.2)--(-7.9, 1.2)  node[left,xshift=-.2cm]{$2$};
    \draw[-, thick]   (-8.1, 2.4)--(-7.9, 2.4)  node[left,xshift=-.2cm]{$4$};
    
    \node (A) at (-4, -3.5){\textcolor{black}{\large (a)} };
    \node (B) at (4.5,-3.5){\textcolor{black}{\large (b)} };
    
    \coordinate (m1lb) at (-4.45,3.25);
    \coordinate (m1rb) at (-3.55,3.25);
    \coordinate (m1l) at (-4.45,3.5);
    \coordinate (m1r) at (-3.55,3.5);
    \draw [thick] (m1l) -- (m1r);
    \draw [thick] (m1lb) -- (m1l);
    \draw [thick] (m1rb) -- (m1r);
    
    \coordinate (m1lzb) at (1.25,3.25);
    \coordinate (m1rzb) at (7.75,3.25);
    \coordinate (m1lz) at (1.25,3.5);
    \coordinate (m1rz) at (7.75,3.5);
    \draw [thick] (m1lz) -- (m1rz);
    \draw [thick] (m1lzb) -- (m1lz);
    \draw [thick] (m1rzb) -- (m1rz);
    
    \coordinate (ave) at ($(m1l)!0.5!(m1r)+(0,0.25)$);
    \coordinate (avez) at ($(m1lz)!0.5!(m1rz)+(0,0.25)$);
    \coordinate (aveb) at ($(m1l)!0.5!(m1r)$);
    \coordinate (avezb) at ($(m1lz)!0.5!(m1rz)$);
    \draw [thick] (ave) -- (aveb);
    \draw [thick] (avez) -- (avezb);
    \draw [thick] (ave) -- (avez);

    \node (11) at (-4,0)  {\includestandalone[width=0.45\textwidth]{2-srd/ortho_m0}};
    \node (12) at (4.5,0) {\includestandalone[width=0.45\textwidth]{2-srd/ortho_m1}};
    \end{tikzpicture}
}
    \caption{In a), we plot the basis functions (solid lines) that are orthogonal with respect to the weighted inner product $\langle \cdot, \cdot \rangle_{\hat K_{\shortminus1}}$ in \eqref{eq:ip1}.
    Although the basis functions in a) appear to be translations and rescaling of the standard Legendre polynomials, a closer inspection reveals that they are not.
    We have plotted the normalized Legendre polynomials (dotted lines) on the same axes to illustrate the difference.
    In b), we plot the orthogonal basis functions with respect to $\langle \cdot, \cdot \rangle_{K_{\shortminus 1}}$, which are simply translations and rescaling of the standard Legendre polynomials.  
    Note that functions in b) are plotted on element $K_{\shortminus1}$, which is the center element in a).}
    \label{fig:basis_functions}
\end{figure}

Similar to the SRD stabilized FV scheme described earlier, we must complete a preprocessing stage, whereby merging neighborhoods are associated to each cell in the base grid.  For the illustration here, we use the same merging neighborhoods as in the finite volume example.
The postprocessing stage described below is more involved for the DG scheme due to the spatial accuracy required and comprises two projections: 
\begin{enumerate}
    \item  The first step projects the provisionally updated numerical solution from the base grid, $\hat U_i$, onto basis functions defined on merging neighborhoods. They are called, $\hat \varphi_{i,k}$ and they span $S^p(\hat K_i)$.  
    \item The second step projects the numerical solution from the merging neighborhoods back onto the basis functions defined on the base grid.  They are called $\varphi_{i,k}$, and they span $S^p(K_i)$.
\end{enumerate}
We emphasize that the neighborhood basis functions are not generally the same as those on the base grid.  In this example, $\varphi_{i,k}$ and $\hat \varphi_{i,k}$ are not identical, and are defined on different domains.
The basis functions on the base grid $\varphi_{i,k}$ are orthogonal with respect to the inner product in \eqref{eq:l2}.
In contrast, the basis functions defined on the merging neighborhoods $\hat \varphi_{i,k}$ are orthogonal with respect to a particular weighted $L_2$ inner product.  
For example, $\hat \varphi_{\shortminus1,k}$ is the $k$th orthonormal basis function of $S^p(\hat K_{\shortminus1})$ with respect to the inner product
\begin{equation} \label{eq:ip1}
\langle f, g \rangle_{\hat K_{\shortminus1}} = \frac{1}{h/2 + \alpha h + h/3} \left( \frac{1}{2}\int_{K_{\shortminus2}} f~g~dx +  \int_{ K_{\shortminus1}} f~g~dx +  \frac{1}{3}\int_{ K_{0}} f~g~dx \right),
\end{equation}
where the cell sizes are divided by their associated overlap counts.
Additionally, each integral is divided by the overlap count associated to its domain of integration: the first integral is divided by two since $N_{\shortminus2} = 2$, the second integral is divided by one since $N_{\shortminus1} = 1$, and the final integral is divided by three since $N_0=3$.
In Figure \ref{fig:basis_functions}, we plot both $\varphi_{\shortminus1,k}$ and $\hat \varphi_{\shortminus1,k}$ for $k = 0,1, \hdots, 4$ and make two observations.  First, $\varphi_{\shortminus1,k}$ and $\hat \varphi_{\shortminus1,k}$ are defined on different domains.  Second, the basis functions $\varphi_{\shortminus1,k}$ on the base grid are simply translations and rescalings of the Legendre polynomials, as expected, while the neighborhood basis functions $\hat \varphi_{\shortminus1,k}$ are not.  

We are now ready to compute the projection $\hat Q_{\shortminus1}(x)$ of the provisionally updated DG solutions, $\hat U_{\shortminus2}$, $\hat U_{\shortminus1}$, $\hat U_0$, onto the merging basis for $\hat K_{\shortminus1}$:
\begin{equation}\label{eq:qhatm1}
\hat Q_{\shortminus1}(x) = \sum^{N_p}_{k = 0} \hat q_{i,k} \hat \varphi_{\shortminus1,k}.
\end{equation}
Due to the orthogonality of the basis functions, the coefficients $\hat q_{i,k}$ can be written
\begin{equation}
\begin{aligned}\label{eq:qm1_dg}
\hat q_{\shortminus1,k} =\frac{1}{h/2 + \alpha h + h/3}\biggl( &\frac{1}{2}\int_{K_{\shortminus2}} \hat U_{\shortminus2} \hat \varphi_{\shortminus1,k} ~dx + \int_{K_{\shortminus1}} \hat U_{\shortminus1}   \hat \varphi_{\shortminus1,k} ~dx \\
&+ \frac{1}{3} \int_{K_{0}} \hat U_{0} \hat \varphi_{\shortminus1,k} ~dx \biggr) ~\text{ for } k = 0 \hdots p,
\end{aligned}
\end{equation}
This operation can be viewed as a local coarsening of the provisionally updated DG solution from the cells $K_{\shortminus2}$, $K_{\shortminus1}$, $K_{0}$ onto the merging neighborhood $\hat K_{\shortminus1}$.
Due to the projection operation, $\hat Q_{\shortminus1}(x)$ defined in \eqref{eq:qhatm1} is the closest polynomial of degree $p$ to the piecewise defined provisionally updated solution, $\hat U_{\shortminus2}$, $\hat U_{\shortminus1}$, and $\hat U_{0}$ on $\hat K_{\shortminus1}$, measured in the weighted norm $\sqrt{\langle \cdot,\cdot \rangle_{\hat K_{\shortminus1}}}$ that follows from \eqref{eq:ip1}.
Also, note that formula \eqref{eq:qm1_dg} is a generalization of the finite volume formula in \eqref{eq:qm1_fv}, and simplifies to \eqref{eq:qm1_fv} when $p=0$ since the provisional solution updates $\hat U_{\shortminus2}$, $\hat U_{\shortminus1}$, $\hat U_{0}$ are constants for finite volume methods.  In contrast, they are polynomials for the DG method in \eqref{eq:qm1_dg}.

Each merging neighborhood has its own weighted inner product, which depends on how cells are merged and their associated overlap counts.
For example, $\hat \varphi_{1,k}$ is the $k$th orthonormal basis function of $S^p(\hat K_{1})$ with respect to the weighted $L_2$ inner product
$$
\langle f, g \rangle_{\hat K_{1}} = \frac{1}{h/2 + \alpha h + h/3} \left( \frac{1}{2}\int_{K_{2}} f~g~dx +  \int_{ K_{1}} f~g~dx +  \frac{1}{3}\int_{ K_{0}} f~g~dx \right),
$$
where again the cell sizes and integrals are divided by their associated overlap counts.
Additionally, the solution coefficients of the projected solution on $\hat K_{1}$ are written
$$
\hat q_{1,k} =\frac{1}{h/2 + \alpha h + h/3}\left( \frac{1}{2}\int_{K_{2}} \hat U_{2} \hat \varphi_{1,k} ~dx + \int_{K_{1}} \hat U_{1}   \hat \varphi_{1,k} ~dx + \frac{1}{3} \int_{K_{0}} \hat U_{0} \hat \varphi_{1,k} ~dx \right).
$$
The merging neighborhood associated to large cells contain only one cell of the base grid.  
On these neighborhoods, we have $\langle f, g \rangle_{\hat K_i} = \langle f, g \rangle_{K_i}$ since
\begin{align*}
\langle f, g \rangle_{\hat K_i} &= \frac{1}{h_i/N_i} \biggl(\frac{1}{N_i}\int_{K_i}  fg ~dx \biggr)\\
&= \frac{1}{h_i} \int_{K_i}  fg ~dx \\
&= \langle f, g \rangle_{K_i} , \text{ for } i = \shortminus 3,\shortminus 2,0,2,3.
\end{align*}
due to cancellation of the overlap counts.
Thus, the solution coefficients of the numerical solution on the merging neighborhoods associated to large cells in the base grid are
\begin{align*}
\hat q_{i,k} &= \frac{1}{h_i} \left( \int_{\hat K_i} \hat U_i \hat \varphi_{i,k}~dx  \right) \\
&= \hat c_{i,k} 
\end{align*}
since $\hat K_i = K_i$ and the neighborhood inner product is the same as the base grid inner product.  This implies that $\hat \varphi_{i,k} = \varphi_{i,k}$ for $k = 0 \hdots p$ and $i = \shortminus 3,\shortminus 2,0,2,3$.
In other words, preprocessing is only required in the neighborhood of cut cells.

The final step is to project the stabilized neighborhood polynomials from the merging neighborhoods, $\hat Q_{i}$, back onto the base grid.  For a given cell, $K_i$, in the base grid, we can average the stabilized neighborhood polynomials that overlap it and then project that average onto the base grid basis functions $\varphi_{i,k}$ associated to $K_i$.
For example, $K_0$ is overlapped by three merging neighborhoods: $\hat K_{\shortminus1}$, $\hat K_{0}$, $\hat K_{1}$. Therefore, the stabilized solution coefficients on $K_{0}$ after one stage of the SSP-RK method are written 
$$
c^{n+1}_{0,k} = \frac{1}{h_0}\int_{K_0} \frac{1}{3}(\hat Q_{\shortminus1} + \hat Q_{0} + \hat Q_{1})\varphi_{0,k} ~dx, \text{ for } k = 0 \hdots p.
$$
The solution coefficients on cells overlapped by two merging neighborhoods, $K_{\shortminus2}$ and $K_{2}$,  can be written
\begin{equation}
c^{n+1}_{\shortminus2,k} = \frac{1}{h_{\shortminus2}}\int_{K_{\shortminus2}} \frac{1}{2}(\hat Q_{\shortminus2} + \hat Q_{\shortminus1}) \varphi_{\shortminus2,k} ~~dx,  \text{ and } ~c^{n+1}_{2,k} = \frac{1}{h_2}\int_{K_2} \frac{1}{2}(\hat Q_{2} + \hat Q_{1} )\varphi_{2,k} ~~dx,
\end{equation}
for $k = 0 \hdots p$.  Finally, the solution coefficients on cells overlapped by one merging neighborhood is 
$$
c^{n+1}_{i,k} = \frac{1}{h_i}\int_{K_{i}} \hat Q_{i}\varphi_{i,k}  ~dx, \text{ for } k = 0 \hdots p,~ i = \shortminus 3,\shortminus 1,1,3.
$$
Note that the projections involved in the coarsening and refining operation are linear and thus can be implemented as matrix vector products.
Convergence tests illustrating the state redistribution method in one dimension have been included in the supplemental materials document.

We have just described how to apply state redistribution to a single forward Euler stage of an SSP Runge Kutta method for a scalar conservation law.
For non-SSP explicit Runge Kutta methods, state redistribution can be applied to the intermediate and final solutions.
Extending this technique to systems of conservation laws is done by applying state redistribution to each conserved variable of the system separately.
In the next section, we will describe the details of the pre- and postprocessing steps of the state redistribution algorithm for the discontinuous Galerkin method in two dimensions.

\section{Preprocessing for state redistribution in two dimensions} \label{sec:preprocessing}
There are three preprocessing steps for DG methods and they are similar to those for FV methods.  The main difference is that we need to compute polynomial bases.
In the first step, merging neighborhoods are associated to each cell in the cut cells (section \ref{sec:mhoods}).  
This operation relies only on mesh information already available in many cut cell codes: cell connectivity and the vertices of cut cells.
Then, basis functions are generated on the base grid (section \ref{sec:ccbasis}), and on the merging neighborhoods (section \ref{sec:wbasis}).
This operation relies on quadrature rules for the evaluation of inner products on arbitrary polygons (section \ref{sec:quadrature}).
Since we have only considered static embedded boundaries in this work, these preprocessing steps are done only once. %
For moving boundaries, the mesh preprocessing would need to be done every time the boundary is modified.

\subsection{Merging neighborhoods} \label{sec:mhoods}

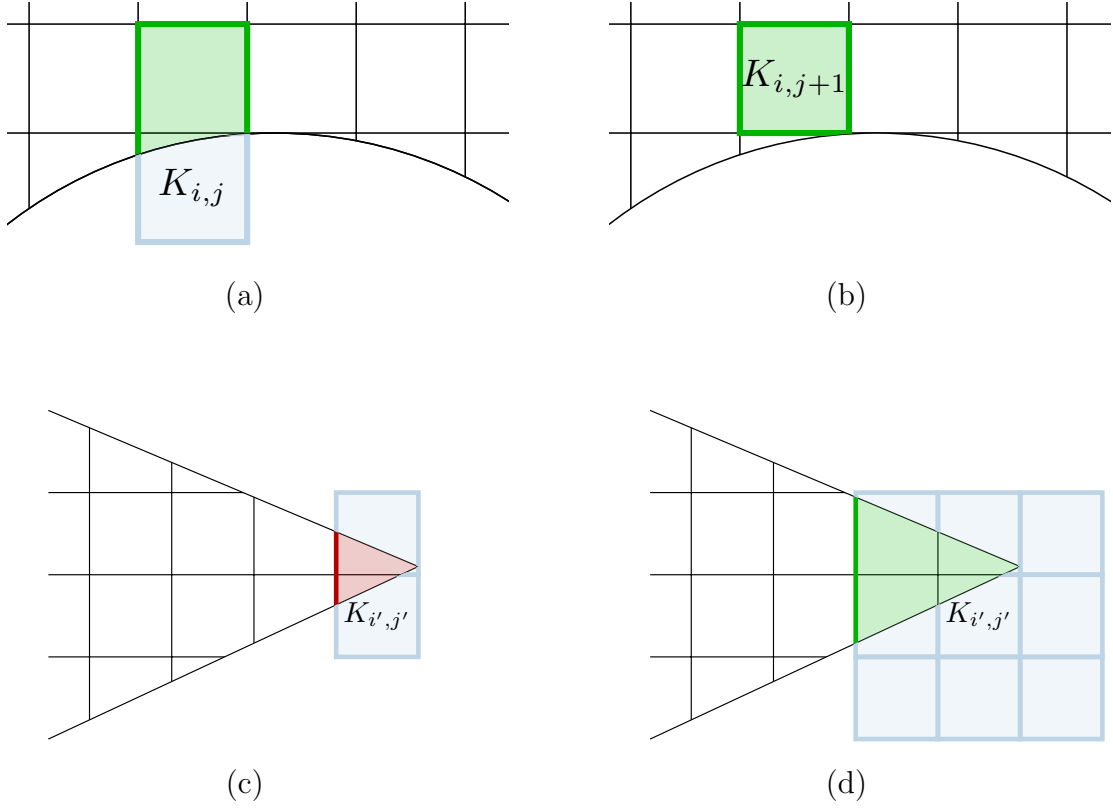
\begin{figure}
	\centering
	\begin{tikzpicture}
	    \node (11) at (0,0) 
	    {\includestandalone[mode=buildnew, width=0.42\linewidth]{2-srd/neigh1}};
	    \node at (0,-2.25) {\textcolor{black}{\large (a)}};
	    \node (12) at (8,0)
	    {\includestandalone[mode=buildnew, width=0.42\linewidth]{2-srd/neigh2}};
	    \node at (8,-2.25) {\textcolor{black}{\large (b)}};

	    \node (11) at (0,-6) 
	    {\includestandalone[mode=buildnew, width=0.42\linewidth]{2-srd/neigh3}};
	    \node at (0,-2.25-6.5) {\textcolor{black}{\large (c)}};
	    \node (12) at (8,-6)
	    {\includestandalone[mode=buildnew, width=0.42\linewidth]{2-srd/neigh4}};
	    \node at (8,-2.25-6.5) {\textcolor{black}{\large (d)}};
	\end{tikzpicture}
	\caption{
	The merging neighborhood of cells $K_{i,j}$ and $K_{i,j+1}$ are highlighted in a) and b), respectively.  These neighborhoods overlap on cell $K_{i,j+1}$. 
	The normal merging neighborhood does not always satisfy the volume constraint in \eqref{eq:vmerge}.
	In c), we show one such neighborhood, associated to $K_{i',j'}$.
	In d), $K_{i',j'}$ is merged instead with cells on the $3 \times 3$ tile centered on $K_{i',j'}$.  This alternative approach to merging results in a neighborhood that is large enough.
	The portions of the grid that are highlighted in blue are solid and do not belong to the computational grid. 
	}\label{fig:nor}

\end{figure}

Merging neighborhoods on two-dimensional cut cells grids are generated following the same procedure as described in \cite{berger2020state}.  
Each cell in the base grid, $K_{i,j}$, is associated to a merging neighborhood $\hat K_{i,j}$.
This merging neighborhood is formed by grouping irregular boundary cells with other nearby cells in the direction closest to the boundary normal.
This is done until the neighborhood volume satisfies
\begin{equation} \label{eq:vmerge}
\sum_{(r,s) \in M_{i,j}} |K_{r,s}|\geq \frac{1}{2}\Delta x \Delta y,
\end{equation}
where $M_{i,j}$ is a set containing the indices of each cell belonging to the merging neighborhood $\hat K_{i,j}$.  For example, in Figures \ref{fig:nor}a) and b) $M_{i,j} = \{(i,j), (i,j+1)\}$ and $M_{i,j+1} = \{ (i,j+1)\}$.
In \eqref{eq:vmerge}, the chosen threshold of $\Delta x \Delta y/2$ is informed by analytical results in \cite{BERGER2015180}.

Note that these neighborhoods can overlap and thus a given element in the base grid can belong to multiple neighborhoods.
Consider the neighborhood $\hat K_{i,j}$ associated to $K_{i,j}$ in Figure \ref{fig:nor}a).  The neighborhood is highlighted in green and is generated when $K_{i,j}$ is merged in the direction closest the boundary normal (upward) with $K_{i,j+1}$.
In contrast, $K_{i,j+1}$ does not need to merge with any neighboring cells, since it is already large enough and satisfies \eqref{eq:vmerge}.  Thus, its merging neighborhood is composed of only itself, $\hat K_{i,j+1} = K_{i,j+1}$, and is highlighted in Figure \ref{fig:nor}b).
Finally, notice that neighborhoods $\hat K_{i,j}$ and $\hat K_{i,j+1}$ overlap on cell $K_{i,j+1}$.

For certain boundary configurations, the normal merging neighborhood will not satisfy \eqref{eq:vmerge}.  This is shown in Figure \ref{fig:nor}c), for the neighborhood associated to $K_{i',j'}$.  
In this case, one could merge $K_{i',j'}$ with cells located on the $3\times 3$ tile centered on cell $K_{i',j'}$ (Figure \ref{fig:nor}d), also called the central merging neighborhood.

\subsection{Cut cell basis generation} \label{sec:ccbasis}
On all cells in the base grid, whole and cut, DG requires a polynomial basis with which to represent the numerical solution.  This basis is precomputed and stored before the time stepping portion of the code.
We prefer to use a basis that is orthogonal with respect to the inner product
\begin{equation} \label{eq:l2_disc}
\langle f,g \rangle_{K_{i,j}} = \frac{1}{|K_{i,j}|} \sum_{q = 0}^{N^{\text{vol}}_{i,j}} w_{i,j,q} f(x_{i,j,q}, y_{i,j,q}) g(x_{i,j,q}, y_{i,j,q}),
\end{equation}
a discrete approximation to \eqref{eq:l2}.  
In the above, $w_{i,j,q},x_{i,j,q},y_{i,j,q}$ are the $q$th volume quadrature weights and points on $K_{i,j}$ (section \ref{sec:quadrature}).
The inner product \eqref{eq:l2_disc} requires the value of the basis functions at the volume and surface quadrature points in order to evaluate the numerical solution in \eqref{eq:sol} for the volume and surface integrals in \eqref{eq:dg}.
The gradient of the basis functions is also needed at the volume quadrature points for the volume integral in \eqref{eq:dg}.
This can be done by taking the monomial basis, and applying the modified Gram-Schmidt (MGS) algorithm \cite{bassi2012flexibility, mascotto2018ill}.  A coordinate transformation can also be applied to the initial non-orthogonal monomial basis to reduce the effects of finite precision computations \cite{bassi2012flexibility}. 
In this work, we describe a similar approach, where we determine these basis function values using the $QR$ factorization of a matrix that we describe below.

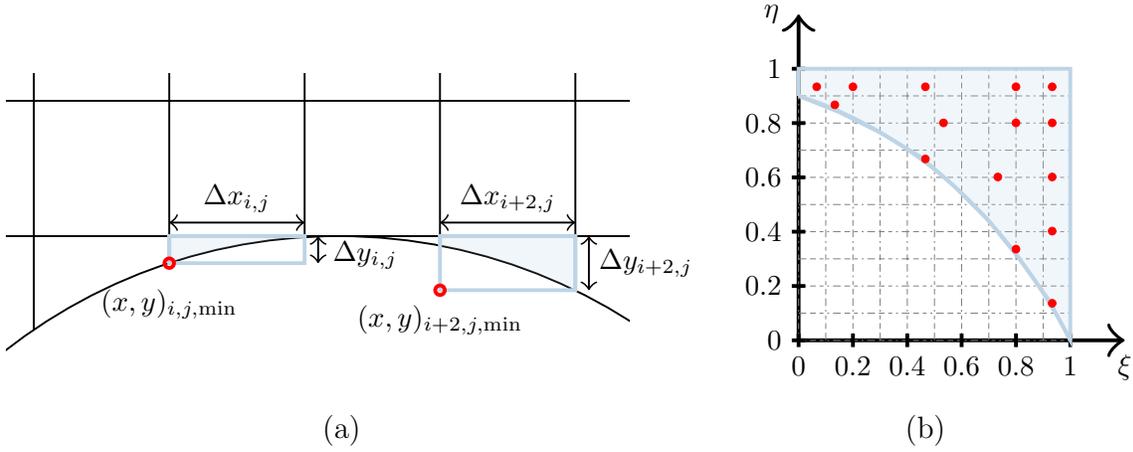
\begin{figure}
	\centering
	\resizebox{\textwidth}{!}{
	\begin{tikzpicture}
	
	    \node (11) at (0,0) 
	    {\includestandalone[mode=buildnew, width=0.6\linewidth]{2-srd/bounding_box}};
	    \node at (0,-2.75) {\textcolor{black}{\large (a)}};
	    \node (12) at (8,0.5)
	    {\includestandalone[mode=buildnew, width=0.35\linewidth]{2-srd/mapping}};
	    \node at (8,-2.75) {\textcolor{black}{\large (b)}};
	\end{tikzpicture}
	}
	\caption{
	The above figures illustrate how the bounding box of a cut cell is used to generate basis functions.
	In a), we indicate the bottom left coordinates ($\textcolor{red}{\bm \circ}$) and dimensions of the bounding box for cells $K_{i,j}$ and $K_{i+2,j}$ used for the basis function calculation on the grid in Figure \ref{fig:nor}a) and b).
	In b), we plot the points $(\xi,\eta)_{i+2,j,q}$ for $q = 0, 2, \hdots$ (indicated by $\textcolor{red}{\bullet}$) at which the basis functions $P_k(\xi, \eta)$ are evaluated to populate the initial Vandermonde-like matrix $V_{i+2,j}$ in \eqref{eq:vander} for element $K_{i+2,j}$.
	}\label{fig:bounding_box}
\end{figure}

The matrix that we factorize is based on a Vandermonde-like matrix, populated with the values of initial basis functions that are not necessarily orthogonal with respect \eqref{eq:l2_disc}.
This initial basis $\{P_k(\xi, \eta)\}_{k = 0,\hdots, N_p}$ is chosen to span $S^p([0,1]^2)$ and is orthonormal with respect to the inner product $\langle \cdot, \cdot \rangle_{[0,1]^2}$.  It is generated symbolically using the Gram-Schmidt procedure, using the monomials ordered as $\{ 1,\xi,\eta,\xi^2,\xi \eta ,\eta^2,\hdots \}$. The first three of these functions are
\begin{equation*}
P_0(\xi, \eta) = 1, \qquad P_1(\xi, \eta) = 2\sqrt{3}(\xi-1/2), \qquad P_2(\xi, \eta) = 2\sqrt{3}(\eta-1/2). 
\end{equation*}
Then, we evaluate $P_k(\xi_{i,j,q}, \eta_{i,j,q})$ at each volume quadrature point $(x_{i,j,q}, y_{i,j,q})$ on cell $K_{i,j}$ using the mapping
\begin{equation}
\label{eq:mapping}
(\xi_{i,j,q}, \eta_{i,j,q}) = \left( \frac{x_{i,j,q} - x_{i,j, \text{min}} }{\Delta x_{i,j}} , \frac{y_{i,j,q} - y_{i,j, \text{min}} }{\Delta y_{i,j}}\right),
\end{equation}
where $(\xi_{i,j,q}, \eta_{i,j,q})$ is the volume quadrature point $(x_{i,j,q}, y_{i,j,q})$ mapped to $[0,1]^2$, $\Delta x_{i,j}$ is the width of the cell $i,j$'s bounding box in $x$ direction, and $x_{i,j,\text{min}}$ is the minimum $x$ component of all the vertices of $K_{i,j}$.  $\Delta y_{i,j}$ and $y_{i,j,\text{min}}$ are similarly defined in the $y$ direction.
These quantities are illustrated in Figure \ref{fig:bounding_box} for the mesh in Figure \ref{fig:nor}a) and b).
Computing the basis function values at the quadrature points uses the $QR$ factorization of a matrix in an oblique inner product \cite{lowery2014stability}, and requires three steps:
\begin{enumerate}
    \item Assemble the initial values of $P_k(\xi_{i,j,q}, \eta_{i,j,q})$ into a Vandermonde-like matrix as follows
\begin{equation} \label{eq:vander}
    V_{i,j} = \begin{bmatrix}
    P_0(\xi_{i,j,0}, \eta_{i,j,0}) & P_1(\xi_{i,j,0}, \eta_{i,j,0}) & \hdots &  P_{N_p}(\xi_{i,j,0}, \eta_{i,j,0})\\ 
    P_0(\xi_{i,j,1}, \eta_{i,j,1}) & P_1(\xi_{i,j,1}, \eta_{i,j,1})  & \hdots &  P_{N_p}(\xi_{i,j,1}, \eta_{i,j,1}) \\ 
     & & \ddots & \\
    P_0(\xi_{i,j,N^{\text{vol}}_{i,j}}, \eta_{i,j,N^{\text{vol}}_{i,j}}) & P_1(\xi_{i,j,N^{\text{vol}}_{i,j}}, \eta_{i,j,N^{\text{vol}}_{i,j}})  & \hdots &  P_{N_p}(\xi_{i,j,N^{\text{vol}}_{i,j}}, \eta_{i,j,N^{\text{vol}}_{i,j}}) 
    \end{bmatrix}. 
\end{equation}
\item Multiply the rows of $V_{i,j}$ by the square root of the quadrature weights and using MATLAB we compute the reduced $QR$ factorization of the product 
$$
W_{i,j}^{1/2}~V_{i,j} = Q_{i,j}R_{i,j},
$$
where $W_{i,j} = \text{diag}(w_{i,j,1}, w_{i,j,2}, \hdots, w_{i,j,N_{i,j}^{\text{vol}}})$.
\item Compute the product $W_{i,j}^{\shortminus 1/2} Q_{i,j} $, which contains the orthogonal basis function values at the quadrature points

\begin{equation*}
\small
    W_{i,j}^{\shortminus 1/2} Q_{i,j} = \begin{bmatrix}
    \varphi_{i,j,0}(x_{i,j,0}, y_{i,j,0}) & \varphi_{i,j,1}(x_{i,j,0}, y_{i,j,0}) & \hdots &  \varphi_{i,j,N_p}(x_{i,j,0}, y_{i,j,0})\\ 
    \varphi_{i,j,0}(x_{i,j,1}, y_{i,j,1}) & \varphi_{i,j,1}(x_{i,j,1}, y_{i,j,1}) & \hdots & \varphi_{i,j,N_p}(x_{i,j,1}, y_{i,j,1}) \\ 
     & & \ddots & \\
    \varphi_{i,j,0}(x_{i,j,N^{\text{vol}}_{i,j}}, y_{i,j,N^{\text{vol}}_{i,j}}) & \varphi_{i,j,1}(x_{i,j,N^{\text{vol}}_{i,j}}, y_{i,j,N^{\text{vol}}_{i,j}}) & \hdots &  \varphi_{i,j,N_p}(x_{i,j,N^{\text{vol}}_{i,j}}, y_{i,j,N^{\text{vol}}_{i,j}}) \\ 
    \end{bmatrix}.
\end{equation*}

\end{enumerate}
The information stored in $R_{i,j}$ can be used to determine the gradient of the orthogonal basis functions.
For example, the $x$ and $y$ components of the gradient of the orthonormal basis functions at the volume quadrature points are given by 
$$\frac{1}{\Delta x_{i,j}}\frac{\partial V_{i,j}}{\partial \xi} R^{\shortminus 1}_{i,j}, \text{ and } \frac{1}{\Delta y_{i,j}}\frac{\partial V_{i,j}}{\partial \eta}  R^{\shortminus 1}_{i,j},$$
respectively.  In the above formulas, the $1/\Delta x_{i,j}$ and $1/\Delta y_{i,j}$ multipliers are due to the mapping \eqref{eq:mapping}.
Similarly, the values of the basis functions at the surface quadrature points are given by
$S_{i,j}R^{\shortminus 1}_{i,j}$, where $S_{i,j}$ is a Vandermonde-like matrix containing the initial basis functions $\{P_k(\xi, \eta) \}_{k = 0, \hdots, N_p}$ evaluated at the Gauss-Legendre quadrature points on each edge of $K_{i,j}$ mapped to $[0,1]^2$ using \eqref{eq:mapping}.

\subsection{Merging neighborhood basis generation} \label{sec:wbasis}
\begin{figure}
	\centering
	\resizebox{\textwidth}{!}{
	\begin{tikzpicture}
	
	    \node (11) at (0,0) 
	    {\includestandalone[mode=buildnew, width=0.6\linewidth]{2-srd/bounding_box_hat}};
	    \node at (0,-2.75) {\textcolor{black}{\large (a)}};
	    \node (12) at (8,0.5)
	    {\includestandalone[mode=buildnew, width=0.35\linewidth]{2-srd/mapping_hat}};
	    \node at (8,-2.75) {\textcolor{black}{\large (b)}};
	\end{tikzpicture}
	}
	\caption{
	The above figures illustrate how the bounding box of a merging neighborhood is used to generate the neighborhood basis functions.
	In a), we indicate the bottom left coordinates ($\textcolor{red}{\bm \circ}$) and dimensions of the bounding box for merging neighborhoods $\hat K_{i,j}$ and $\hat K_{i+2,j}$ used for the merging basis function calculation on the grid on the top row of Figure \ref{fig:nor}.
	In b), we plot the points $(\xi,\eta)_{i+2,j,q}$ for $q = 0, 2, \hdots$ (indicated by $\textcolor{red}{\bullet}$) at which the basis functions $P_k(\xi, \eta)$ are evaluated to populate the initial Vandermonde-like matrix for merging neighborhood $\hat K_{i+2,j}$.
	Note that the quadrature points on each subelement of the merging neighborhoods are the same as the ones used in section \ref{sec:ccbasis}.
	}\label{fig:bounding_box_hat}
\end{figure}
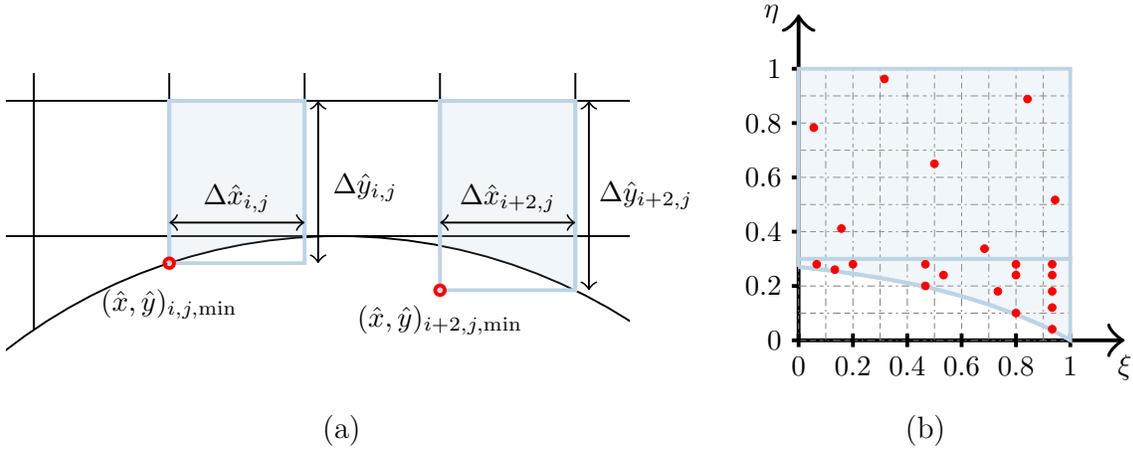

We also need to generate an orthogonal polynomial basis on each merging neighborhood that is orthogonal with respect to the weighted inner product
\begin{equation} \label{eq:wL2}
\langle f,g \rangle_{\hat K_{i,j}} = \frac{1}{|\hat{K}_{i,j}|}\sum_{(r,s) \in M_{i,j}} \frac{1}{N_{r,s}}\int_{K_{r,s}}  f g ~ dx ~dy,
\end{equation}
where $N_{r,s}$ is the number of overlapping neighborhoods on cell $(r,s)$, and 
$$
|\hat{K}_{i,j}| = \sum_{(r,s) \in M_{i,j}} \frac{|K_{r,s}|}{N_{r,s}},
$$
is the weighted volume of the merging neighborhood.
The quadrature points at which the merging neighborhood basis is computed are the same quadrature points as the ones used on the cells in base grid (section \ref{sec:quadrature}).
Similar to the basis functions on the base grid, the merging basis functions are precomputed and stored before the time stepping portion of the code.
These basis function values are computed by adapting \eqref{eq:mapping} to the bounding box of the merging neighborhood.  We evaluate $P_k(\xi_{r,s,q}, \eta_{r,s,q})$ at each of the volume quadrature points $(x_{r,s,q}, y_{r,s,q})$ on the merging neighborhood $\hat K_{i,j}$, $\forall (r,s) \in M_{i,j}$, using the mapping
\begin{equation}
\label{eq:mapping_hat}
(\xi_{r,s,q}, \eta_{r,s,q}) = \left( \frac{x_{r,s,q} - \hat x_{i,j, \text{min}} }{\Delta \hat x_{i,j}} , \frac{y_{r,s,q} -\hat  y_{i,j, \text{min}} }{\Delta \hat y_{i,j}}\right)
\end{equation}
where $(\xi_{r,s,q}, \eta_{r,s,q})$ is again the volume quadrature point $(x_{r,s,q}, y_{r,s,q})$ mapped to $[0,1]^2$, $\Delta \hat x_{i,j}$ is the width of merging neighborhood $i,j$'s bounding box in the $x$ direction, and $\hat x_{i,j,\text{min}}$ is the minimum $x$ component of all the vertices of $\hat K_{i,j}$.  $\Delta \hat y_{i,j}$ and $\hat y_{i,j,\text{min}}$ are similarly defined in the $y$ direction.
These quantities are illustrated in Figure \ref{fig:bounding_box_hat} for two merging neighborhoods on the mesh on the top row of Figure \ref{fig:nor}.
We can quickly precompute the merging basis function values using a discrete inner product and a QR factorization, as described in section \ref{sec:ccbasis}.  
The discrete weighted inner product is given by
\begin{equation} \label{eq:wL2_disc}
\langle f,g \rangle_{ \hat K_{i,j} } = \frac{1}{|\hat{K}_{i,j}|}\sum_{(r,s) \in M_{i,j}} \frac{|K_{r,s}|}{N_{r,s}}   \sum_{q = 0}^{ N^{\text{vol}}_{i,j}} w_{r,s,q} f(r_{r,s,q},s_{r,s,q}) g(r_{r,s,q},s_{r,s,q}) ,
\end{equation}
where $(w_{r,s,q},x_{r,s,q},y_{r,s,q})$ are the quadrature weights and points on each cell in the merging neighborhood $(r,s) \in M_{i,j}$, and again $\sum_q w_{r,s,q} = |K_{r,s}|$ and $w_{r,s,q} > 0$ for $q=0\hdots  N^{\text{vol}}_{i,j}$.
This basis will be used in the projection operations that stabilize the numerical solution, described in the next section.

\section{State redistribution in two dimensions} \label{sec:srd2d}
We are now in a position to describe the state redistribution method on two-dimensional cut cell grids. 
State redistribution is applied after each forward Euler step in a high order SSP-RK time stepper.  For non SSP-RK time steppers, the SRD operator can be applied to the intermediate and final solutions.
The algorithm follows the operations described in one dimension (section \ref{sec:srd1d}), where the numerical solution is provisionally updated, then stabilized by two projection operations.
The first projection moves the provisionally updated numerical solution from the base grid onto the merging neighborhoods (section \ref{sec:srd1d}, step 1).
The second projection moves the merging neighborhood solution back onto the base grid in a conservative manner (section \ref{sec:srd1d}, step 2).

To begin, we compute a forward Euler step or an intermediate solution of the RK scheme using the same $\Delta t$ on all cells, both whole and cut.  For example, a forward Euler step for a scalar conservation law is written
\begin{equation} \label{eq:dg2d}
    \hat{c}_{i,j,k} = c_{i,j,k}^n - \frac{\Delta t}{|K_{i,j}|} \left[ \int_{\partial K_{i,j}} \varphi_{i,j,k}\mathbf{F}^*(U^-_{i,j}, U_{i,j}^+)\cdot \mathbf{n} ~dl  - \int_{K_{i,j}}\mathbf{F}(U_{i,j})\cdot \nabla \varphi_{i,j,k} ~dx~dy\right],
\end{equation}
    where $\hat c_{i,j,k}$ for $k = 0, \hdots N_p$, are the provisionally updated, possibly unstable degrees of freedom, $\mathbf{F}^*(U^-_{i,j}, U_{i,j}^+)$ is the numerical flux evaluated on the cell interface, and $\Delta t$ is the stable time step given by the CFL condition in \eqref{eq:cfl}.
    
Next, we project the provisionally updated numerical solution onto the weighted merging basis precomputed in section \ref{sec:wbasis}, using the discrete weighted inner product \eqref{eq:wL2_disc}.  The solution on the merging basis is
$$
\hat Q_{i,j}(x) = \sum_{k = 0}^{N_p} \hat q_{i,j,k} \hat \varphi_{i,k}(x,y)
$$
where the solution coefficients are
\begin{equation} \label{eq:qcoeffs}
\hat q_{i,j,k} = \frac{1}{|\hat{K}_{i,j}|}  \sum_{(r,s) \in M_{i,j}} \frac{1}{N_{r,s}}\int_{K_{r,s}}\hat{\varphi}_{i,j,k} \hat{U}_{r,s} ~dx ~dy,
\end{equation}
due to the orthogonality of $\hat \varphi_{i,j,k}$, and 
$$
\hat{U}_{r,s}(x) = \sum_{k = 0}^{N_p} \hat{c}_{i,j,k} \varphi_{i,j,k}(x,y).
$$

Finally, the numerical solution on the neighborhoods is projected back onto the base grid by averaging the overlapping neighborhood solutions:
\begin{equation} \label{eq:step3}
c_{i,j, k}^{n+1} = \frac{1}{|K_{i,j}|}  \sum_{(r,s) \in W_{i,j}} \frac{1}{N_{i,j}}\int_{K_{i,j}}\varphi_{i,j,k} \hat{Q}_{r,s} ~dx ~dy,
\end{equation}
where $c_{i,j, k}^{n+1}$ is the stabilized solution coefficient on the base grid, and $W_{i,j}$ is the set of neighborhood indices that overlap cell $(i,j)$.  
In the above formula, we average the projection from different, overlapping neighborhoods that contain the same cell in the base grid.
As an example, the set of neighborhood indices that overlap $K_{i,j+1}$ in Figure \ref{fig:nor}a) and b) is $W_{i,j+1} = \{ (i,j), (i,j+1) \}$.  This set does not in fact need to be precomputed and the action of \eqref{eq:step3} can be implemented using a nested for loop \cite{berger2020state}.  More information on implementation details are provided in the supplemental materials document.  We also note that the coarsening and refining operations in the state redistribution method are linear operations, so they can be implemented with simple matrix vector products.  
To conclude this section, we state two claims
\begin{itemize}
    \item[] \textbf{Claim 1}: state redistribution is \textit{p}-exact.  This is a result of the projections used in the coarsening and refining operations.
    \item[] \textbf{Claim 2}: state redistribution is conservative.  This follows from the analysis in \cite{berger2020state}.
\end{itemize}
Proofs of the above statements are provided in the supplementary material document that accompanies this paper.

\section{Numerical experiments} \label{sec:numex}
In this section, we demonstrate the performance of state redistribution method on a number of test problems.  We solve steady state and time-dependent problems governed by the scalar advection equation and the Euler equations of gas dynamics.
In all numerical experiments, we use a time step that is proportional to the size of cells in the background mesh, given by
$$
\Delta t \biggl(\frac{|a|}{\Delta x} + \frac{|b|}{\Delta y} \biggr)= \frac{0.9}{2p+1},
$$
where $|a|$, $|b|$ are the maximum wave speeds in the $x$ and $y$ directions \cite{COCKBURN1998199, qinkrivo}.  These numerical examples were completed using SSP-RK schemes for $p=1$, and $p=2$, and non-SSP schemes for $p \geq 3$.
On all embedded boundary grids, $\Delta x = \Delta y$ and $N$ is the number of cells in each direction on the background Cartesian grid.
We compute the discrete $L_1$ and $L_\infty$ errors
\begin{equation}
\begin{aligned}
L_1 = \sum_{i,j}\sum^{N^{\text{vol}}_{i,j}}_{q=0} ~ w_{i,j,q}|U_{i,j}(x_{i,j,q},y_{i,j,q}) - u(x_{i,j,q},y_{i,j,q})|, \\
L_\infty = \max_{i,j} \max_{q}|U_{i,j}(x_{i,j,q},y_{i,j,q}) - u(x_{i,j,q},y_{i,j,q})|,
\end{aligned} \label{eq:l1_linf_standard}
\end{equation}
where $U_{i,j}(x,y)$ and $u(x,y)$ are respectively the numerical and exact solutions.  Note that the above formulae only approximate the true $L_1$ and $L_\infty$ errors.

\subsection{Solid body rotation} \label{sec:solid}
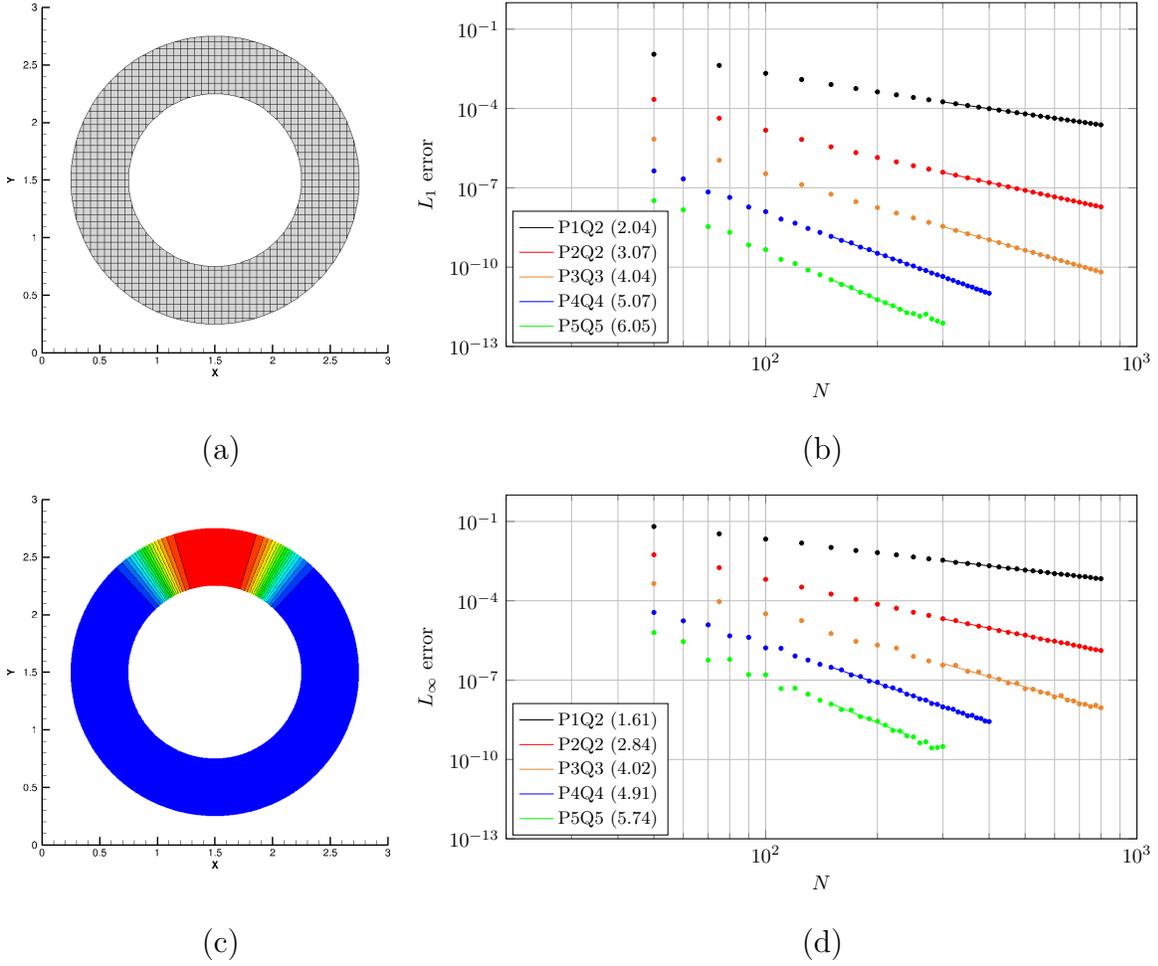
\begin{figure}

\centering
\begin{tikzpicture}
\node (11) at (0, 0.175){\includegraphics[width = 0.32\textwidth]{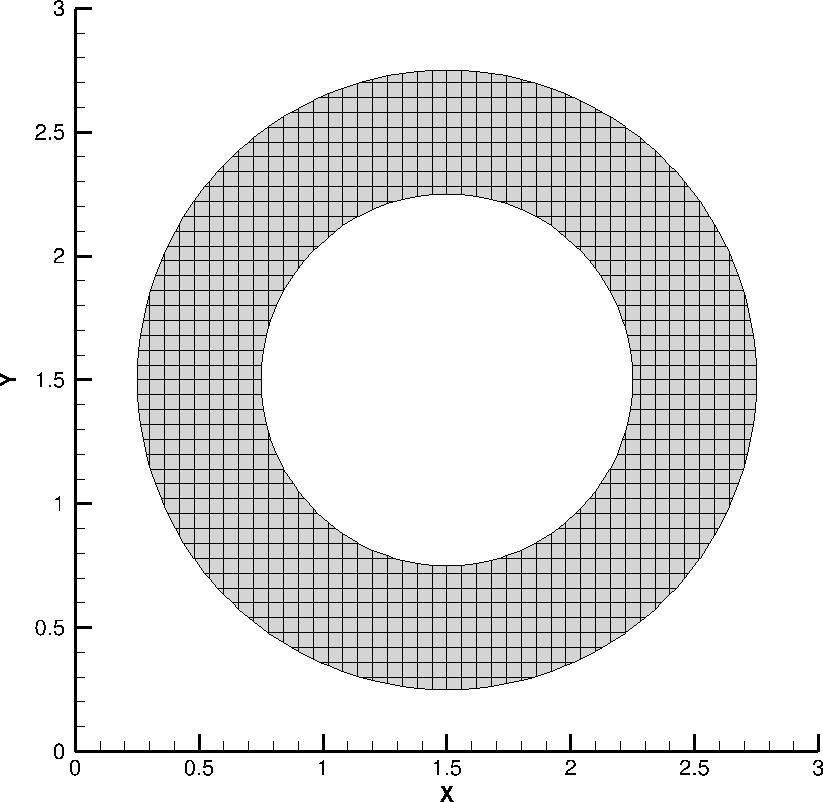}};
\node at (0.3,-3.3) {\textcolor{black}{\large (a)}};
\node (12) at (7.8,0.){\includestandalone[width = 0.625\textwidth]{3-rh/l1}};
\node at (8.3,-3.3) {\textcolor{black}{\large (b)}};

\node (21) at (0, -6.375){\includegraphics[width = 0.32\textwidth]{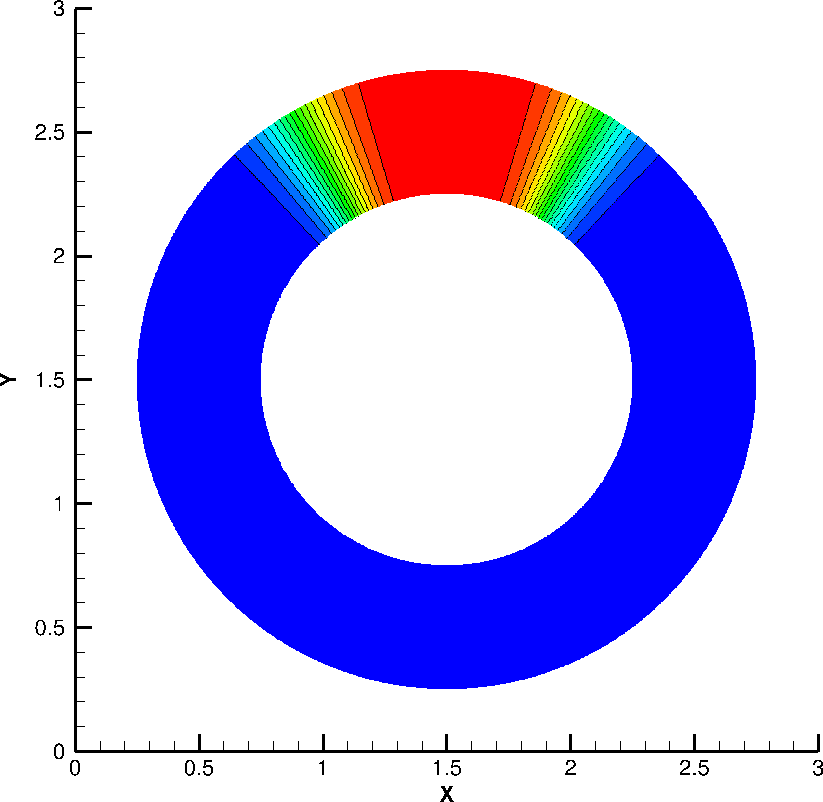}};
\node at (0.3,-9.85) {\textcolor{black}{\large (c)}};
\node (22) at (7.8,-6.55){\includestandalone[width = 0.625\textwidth]{3-rh/linf}};
\node at (8.3,-9.85) {\textcolor{black}{\large (d)}};
    \end{tikzpicture}
    \caption{In a) we plot the annulus domain superimposed on the $50\times50$ background grid. In c) we plot the isolines of exact solution at the initial and final time. In b) and d) we plot $L_1$ and $L_\infty$ errors of the solid body rotation problem, respectively.  
    In the legend, we indicate the degree of polynomial approximation on the cells and the boundary segments by $P$ and $Q$, respectively.
    The rate of convergence quoted in parentheses is computed by a least squares fit, plotted with a straight line.
    We observe the expected $p+1$ rate of convergence in the $L_1$ norm, and a rate of convergence between $p$ and $p+1$ in the $L_\infty$ norm.
    Additionally, the error in the $L_\infty$ norm is noisy, due to the irregularity of the grid at the cut cells.
    In figures b) and d), $N$ refers to the number of cells in the $x$ and $y$ directions of the background grid.
    } \label{fig:rotating}
\end{figure}

In this example, we solve for the solid body rotation of a pulse around an annulus described in \cite{Berger2012ASH,doi:10.1137/S106482750343028X}.  The conservation law
$$
u_t - 0.4\pi [\left( y - 1.5 \right)u]_x + 0.4\pi[ \left( x - 1.5 \right) u]_y = 0,
$$
rotates the initial condition
$$
u(x,y,0) = w(\theta - \pi/2),
$$
where
$$
w(\theta) = \frac{1}{2}\biggl[ \text{erf}\left( 5(\pi/6 - \theta) \right) + \text{erf}\left( 5(\pi/6 + \theta) \right)\biggr].
$$
The initial condition is plotted in Figure \ref{fig:rotating}c) on the domain enclosed by two concentric discs of radii $R_1 = 0.75$ and $R_2 = 1.25$.
We use the upwind numerical flux and impose zero flow at the boundary, i.e., $\mathbf{F}^{*}\cdot \mathbf{n} = 0$ at each boundary quadrature point. 
Figure \ref{fig:rotating}a) shows the computational domain embedded on a $50 \times 50$ grid.  
The exact solution is the initial condition rotated about the point $(1.5,1.5)$, where the final time $T=5$ corresponds to one solid body rotation.  The annulus is embedded on the domain $[0,3.0001]^2$, where the $x$ and $y$ domain length is slightly perturbed from 3 to prevent cell degeneracies.
We aim to avoid cases where the background grid intersects the annulus with radii $R_1 = 0.75$ and $R_2 = 1.25$ in a degenerate fashion, which would create edges of zero length.
The minimum volume fraction, $\min_{i,j} |K_{i,j}|/(\Delta x \Delta y)$, in this convergence test is 6.84e-10.

In Figure \ref{tab:error_comparison}d) and e), we plot the DG-$P5Q5$ numerical solution on the inner and outer embedded boundaries on a coarse $10 \times 10$ grid at the final time.  Segments corresponding to different elements are plotted in different colors. 
We see that there are no glitches or instabilities in the numerical solution on the cut cells.

\subsubsection*{Discussion of the accuracy}
In Figures \ref{fig:rotating}b) and d), we provide the $L_1$ and $L_\infty$ errors in \eqref{eq:l1_linf_standard}, respectively.  The expected $p+1$ rate of convergence in the $L_1$ norm is observed, but the rate of convergence in the $L_\infty$ norm is between $p$ and $p+1$.
This is in contrast to the one-dimensional results in the supplementary material document, where the full rate of convergence in both norms was observed.
This is the price of automatic mesh generation that we are willing to pay.
The nonsmoothness of the error in the $L_\infty$ norm has been observed before and is due to the irregularity of the grid at the cut cells \cite{berger2020state,engwer2020stabilized}.
In \cite{engwer2020stabilized}, it is found that the rate of convergence in $L_\infty$ depends on the geometry of the domain.
Here, we find that this rate also depends on the polynomial degree of approximation and varies between $p$ and $p+1$.

\subsubsection*{Comparison with finite volume schemes}
We compare the numerical solution obtained with our DG scheme to the one obtained with a second order rotated grid $h$-box method for finite volumes described in \cite{doi:10.1137/S106482750343028X}.

In \cite{doi:10.1137/S106482750343028X}, the $L_1$ error at the final time for second order accurate finite volume methods is measured on the domain interior and boundary using
\begin{equation}\label{eq:fv_error}
E_d = \frac{\sum_{i,j}|\overline{U}_{i,j} - u(x_{i,j},y_{i,j})|~|K_{i,j}|}{\sum_{i,j}|u(x_{i,j},y_{i,j})|~|K_{i,j}|}~ \text{ and } ~E_b = \frac{\sum_{ (i,j) \in B } |\overline{U}_{i,j} - u(x_{i,j},y_{i,j})| ~|b_{i,j}| }{\sum_{ (i,j) \in B } |u(x_{i,j},y_{i,j})| ~|b_{i,j}| },
\end{equation}
respectively, where $(x_{i,j},y_{i,j})$ are the cell centroids, $B$ is the set of indices of all cells that lie on the domain boundary, and $|b_{i,j}|$ is length of the cut cell's boundary segment.
Since second order finite volume methods typically represent curved boundaries with linear segments, we compute the $L_1$ errors in \eqref{eq:fv_error} for both DG-$P1Q1$ and DG-$P1Q2$ methods.

We find that the error is smaller on the outer boundary than on the inner boundary in Figures \ref{tab:error_comparison}b) and c).
This agrees with the results in \cite{doi:10.1137/S106482750343028X}, where this phenomenon is attributed to the outer boundary having more cells to resolve the sloped regions of the pulse.
As expected, the DG-$P1Q2$ scheme is more accurate than DG-$P1Q1$ scheme due to the higher fidelity representation of the boundary.
Additionally, the DG-$P1Q1$ scheme is approximately an order of magnitude more accurate than the FV $h$-box scheme in Figure \ref{tab:error_comparison}a).
The boundary errors of the DG schemes are smaller than the FV ones, although the difference between FV and DG schemes here is not as pronounced.
Finally, we observe that the FV scheme is second order accurate at the boundary, while the DG-$P1Q1$ and DG-$P1Q2$ schemes are not.

\begin{figure}
\centering

\resizebox{\columnwidth}{!}{
\begin{tikzpicture}

\node (12) at (9.5,-1.8){\includestandalone[width = 0.5\textwidth]{3-rh/boundary1}};
\node (13) at (9.9,-3) {\textcolor{black}{\large (d) Outer bdry. }};

\node (14) at (9.5,-7.9){\includestandalone[width = 0.5\textwidth]{3-rh/boundary2}};
\node (15) at (9.9,-9.1) {\textcolor{black}{\large (e) Inner bdry. }};

\node (1) at (0,0) {
\begin{tabular}{|c||c|c|c|}                             
\hline
$N$ & FV (\textit{h}-box \cite{doi:10.1137/S106482750343028X}) & DG-$P1Q1$ (SRD)         & DG-$P1Q2$  (SRD)        \\ \hline \hline
100 & 7.50e-02  (\hspace{7pt}-\hspace{7pt}) & 4.43e-03  (\hspace{7pt}-\hspace{7pt})& 4.24e-03  (\hspace{7pt}-\hspace{7pt}) \\ \hline
200 & 2.03e-02 (1.88) & 9.17e-04 (2.27)& 8.36e-04 (2.34) \\ \hline
400 & 5.14e-03 (1.98) & 2.16e-04 (2.08)& 1.90e-04 (2.13) \\ \hline
800 & 1.29e-03 (1.99) & 5.35e-05 (2.01)& 4.60e-05 (2.04) \\ \hline
\end{tabular}
};
\node (11) at (0,-1.5) {\textcolor{black}{\large (a) Domain $E_d$ }};

\node (2) at (0,-4.35) {
\begin{tabular}{|c||c|c|c|}
\hline
$N$ & FV (\textit{h}-box \cite{doi:10.1137/S106482750343028X}) & DG-$P1Q1$       (SRD)  & DG-$P1Q2$    (SRD)       \\ \hline \hline
100 & 3.87e-02  (\hspace{7pt}-\hspace{7pt}) &  6.73e-03  (\hspace{7pt}-\hspace{7pt}) & 3.29e-03  (\hspace{7pt}-\hspace{7pt})\\ \hline
200 & 1.01e-02 (1.93) &  2.54e-03 (1.40) & 1.14e-03 (1.52)\\ \hline
400 & 2.56e-03 (1.98) &  1.01e-03 (1.33) & 4.57e-04 (1.31)\\ \hline
800 & 6.45e-04 (1.99) &  4.29e-04 (1.23) & 1.77e-04 (1.36)\\ \hline
\end{tabular}
};
\node (22) at (0,-5.85) {\textcolor{black}{\large (b) Outer boundary $E_b$ }};

\node (3) at (0,-8.55) {
\begin{tabular}{|c||c|c|c|}
\hline
$N$                       & FV (\textit{h}-box \cite{doi:10.1137/S106482750343028X})& DG-$P1Q1$ (SRD)       & DG-$P1Q2$     (SRD)     \\ \hline \hline
100 & 1.39e-01  (\hspace{7pt}-\hspace{7pt}) & 1.59e-02  (\hspace{7pt}-\hspace{7pt}) & 1.49e-02  (\hspace{7pt}-\hspace{7pt})\\ \hline
200 & 4.01e-02 (1.79) & 5.37e-03 (1.56) & 4.39e-03 (1.76)\\ \hline
400 & 1.02e-02 (1.97) & 2.10e-03 (1.35) & 1.40e-03 (1.64)\\ \hline
800 & 2.56e-03 (1.99) & 8.78e-04 (1.25) & 4.69e-04 (1.58)\\ \hline
\end{tabular}
};
\node (11) at (0,-10.15) {\textcolor{black}{\large (c) Inner boundary $E_b$ }};

\end{tikzpicture}
}
\caption{On the left, we compare solutions to the solid body rotation problem in section \ref{sec:solid} computed with DG-$P1Q1$, DG-$P1Q2$, and FV \textit{h}-box \cite{doi:10.1137/S106482750343028X} methods.  
In a), b), and c) the $L_1$ errors on the domain, and on the inner and outer boundaries ($E_d$, $E_b$) of the numerical solution computed by all three numerical methods are provided, respectively.
The errors in the FV column are taken from Table 3.1 of \cite{doi:10.1137/S106482750343028X}, and the numbers in parentheses are the estimated order of convergence.
On the right in d) and e), we plot the DG-$P5Q5$ numerical solution on the outer and inner embedded boundary of a $10 \times 10$ grid, respectively.
Segments corresponding to different elements are plotted in different colors.
}
\label{tab:error_comparison}
\end{figure}

\subsection{Ringleb flow}
\begin{figure}
\centering
\begin{tikzpicture}[spy using outlines={circle, magnification=3, connect spies}]
\node (11) at (0, 0.175){\includegraphics[width = 0.32\textwidth]{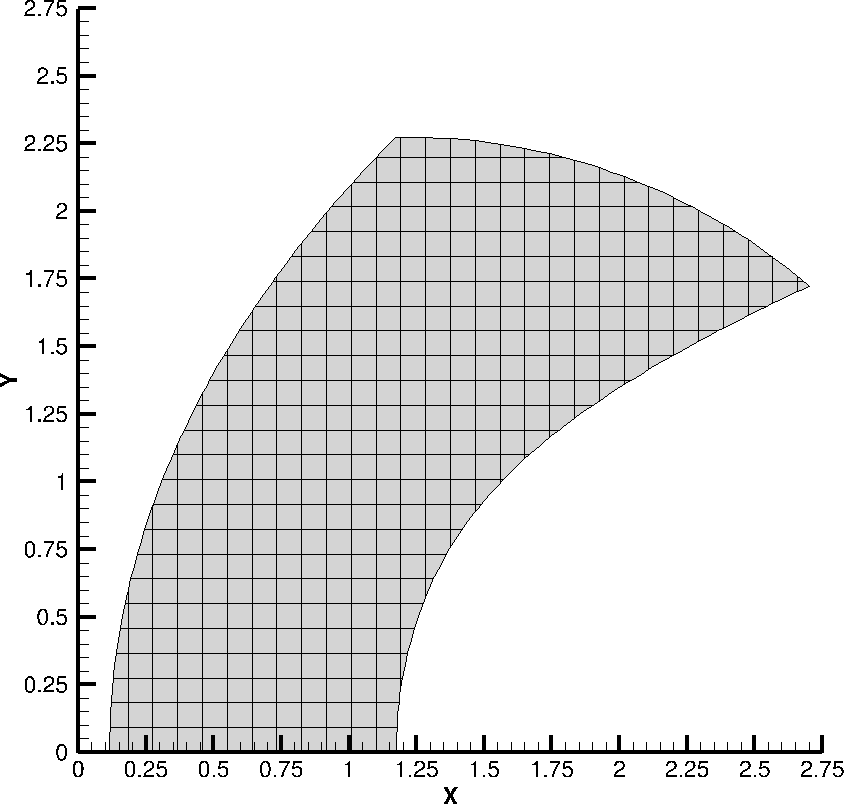}};
\node at (0.3,-3.3) {\textcolor{black}{\large (a)}};
\node (12) at (7.8,0.){\includestandalone[width = 0.625\textwidth]{6-ringleb/l1}};
\node at (8.3,-3.3) {\textcolor{black}{\large (b)}};

\node (21) at (0, -6.375){\includegraphics[width = 0.32\textwidth]{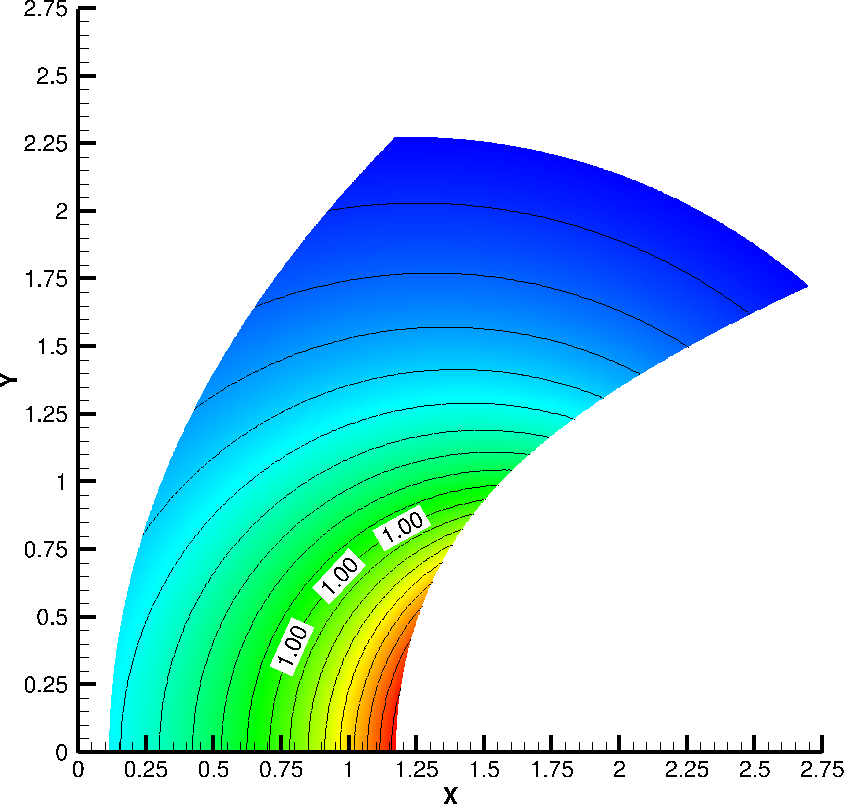}};
\node at (0.3,-9.85) {\textcolor{black}{\large (c)}};
\node (22) at (7.8,-6.55){\includestandalone[width = 0.625\textwidth]{6-ringleb/linf}};
\node at (8.3,-9.85) {\textcolor{black}{\large (d)}};
    \end{tikzpicture}
    \caption{In a) we plot the Ringleb domain superimposed on a $30\times30$ background grid. In c) we plot the Mach isolines of the exact solution, where the sonic line is labeled.  In b) and d) we plot $L_1$ and $L_\infty$ entropy errors, respectively.  
    In the legend, we indicate the degree of polynomial approximation on the cells and the boundary segments by $P$ and $Q$, respectively.
    The rate of convergence quoted in parentheses is computed by a least squares fit, plotted with a straight line.
    In figures b) and d), $N$ refers to the number of cells in the $x$ and $y$ directions of the background grid.
    } \label{fig:ringleb}
\end{figure}
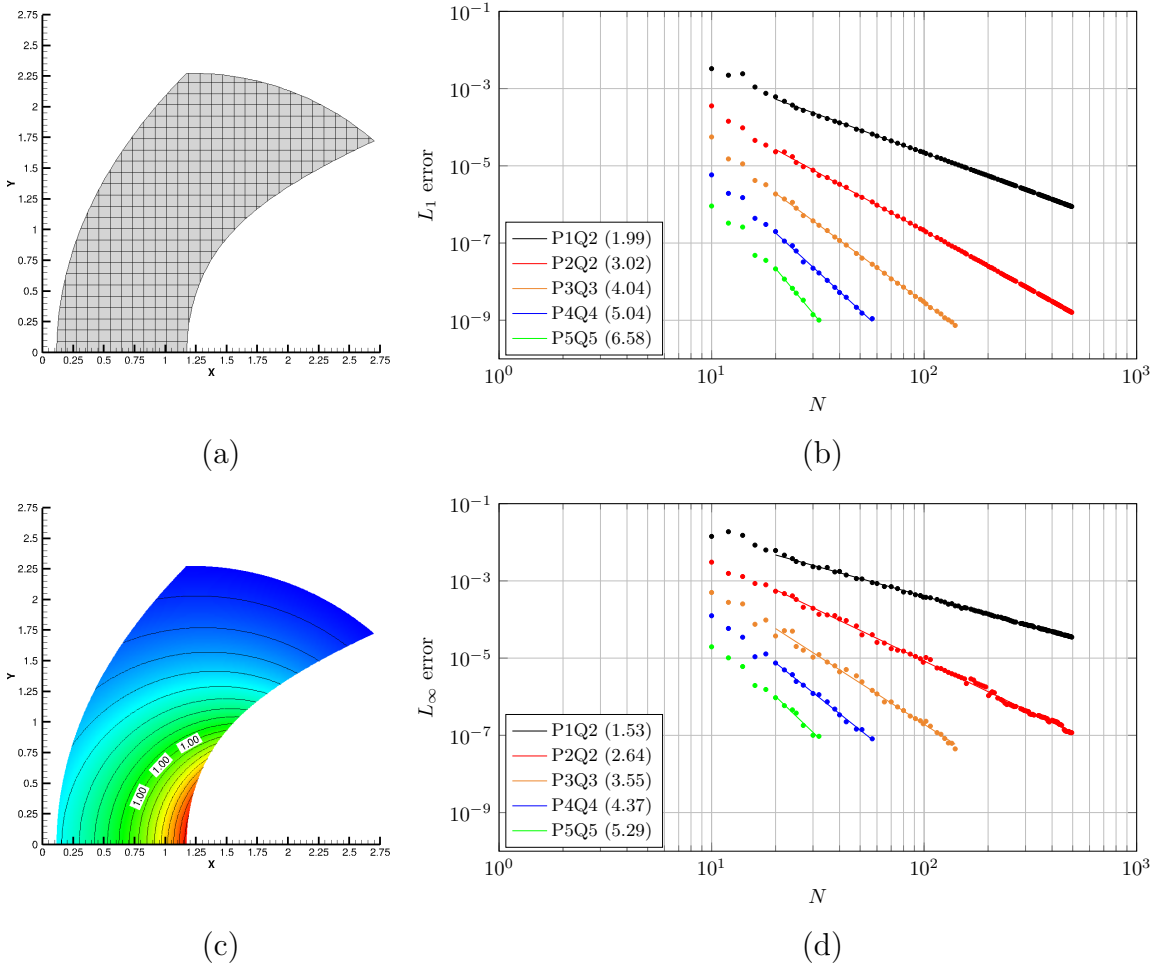

In this example, we solve the Euler equations
\begin{equation} \label{eq:euler}
\frac{\partial}{\partial t}\begin{pmatrix}
\rho \\
\rho u \\
\rho v \\
E
\end{pmatrix}  + \frac{\partial}{\partial x}\begin{pmatrix}
\rho u\\
 \rho u^2 + P\\
 \rho uv \\
 u(E+P)
\end{pmatrix}  + \frac{\partial}{\partial y}\begin{pmatrix}
\rho v\\
 \rho u v \\
 \rho v^2 + P\\
 v(E+P) 
\end{pmatrix}  = 0  \text{ with } P = (\gamma-1)[E - \frac{1}{2}\rho(u^2 + v^2)]
\end{equation}
for smooth transonic flow through a curved channel by time stepping to steady state. 
This test problem is useful to evaluate the accuracy of high order methods with curved boundaries \cite{WANG2006154,DUMBSER2007204,houston2018automatic,wang2013high,IVAN2014830}.
The exact solution at a given point $(x,y)$ is computed by solving
\begin{equation*}
\left( x - \frac{J}{2} -s_x \right)^2 + (y-s_y)^2 = \frac{1}{4 \rho^2 k^4},
\end{equation*}
for the speed of sound $c$ using the method of bisection, where
\begin{equation*}
\begin{aligned}
    J &= \frac{1}{c} + \frac{1}{3c^3} + \frac{1}{5c^5} - \frac{1}{2}\ln \biggl(\frac{1+c}{1-c} \biggr), \\
    q &= \sqrt{\frac{2(1-c^2)}{\gamma-1}}, ~ k = \sqrt{ \frac{2}{1/ q^2 - 2 \rho (x - J / 2 - s_x) } },\\
    \rho &= c^{2/(\gamma-1)},
\end{aligned}
\end{equation*}
$\gamma = 1.4$, and $\rho$, $q$ are the density and speed of the gas, respectively.  
Additionally, $(s_x,s_y) = (1.5,0)$ are shifts to center the channel on the domain.
The pressure, and $x,y$ components of the flow velocity are given by
$$
P = \frac{c^2\rho}{\gamma}, ~u = \sqrt{q^2-v^2}, ~\text{ and } ~v= \frac{q^2}{k},
$$
respectively.
The boundaries of the channel are embedded on $[0,2.75]^2$, and are given by the curves
\begin{align*}
x(k,q) &= \frac{1}{2\rho}\biggl(\frac{1}{q^2}-\frac{2}{k^2} \biggr) + \frac{J}{2} + s_x, \\
y(k,q) &= \frac{1}{kq\rho}\sqrt{1-\frac{q^2}{k^2} } + s_y.
\end{align*}
The reflecting walls are described by $k=0.7$ and $k = 1.2$, while the top boundary is described by $q = 0.5$.
The domain and exact solution are shown in Figure \ref{fig:ringleb}a) and c).
The exact solution is used as the initial condition, and as the ghost states at the top and bottom boundaries.  Roe's numerical flux is used on the interior element interfaces, and the flux on the reflecting walls is given by 
\begin{equation} \label{eq:refl}
\mathbf{F}^*\cdot \mathbf{n} = \begin{pmatrix}
0 \\
P n_x\\
P n_y \\
0
\end{pmatrix},
\end{equation}
where $\mathbf{n} = (n_x, n_y)$ is the outward facing normal and $P$ is the interior pressure at the boundary at each quadrature point.
The solver is run until the stopping criterion 
\begin{equation}\label{eq:stopping}
\max_{i,j} || \mathbf{c}_{i,j}^{n+1} - \mathbf{c}_{i,j}^{n} ||_{\infty} < 10^{\shortminus 13}
\end{equation}
is satisfied.
The minimum volume fraction, $\min_{i,j} |K_{i,j}|/(\Delta x \Delta y)$, in this convergence test is 1.25e-09.

\begin{figure}[htpb]
\centering

\begin{tikzpicture}[spy using outlines={circle, magnification=2.5, connect spies}]
\def\sy{-0.35}
\node (31) at (0, 0.175){\includegraphics[width = 0.5\textwidth]{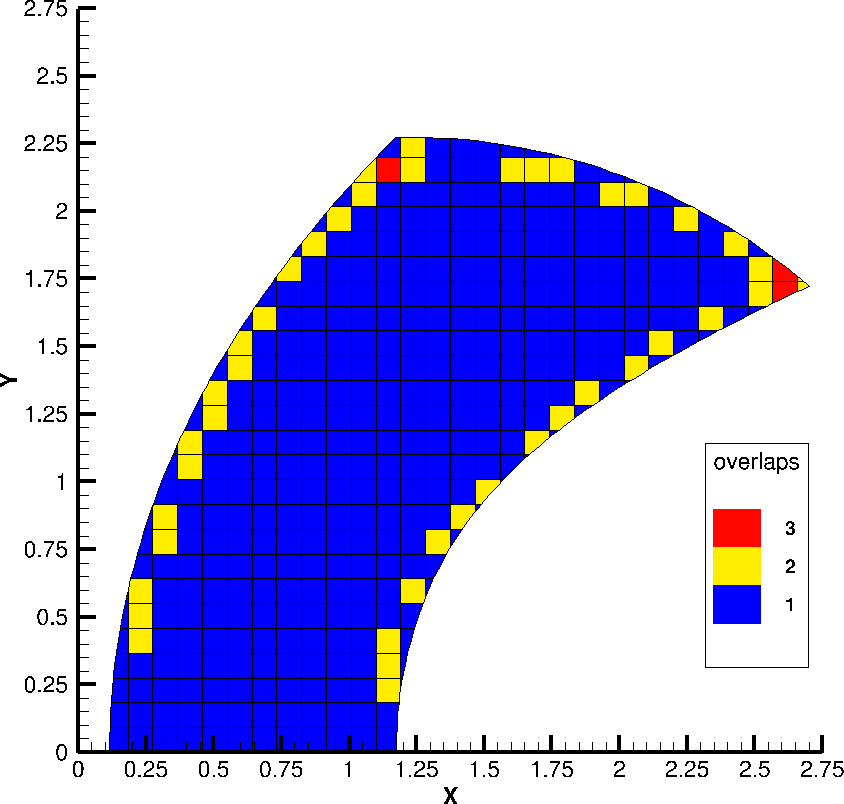}};
\node at (0.3,-4.2) {\textcolor{black}{\large (a)}};

\node (32) at (7.35,3.15+\sy){\includestandalone[width = 0.37\textwidth]{6-ringleb/overlaps1_ringleb}};
\node at (7.9,0.3+\sy) {\textcolor{black}{\large (b)}};
\node (33) at (7.35,-1.25+\sy){\includestandalone[width = 0.37\textwidth]{6-ringleb/overlaps2_ringleb}};
\node at (7.9,-4.2+\sy) {\textcolor{black}{\large (c)}};

\spy [black, every spy on node/.append style={thick}, size=2.cm] on (-0.3,2.45) in node[fill=white] at (-2,2.75);
\end{tikzpicture}

\caption{In a), we plot the overlap counts of the Ringleb domain on the $30 \times 30$ grid.  On most cut cells, we use the normal merging neighborhood.  When the volume constraint \eqref{eq:vmerge} cannot be satisfied, or there is a corner on the cut cell, we use the $3 \times 3$ central merging neighborhood.  This is illustrated with the neighborhood associated to $ K_{i,j}$ in b), which overlaps with the neighborhood associated to $K_{i-1,j-1}$ shown in c).
In both b) and c), the cells highlighted in green are the merging neighborhoods and the portions that are highlighted in blue are solid and do not belong to the computational grid.  
The highlighted red cell in a), $K_{i,j-1}$, is overlapped by three neighborhoods: $\hat K_{i,j}$, $\hat K_{i-1,j-1}$, and $\hat K_{i,j-1}$.
}\label{fig:overlaps_ringleb}
\end{figure}
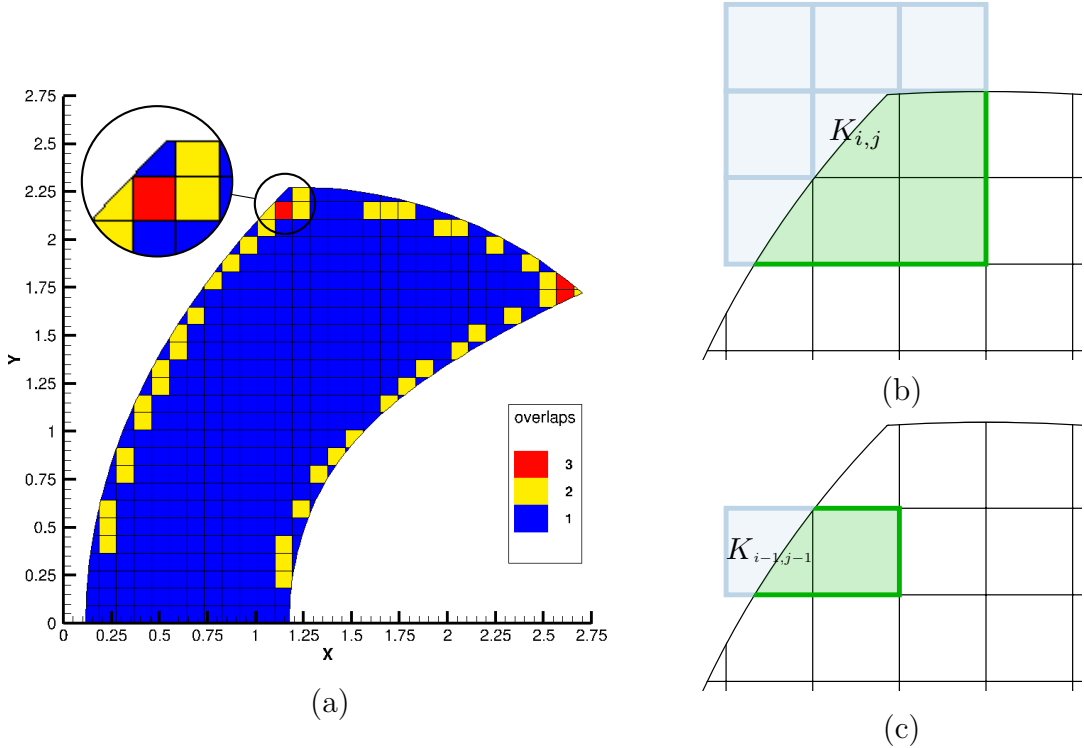

This test case is interesting for two reasons.  First, the flow smoothly transitions from subsonic to supersonic and the Mach number reaches an approximate maximum of $1.42$ near the bottom right corner of the domain. Second, the embedded boundary presents sharp corners.
The cut cells on which these corners lie may not have a unique normal merging neighborhood.  For example, in Figure \ref{fig:overlaps_ringleb}a), the cut cell on the top left corner of the domain can either be merged downward or to the right.
Instead of choosing one of the two possible normal merging neighborhoods, we use the $3 \times 3$ central merging neighborhood, shown in Figure \ref{fig:overlaps_ringleb}b). 
We do this for the cut cells at the top left and right corners of the domain, which result in higher overlap counts.

We provide the $L_1$ and $L_\infty$ errors in the entropy in Figure \ref{fig:ringleb}b) and d), where the entropy is given by $P/\rho^\gamma$.  
We observe the expected $p+1$ rate of convergence in the $L_1$ norm and a rate between $p$ and $p+1$ in the $L_\infty$ norm.
This test problem was solved in \cite{WANG2006154}, where their spectral volume scheme on unstructured triangles may form transient shocks that do not allow convergence of the method.  
For all meshes considered here, our solver converged and the stopping criterion \eqref{eq:stopping} was attained.
Finally, in \cite{coirier1995accuracy} a second order finite volume method solved a similar Ringleb flow problem on cut cell grids, and observed a rate of convergence of 2.02 in $L_1$ and 1.40 in $L_\infty$.  

\subsection{Scattering of a smooth pressure pulse}
In this example, we solve the Euler equations \eqref{eq:euler} for the scattering of a pressure pulse surrounded by five circular obstacles.  The computational domain is $[0,20]^2$ and the obstacles are centered at $(10 + 5 \cos(\varphi_k), 10 + 5 \sin(\varphi_k) )$ for $\varphi_k = 2 \pi k/5 - 2\pi/3$.
The initial condition is given by
$$
\begin{pmatrix}
\rho\\
u \\
v \\
P
\end{pmatrix} = \begin{pmatrix}
1-1/\gamma + P\\
0\\
0\\
1/\gamma + 10^{-4} \exp( - b\left[ (x-10)^2 + (y-10)^2 \right] )
\end{pmatrix},
$$
where $\gamma = 1.4$, $b = \log(2)/0.2^2$, and the disturbance in the pressure is a Gaussian located at $(10,10)$.
The domain is discretized into a $53 \times 53$ cut cell grid, where each boundary segment is a polynomial of degree five ($q = 5$).
On this grid, the minimum volume fraction, $\min_{i,j} |K_{i,j}|/(\Delta x \Delta y)$, is 7.28e-06.
Roe's numerical flux is used on the interior element interfaces, the flux along the reflecting obstacles is given by \eqref{eq:refl}, and the ghost state on the outer boundaries of the domain is set to quiescent gas, $\rho = 1$, $p = 1/\gamma$, and $u = v = 0$.

\begin{figure}
\centering
\begin{tikzpicture}

\node (11) at (0, 0.175){\includegraphics[width = 0.36\textwidth]{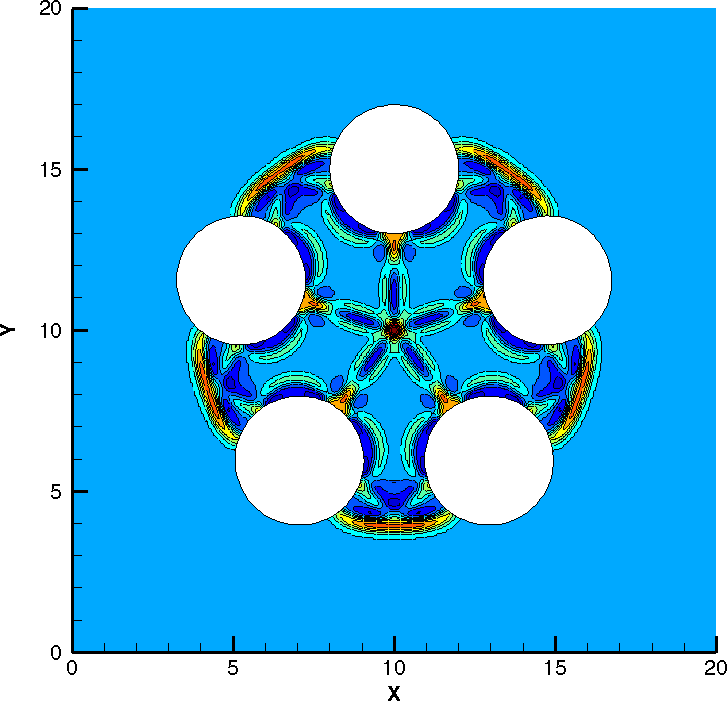}};
\node at (0.3,-3.3) {\textcolor{black}{\large (a)}};
\node (12) at (7.8,0.15){\includestandalone[width = 0.55\textwidth]{9-ap/boundary1}};
\node at (8.3,-3.3) {\textcolor{black}{\large (b)}};

\node at (-0.95,0.7) {\textcolor{black}{\large $\bm{\times}$}};
    \end{tikzpicture}
    \caption{
    In a), we plot the scaled deviation from the mean pressure, $(P-1/\gamma)/10^{-4}$, at the final time.
    In b), we provide this pressure deviation along the obstacle indicated by a $\bm{\times}$.
    Segments corresponding to different elements are plotted in different colors.
    We observe that the numerical solution retains the five-fold symmetry of the domain.
    On the embedded boundary, we observe the complex variation of the pressure in the scattered pulse.
    } \label{fig:ap}
\end{figure}
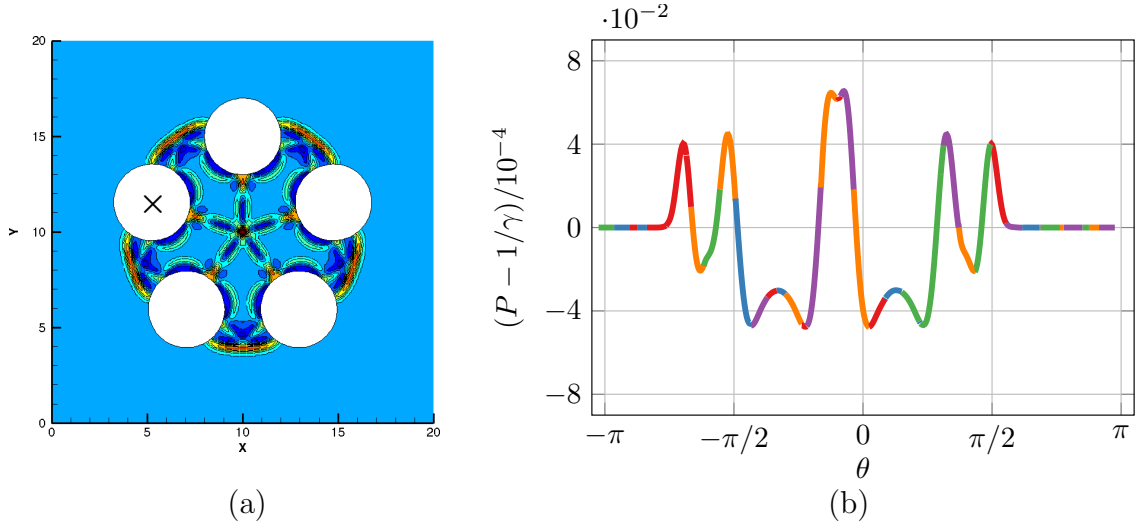
We provide the DG-$P5Q5$ numerical solution at the final time, $T = 6$ in Figure \ref{fig:ap}a) and c), respectively.  The solution presents the same five-fold symmetry as the obstacles.
The solution along the obstacle indicated by $\bm{\times}$ is plotted in Figures \ref{fig:ap}b) and d), where we observe sub-cell resolution of the complex scattering phenomenon.
Here, sub-cell resolution means that complex features of the solution are resolved entirely within an element due to high polynomial orders of approximation.
\subsection{Double Mach reflection}

In this example, we solve a shock-reflection problem where a vertical rightward moving Mach 10 shock reflects obliquely off of a solid wall.  
We solve \eqref{eq:euler} on $[0,2.5]\times [0,1.75]$, where a reflecting wedge that forms a $30^{\circ}$ angle with the horizontal begins at $(1/6, 0)$.
The initial position of the shock is $x = 1/6$ and its left and right states are
$$
\begin{pmatrix}
\rho\\
u \\
v \\
p
\end{pmatrix}_{\text{left}} = \begin{pmatrix} 8 \\ 8.25 \\ 0 \\ 116.5 \end{pmatrix} \text{ and } \begin{pmatrix}
\rho\\
u \\
v \\
p
\end{pmatrix}_{\text{right}} = \begin{pmatrix} 1.4 \\ 0 \\ 0 \\ 1 \end{pmatrix}.
$$
On the left and right boundaries, we apply the post and pre-shock states as boundary conditions, respectively.
On the top boundary, we use constant extrapolation.
The flux at the reflecting bottom boundary is given by \eqref{eq:refl} and the local Lax-Friedrichs numerical flux is used to compute the flux at the interior element interfaces.  
We use a formally second order DG method ($p=1$) to compute the solution at the final time $T=0.2$ on a sequence of grids where $\Delta x = \Delta y = 1/120, 1/240, 1/480$.
On these three grids, the minimum volume fraction ($K_{i,j}/\Delta x \Delta y$) was 4.11e-06, 4.11e-06, and 2.95e-07, respectively.
Since shocks are present in the solution, we must limit the gradient on each cell in characteristic variables. 
On the cells with a regular stencil, we use the MC limiter.  For the $x$-component of the gradient, the limited slope is
$$
\tilde{\mathbf{U}}_{x,i,j} = R_{i,j} ~\text{minmod}(2R_{i,j}^{-1}[\mathbf{c}_{i+1,j,0}-\mathbf{c}_{i,j,0}]/\Delta x, ~R_{i,j}^{-1}\mathbf{U}_{x,i,j}, ~2R_{i,j}^{-1}[\mathbf{c}_{i,j,0}-\mathbf{c}_{i-1,j,0}]/\Delta x),
$$
where $R_{i,j}$ is matrix of right eigenvectors of the flux Jacobian in the $x$-direction ($\partial \mathbf{F}/\partial \mathbf{U}\cdot[1,0]$) evaluated at the cell average of $K_{i,j}$.  An analogous formula is used to limit the slope in the $y$-direction.
On cells with an irregular stencil, we use a characteristic Barth-Jespersen type limiter \cite{char_lim} to compute a limited gradient
$$
\nabla \tilde{\mathbf{U}}_{i,j} = R_{i,j} \Phi[R^{-1}_{i,j} \nabla \mathbf{U}_{i,j}],
$$
where $R_{i,j}$ here is the matrix of right eigenvectors of the flux Jacobian evaluated at the cell average in the direction parallel the boundary: $\partial \mathbf{F}/\partial \mathbf{U} \cdot [\sqrt{3}/2, 1/2] $ for $x>1/6$, $\partial \mathbf{F}/\partial \mathbf{U} \cdot [1, 0] $ otherwise.  Finally, $\Phi$ applies Barth-Jespersen componentwise to the characteristic variables.  This limiter modifies the gradient of cells on the base grid and merging neighborhoods such that when the DG solution in characteristic variables is evaluated at neighboring centroids, the solution lies within the range defined by the cell averages of the neighbors \cite{berger2020state}. 
We have tested limiting in conserved variables and found that this produced unsatisfactory and oscillatory results.
In conjunction with these limiters, we implement the limiter in \cite{ZHANG20108918} to prevent negative densities and pressures at the quadrature points on the base grid and neighborhoods.

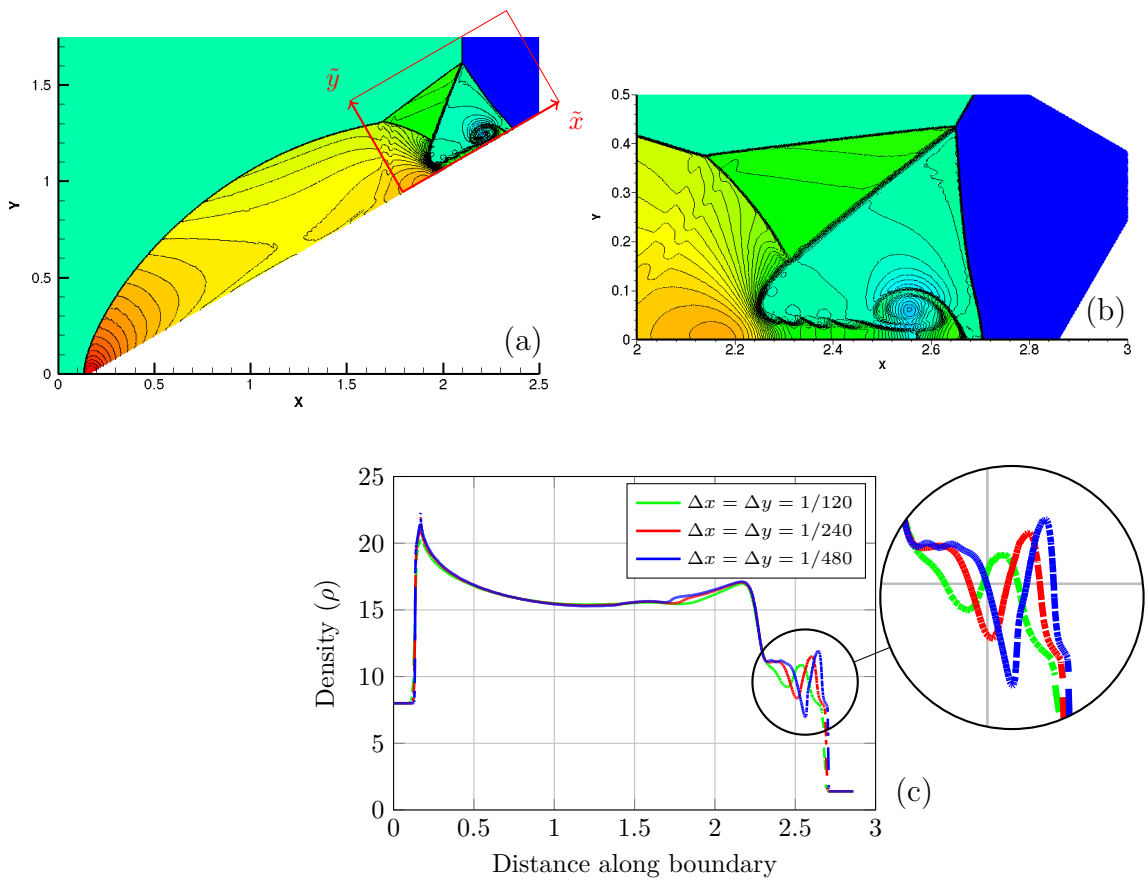
\begin{figure}
\pgfdeclarelayer{bg}    %
\pgfsetlayers{main,bg}  %

\begin{tikzpicture}[spy using outlines={circle, magnification=2.5, connect spies}]
\node (12) at (3.75,0) {\includegraphics[width = 0.45\linewidth]{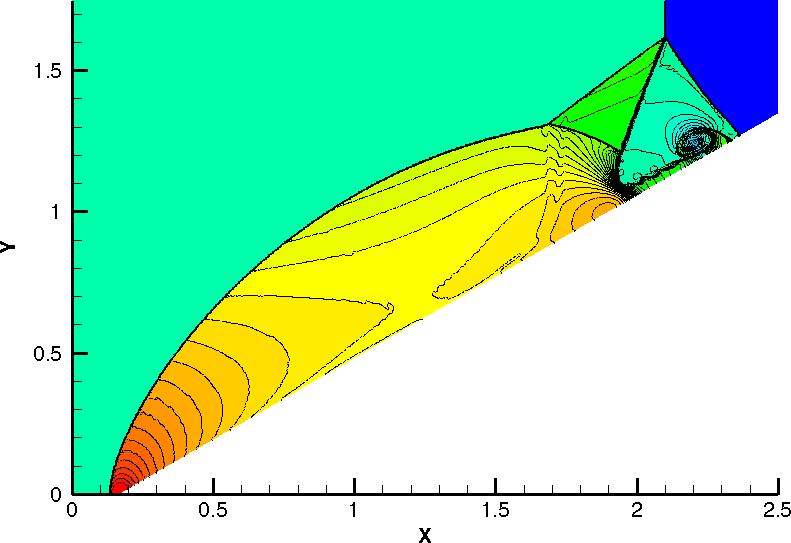}}; %
\node at (7, -1.6) {\textcolor{black}{\large (a)}};

\draw[red, line width=0.05mm,rotate around={30:(3.75+1/6, 0)}] (5.4,-0.4) rectangle (7.8,1);
\node (22) at (11.5,-0.5+0.4) {\includegraphics[width = 0.45\linewidth]{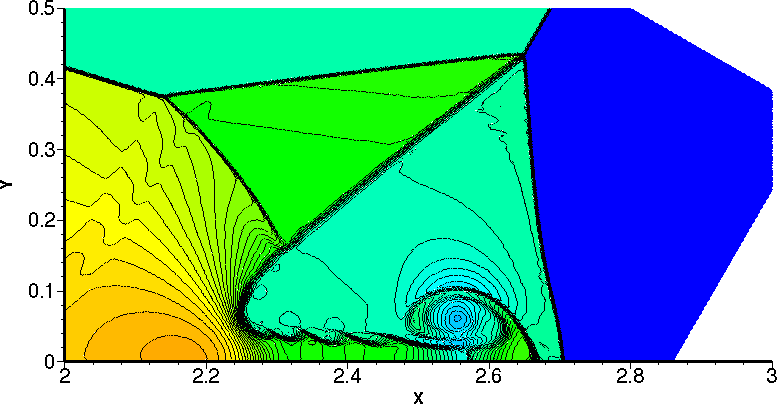}};
\node at (14.75,-1.6+0.4) {\textcolor{black}{\large (b)}};
    
  \begin{pgfonlayer}{bg} 
  \draw[->,rotate around={30:(3.75+1/6, 0)}, red, thick] (5.4,-0.4) -- (7.8,-0.4) coordinate[label=below right:$\tilde x$] ;
  \draw[->,rotate around={30:(3.75+1/6, 0)}, red, thick] (5.4,-0.4) -- (5.4,1) coordinate[label=above left:$\tilde y$] ;
  \end{pgfonlayer}

\node (11) at (8,-6.){\includegraphics[width = 0.49\textwidth]{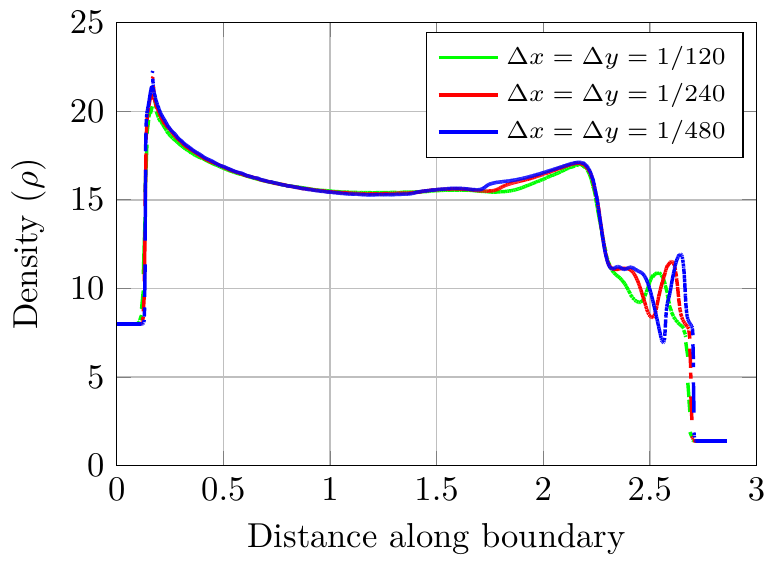} };
\spy [black, every spy on node/.append style={thick}, size=3.5cm] on (10.75,-6.125) in node[fill=white] at (13.5,-5);
\node at (12.2,-7.6) {\textcolor{black}{\large (c)}};
  \end{tikzpicture}
  \caption{Figures a) and b) show the density of the solution at the final time when $\Delta x = \Delta y = 1/480$, with 60 isolines between 2 and 20. In c), the solution along the embedded boundary measured from the origin are plotted for various grid resolutions.  Despite the highly irregular nature of grid at the boundary and the presence of shocks, the numerical solution is smooth.}\label{fig:dm}
\end{figure}

The numerical solution is a rich reflection pattern (Figure \ref{fig:dm}), comprised of incident and reflected shocks, a Mach stem, and contact discontinuity. 
This is a difficult test problem due to the strong shock and because the carbuncle phenomenon may occur at the boundary \cite{kemm2018heuristical}.
We obtain comparable results to those computed using a DG scheme on quadrilaterals and triangles \cite{COCKBURN1998199} and a finite volume scheme on cut cell grids \cite{berger2020state}.
The solution along the boundary is smooth and nonoscillatory (Figure \ref{fig:dm}c).

\section{Conclusions}
We have presented a state redistribution method for discontinuous Galerkin methods on curvilinear embedded boundary grids.
This technique preserves the formal order of accuracy of the base DG scheme and is conservative.
The advantage of state redistribution is that it only uses basic mesh information that is already available in many cut cell codes and does not require complex geometric computations.
Numerical examples reveal that the order of accuracy of the DG scheme is $p+1$ in the $L_1$ norm and between $p$ and $p+1$ in $L_\infty$.
We find that the rate in $L_\infty$ depends on the problem, the geometry of the embedded boundary, and the polynomial degree of approximation $p$.
This is a price that we are willing to pay in exchange for automatic mesh generation.
Future work includes studying more closely the loss of accuracy at the boundary and attempting to improve it.
We will also extend the method to three dimensions and more realistic engineering geometries.
This will require dealing with severely ill-shapen cut cells that commonly occur in three dimensions.
Finally, we will extend our method to viscous and incompressible flows.

\section{Code availability} \label{sec:codes}
The one-dimensional DG-SRD implementation in Python is available at \url{https://github.com/andrewgiuliani/PyDGSRD1D}.  This code will reproduce the convergence tests in the supplementary material document and can be applied to other conservation laws on nonuniform grids.
We have also made our cut cell mesh generation code written in Python available at \url{https://github.com/andrewgiuliani/PyGrid2D}.  
This code can be used to reproduce all the grids used in the numerical experiments in section \ref{sec:numex}, as well as generate other high order embedded boundary grids in two dimensions.

\section{Supplementary convergence tests and analysis}
In the supplementary material of this paper, we include 
\begin{itemize}
    \item[-] \textit{one-dimensional convergence tests on nonuniform grids}: two different sequences of nonuniform grids are considered, where the optimal $p+1$ rate of convergence is observed in both $L_1$ and $L_\infty$ errors.  For one-dimensional problems, optimal rates of convergence in $L_\infty$ have been observed before \cite{berger2020state, dod2}.  This is in contrast to two-dimensional methods which often suffer from a loss of accuracy at the embedded boundary as illustrated in this work.
    \item[-] \textit{proofs of claims 1 and 2}: we prove that the scheme is both \textit{p}-exact and conservative.  \textit{P}-exactness follows from the projections used in state redistribution.  The conservation proof in \cite{berger2020state} is adapted to DG methods.
    \item[-] \textit{implementation details to avoid precomputing the set $W_{i,j}$}: we adapt the approach in \cite{berger2020state} to DG methods.
\end{itemize}

\section{Acknowledgements}
The author would like to thank Marsha Berger, Sandra May, Georg Stadler, and Lilia Krivodonova for their reading of the manuscript and helpful comments.
We would also like to thank Marco Vianello for the helpful discussions on using the quadrature rule generation code.
AG is partially supported by the Natural Sciences and Engineering Research Council of Canada (NSERC) under Award No.\ PDF-546085-2020 and by the Simons Foundation/SFARI (560651, AB).

\FloatBarrier
\bibliography{references}
\bibliographystyle{ieeetr}

\end{document}

%% file: 1-dg/example_mesh.tex
  \begin{tikzpicture}
  
  \def\xx{3.5}
  \def\yy{-6}
  \def\R{8.5}
  \tikzset{
  reverseclip/.style={insert path={(-4,-4) rectangle (4cm,4cm)}}
  }

    \tikzdeclarepattern{
        name=new lines,
        parameters={
            \pgfkeysvalueof{/pgf/pattern keys/size},
            \pgfkeysvalueof{/pgf/pattern keys/angle},
            \pgfkeysvalueof{/pgf/pattern keys/line width},
        },
        bounding box={
            (
                0,
                -0.5*\pgfkeysvalueof{/pgf/pattern keys/line width}
            ) and (
                \pgfkeysvalueof{/pgf/pattern keys/size},
                0.5*\pgfkeysvalueof{/pgf/pattern keys/line width}
            )
        },
        tile size={
            (
                \pgfkeysvalueof{/pgf/pattern keys/size},
                \pgfkeysvalueof{/pgf/pattern keys/size}
            )
        },
        tile transformation={
            rotate=\pgfkeysvalueof{/pgf/pattern keys/angle}
        },
        defaults={
          size/.initial=4pt,
          angle/.initial=45,
          line width/.initial=0.5pt,
        },
        code={
            \draw [
                line width=\pgfkeysvalueof{/pgf/pattern keys/line width},
            ] (0,0) --
              (
                \pgfkeysvalueof{/pgf/pattern keys/size},
                0
              )
            ;
        },
    }

\clip (-2.1,-0.1) rectangle (1.1,3.1);

\begin{scope}

\clip[draw,reverseclip] (\xx,\yy) circle (\R);
	
\definecolor{beaublue}{rgb}{0.74, 0.83, 0.9};
\draw (-2,3) rectangle node [text=black, opacity=1, scale = 0.75]{$K_{i,j+1}$}  (-1,2) ;
\draw [draw=beaublue, fill=beaublue, fill opacity=  0.5] (-2,1) rectangle node [text=black, opacity=1]{}  (-1,0) ;
\draw [draw=beaublue, fill=beaublue, fill opacity=  0.5] (-3,1) rectangle node [text=black, opacity=1]{}  (-2,0) ;
\draw [draw=beaublue, fill=beaublue, fill opacity=  0.5] (-2,2) rectangle node [text=black, opacity=1, scale = 0.75]{$K_{i,j}$}  (-1,1) ;
\draw [draw=beaublue, fill=beaublue, fill opacity=  0.5] (-1,2) rectangle (0,1) ;
\node [text=black, opacity=1, scale = 0.6] at (-0.6,1.75)  {$K_{i+1,j}$} ;
\draw [draw=beaublue, fill=beaublue, fill opacity=  0.5] (0,2) rectangle node [text=black, opacity=1]{} (1,1) ;
\draw [draw=beaublue, fill=beaublue, fill opacity=  0.5] (0,3) rectangle node [text=black, opacity=1]{} (1,2) ;
\draw [draw=beaublue, fill=beaublue, fill opacity=  0.5] (1,3) rectangle node [text=black, opacity=1]{} (2,2) ;
\draw (-4,-4) grid (4,4);

\draw [black,ultra thick,domain=90:180, samples = 100] plot ({\R*cos(\x)+\xx}, {\R*sin(\x)+\yy});
\end{scope}
\begin{scope}
\clip[draw] (\xx,\yy) circle (\R);
\draw[pattern={new lines[size=3pt,line width=0.35pt,angle=125]}, pattern color=gray] (-2.1,-0.1) rectangle (1.1,3.1);
\draw [black,ultra thick,domain=90:180, samples = 100] plot ({\R*cos(\x)+\xx}, {\R*sin(\x)+\yy});

\end{scope}
  \end{tikzpicture}

%% file: 1-dg/cc.tex
\begin{tikzpicture}

\coordinate (A) at (0.75,0);
\coordinate (B) at (0.8609,0.1391);
\coordinate (C) at (1,0.25);

\definecolor{beaublue}{rgb}{0.74, 0.83, 0.9};
\draw [fill = beaublue, fill opacity = 0.5, ultra thin] (A) -- (0,0) -- (0,1) -- (1,1) -- (C) -- (B) ;

\draw [ultra thin,fill=black] (0,0) ;
\draw [ultra thin,fill=black] (0,1) ;
\draw [ultra thin,fill=black] (1,1) ;
\draw [ultra thin,fill=black] (A) circle (0.4pt);
\draw [ultra thin,fill=black] (B) circle (0.4pt);
\draw [ultra thin,fill=black] (C) circle (0.4pt);
\draw plot [smooth, tension=1, thick] coordinates {(A) (B) (C)};

\draw [ultra thin, fill=red,red,fill opacity=1] (0.100000000000000,0.100000000000000) circle (0.4pt);
\draw [ultra thin,fill=red,red, fill opacity=1] (0.100000000000000,0.500000000000000) circle (0.4pt);
\draw [ultra thin,fill=red,red, fill opacity=1] (0.100000000000000,0.800000000000000) circle (0.4pt);
\draw [ultra thin,fill=red,red, fill opacity=1] (0.100000000000000,0.900000000000000) circle (0.4pt);
\draw [ultra thin,fill=red,red, fill opacity=1] (0.200000000000000,0.100000000000000) circle (0.4pt);
\draw [ultra thin,fill=red,red, fill opacity=1] (0.200000000000000,0.900000000000000) circle (0.4pt);
\draw [ultra thin,fill=red,red, fill opacity=1] (0.400000000000000,0.200000000000000) circle (0.4pt);
\draw [ultra thin,fill=red,red, fill opacity=1] (0.500000000000000,0.100000000000000) circle (0.4pt);
\draw [ultra thin,fill=red,red, fill opacity=1] (0.500000000000000,0.500000000000000) circle (0.4pt);
\draw [ultra thin,fill=red,red, fill opacity=1] (0.500000000000000,0.900000000000000) circle (0.4pt);
\draw [ultra thin,fill=red,red, fill opacity=1] (0.700000000000000,0.400000000000000) circle (0.4pt);
\draw [ultra thin,fill=red,red, fill opacity=1] (0.800000000000000,0.200000000000000) circle (0.4pt);
\draw [ultra thin,fill=red,red, fill opacity=1] (0.800000000000000,0.800000000000000) circle (0.4pt);
\draw [ultra thin,fill=red,red, fill opacity=1] (0.900000000000000,0.600000000000000) circle (0.4pt);
\draw [ultra thin,fill=red,red, fill opacity=1] (0.900000000000000,0.900000000000000) circle (0.4pt);

  \end{tikzpicture}

%% file: 2-srd/nonuniform1.tex
  \begin{tikzpicture}

   \def \a{0.6}
   \def \h{1}
   
   \draw (-5*\h -\h*\a, 0.0) -- (5*\h +\h*\a,0.0);
   \draw[thick] (-2.5*\h -\h*\a, 0.0) -- (2.5*\h +\h*\a,0.0);

   \draw[thick] (-2.5*\h -\h*\a,-0.1) -- (-2.5*\h -\h*\a,0.1);
   \draw[thick] (-1.5*\h -\h*\a,-0.1) -- (-1.5*\h -\h*\a,0.1);
   \draw[thick] (-0.5*\h -\h*\a,-0.1) -- (-0.5*\h -\h*\a,0.1);
   \draw[thick] (-0.5*\h,-0.1) -- (-0.5*\h,0.1);
   \draw[thick] ( 0.5*\h,-0.1) -- ( 0.5*\h,0.1);
   \draw[thick] (0.5*\h +\h*\a,-0.1) -- (0.5*\h +\h*\a,0.1);
   \draw[thick] (1.5*\h +\h*\a,-0.1) -- (1.5*\h +\h*\a,0.1);
   \draw[thick] (2.5*\h +\h*\a,-0.1) -- (2.5*\h +\h*\a,0.1);
   
   \node at (-2*\h -\h*\a,0.4){$K_{ \shortminus 3}$};
   \node at (-1*\h -\h*\a,0.4){$K_{\shortminus 2}$};
   \node at (-0.15*\h-\h*\a,0.4){$K_{\shortminus 1}$};
   \node at (0, 0.4){$K_{0}$};
   \node at (0.15*\h+\h*\a,0.4){$K_{1}$};
   \node at (\h + \h*\a,0.4){$K_{2}$};
   \node at (2*\h + \h*\a,0.4){$K_{3}$};

   \node at (-2*\h -\h*\a,-0.4){$h$};
   \node at (-1*\h -\h*\a,-0.4){$h$};
   \node at (-0.15*\h-\h*\a,-0.4){$\alpha h$};
   \node at (0, -0.4){$h$};
   \node at (0.15*\h+\h*\a,-0.4){$\alpha h$};
   \node at (\h + \h*\a,-0.4){$h$};
   \node at (2*\h + \h*\a,-0.4){$h$};

    \node at (4.25*\h + \h*\a,-0.4){cell sizes};
  \end{tikzpicture}

%% file: 2-srd/nonuniform2.tex
  \begin{tikzpicture}

  \definecolor{airforceblue}{rgb}{0.36, 0.54, 0.66}
   
   \def \a{0.6}
   \def \h{1}

   \draw (-5*\h -\h*\a, 0.0) -- (5*\h +\h*\a,0.0);
   \draw[thick] (-2.5\h -\h*\a, 0.0) -- (2.5\h +\h*\a,0.0);
   
   \draw[thick,<->, red!70!black] (-1.5\h -\h*\a, -0.8)  --  (0.5\h,-0.8) node[draw=none,fill=none,midway,below] {$\hat K_{\shortminus 1}$};
   \draw[thick,<->, red!70!black] (-0.5\h, -1.5) -- (1.5\h +\h*\a,-1.5)  node[draw=none,fill=none,midway,below] {$\hat K_{1}$};

   \draw[thick,<->, airforceblue] (-2.5*\h -\h*\a, 0.8)  --  (-1.5*\h -\h*\a-0.05*\h,0.8) node[draw=none,fill=none,midway,above] {$\hat K_{\shortminus 3}$};
   \draw[thick,<->, airforceblue] (-1.5*\h -\h*\a+0.05*\h, 0.8)  --  (-0.5*\h -\h*\a,0.8) node[draw=none,fill=none,midway,above] {$\hat K_{\shortminus 2}$};
   \draw[thick,<->, airforceblue] (-0.5*\h, 0.8)  --  ( 0.5*\h,0.8) node[draw=none,fill=none,midway,above] {$\hat K_{0}$};
   \draw[thick,<->, airforceblue] (1.5*\h +\h*\a-0.05*\h, 0.8)  --  (0.5*\h +\h*\a,0.8) node[draw=none,fill=none,midway,above] {$\hat K_{2}$};
   \draw[thick,<->, airforceblue] (2.5*\h +\h*\a, 0.8)  --  (1.5*\h +\h*\a+0.05*\h,0.8) node[draw=none,fill=none,midway,above] {$\hat K_{3}$};
   
   \draw[thick] (-2.5*\h -\h*\a,-0.1) -- (-2.5*\h -\h*\a,0.1);
   \draw[thick] (-1.5*\h -\h*\a,-0.1) -- (-1.5*\h -\h*\a,0.1);
   \draw[thick] (-0.5*\h -\h*\a,-0.1) -- (-0.5*\h -\h*\a,0.1);
   \draw[thick] (-0.5*\h,-0.1) -- (-0.5*\h,0.1);
   \draw[thick] ( 0.5*\h,-0.1) -- ( 0.5*\h,0.1);
   \draw[thick] (0.5*\h +\h*\a,-0.1) -- (0.5*\h +\h*\a,0.1);
   \draw[thick] (1.5*\h +\h*\a,-0.1) -- (1.5*\h +\h*\a,0.1);
   \draw[thick] (2.5*\h +\h*\a,-0.1) -- (2.5*\h +\h*\a,0.1);
   
   \node at (-2*\h -\h*\a,0.4){$K_{\shortminus 3}$};
   \node at (-1*\h -\h*\a,0.4){$K_{\shortminus 2}$};
   \node at (-0.15*\h-\h*\a,0.4){$K_{\shortminus 1}$};
   \node at (0, 0.4){$K_{0}$};
   \node at (0.15*\h+\h*\a,0.4){$K_{1}$};
   \node at (\h + \h*\a,0.4){$K_{2}$};
   \node at (2*\h + \h*\a,0.4){$K_{3}$};

   \node at (-2*\h -\h*\a,-0.4){$1$};
   \node at (-1*\h -\h*\a,-0.4){$2$};
   \node at (-0.15*\h-\h*\a,-0.4){$1$};
   \node at (0, -0.4){$3$};
   \node at (0.15*\h+\h*\a,-0.4){$1$};
   \node at (\h + \h*\a,-0.4){$2$};
   \node at (2*\h + \h*\a,-0.4){$1$};
   
   \node at (3.75*\h + \h*\a,-0.4){overlap counts};

  \end{tikzpicture}

%% file: 2-srd/ortho_m0.tex
  \begin{tikzpicture}
  
  \pgfplotstableread[col sep=comma]{2-srd/m1.txt}\mone
  \pgfplotstableread[col sep=comma]{2-srd/legendreN.txt}\leg
  
\pgfplotsset{ every non boxed x axis/.append style={x axis line style=-},
     every non boxed y axis/.append style={y axis line style=-}}
     
\tikzfading[name=fade down,
    right color=transparent!0,
    left color=transparent!50]     

   \begin{axis}[xmin=-1.25,xmax=1.25,
     ymin=-5,ymax=5,
     axis x line=center,
     xtick={-1.15,-0.15,0.15,1.15},
     xticklabels={,,},
     yticklabels={,,},
     hide y axis,
     axis line style = thick,
     every major tick/.append style={very thick, major tick length=10pt, black},
     grid=major,
     major grid style={gray, thick},
     cycle list/Dark2,
     cycle list={[indices of colormap={0,1,2,3,4,1,2,3,4 of Dark2}]}
   ]
        \addplot+[mark=none, thick] table[x=x, y=q1] {\mone};
        \addplot+[mark=none, thick] table[x=x, y=q2] {\mone};
        \addplot+[mark=none, thick] table[x=x, y=q3] {\mone};
        \addplot+[mark=none, thick] table[x=x, y=q4] {\mone};
        \addplot+[mark=none, thick] table[x=x, y=q5] {\mone};

        \addplot+[mark=none,dotted, thick,path fading=fade down] table[x=x, y=q2] {\leg};
        \addplot+[mark=none,dotted, thick,path fading=fade down] table[x=x, y=q3] {\leg};
        \addplot+[mark=none,dotted, thick,path fading=fade down] table[x=x, y=q4] {\leg};
        \addplot+[mark=none,dotted, thick,path fading=fade down] table[x=x, y=q5] {\leg};
        
        \draw [stealth-,thick] (axis cs:0,1) -- (axis cs:0,1.5)  node[above]{\small $\hat \varphi_{\shortminus 1,0}$} ;
        \draw [stealth-,thick] (axis cs:-0.95,-1.4) -- (axis cs:-0.95,-1.9)  node[below]{\small $\hat \varphi_{\shortminus 1,1}$} ;
        \draw [stealth-,thick] (axis cs:-0.97,1.3) -- (axis cs:-0.97,1.8)  node[above]{\small $\hat \varphi_{\shortminus 1,2}$} ;
        \draw [stealth-,thick] (axis cs:-0.535,1.24) -- (axis cs:-0.535,1.74)  node[above]{\small $\hat \varphi_{\shortminus 1,3}$} ;
        \draw [stealth-,thick] (axis cs:-0.65,-1.15) -- (axis cs:-0.65,-1.65)  node[below]{\small $\hat \varphi_{\shortminus 1,4}$} ;

        \draw  (axis cs:-0.6,4) node[above]{$K_{\shortminus 2}$} ;
        \draw  (axis cs:0,4) node[above]{ $K_{\shortminus 1}$} ;
        \draw  (axis cs:0.6,4) node[above]{$K_{0}$} ;
        \draw[thick,<->, red!70!black] (axis cs:-1.15, -3.75)  --  (axis cs:1.15,-3.75) node[draw=none,fill=none,midway,below] {$\hat K_{\shortminus 1}$};

   \end{axis}

  \end{tikzpicture}

%% file: 2-srd/ortho_m1.tex
  \begin{tikzpicture}
  
  \pgfplotstableread[col sep=comma]{2-srd/m0.txt}\mone
\pgfplotsset{ every non boxed x axis/.append style={x axis line style=-},
     every non boxed y axis/.append style={y axis line style=-}}
   \begin{axis}[xmin=-1.1,xmax=1.1,
     ymin=-5,ymax=5,
     axis x line=center,
     xtick={-1,1},
     xticklabels={,,},
     yticklabels={,,},
     hide y axis,
     axis line style = thick,
     every major tick/.append style={very thick, major tick length=10pt, black},
     grid=major,
     major grid style={gray, thick},
     cycle list/Dark2
   ]
        \addplot+[mark=none, thick] table[x=x, y=q1] {\mone};
        \addplot+[mark=none, thick] table[x=x, y=q2] {\mone};
        \addplot+[mark=none, thick] table[x=x, y=q3] {\mone};
        \addplot+[mark=none, thick] table[x=x, y=q4] {\mone};
        \addplot+[mark=none, thick] table[x=x, y=q5] {\mone};

        \draw  (axis cs:0,4) node[above]{$K_{\shortminus 1}$} ;

        \draw [stealth-,thick] (axis cs:-0.175,1) -- (axis cs:-0.175,1.5)  node[above]{\small $\hat \varphi_{\shortminus 1,0}$} ;
        \draw [stealth-,thick] (axis cs:-0.95+0.15,-1.4) -- (axis cs:-0.95+0.15,-1.9)  node[below]{\small $\hat \varphi_{\shortminus 1,1}$} ;
        \draw [stealth-,thick] (axis cs:-0.97+0.125,1.3) -- (axis cs:-0.97+0.125,1.8)  node[above]{\small $\hat \varphi_{\shortminus 1,2}$} ;
        \draw [stealth-,thick] (axis cs:-0.535+0.1,1.24) -- (axis cs:-0.535+0.1,1.74)  node[above]{\small $\hat \varphi_{\shortminus 1,3}$} ;
        \draw [stealth-,thick] (axis cs:-0.65+0.1,-1.15) -- (axis cs:-0.65+0.1,-1.65)  node[below]{\small $\hat \varphi_{\shortminus 1,4}$} ;

   \end{axis}

  \end{tikzpicture}

%% file: 2-srd/neigh1.tex
  \begin{tikzpicture}
  \tikzset{
  reverseclip/.style={insert path={(-4,-4) rectangle (4cm,4cm)}}
}
\definecolor{beaublue}{rgb}{0.74, 0.83, 0.9};

\clip (-2.2,0.9) rectangle (2.4,3.2);
\begin{scope}

\clip[draw,reverseclip] (.25,-2) circle (4);
\draw (-5,-5) grid (5,5);

\draw [ultra thick, draw=green!70!black, fill=green!70!black, fill opacity=0.2] (-1,3) rectangle  (0,1);

\draw [ultra thick, draw=green!70!black, fill=green!70!black, fill opacity=0.2] (0,-3) rectangle (1,1);

\end{scope}
\begin{scope}
\clip[draw] (0.25,-2) circle (4);

\draw [ultra thick, draw=beaublue, fill=beaublue, fill opacity=0.2] (-1,2) rectangle node[text=black, opacity=1]{$K_{i,j}$} (0,1) ;

\end{scope}
  \end{tikzpicture}

%% file: 2-srd/neigh2.tex
  \begin{tikzpicture}
  \tikzset{
  reverseclip/.style={insert path={(-4,-4) rectangle (4cm,4cm)}}
}

\clip (-2.2,0.9) rectangle (2.4,3.2);
\begin{scope}
\clip[draw,reverseclip] (.25,-2) circle (4);
\draw (-5,-5) grid (5,5);

\draw [ultra thick, draw=green!70!black, fill=green!70!black, fill opacity=0.2] (-1,3) rectangle node[text=black, opacity=1]{$K_{i,j+1}$}  (0,2);
\end{scope}

  \end{tikzpicture}

%% file: 2-srd/neigh3.tex
  \begin{tikzpicture}
  \tikzset{
  reverseclip/.style={insert path={(-5,-2.1) rectangle (1.1,1.1)}}
}
\definecolor{beaublue}{rgb}{0.74, 0.83, 0.9};

\begin{scope}
\clip [draw](-4.5,-2)--(0,0.1)--(-4.5,2); 

\draw (-5,-5) grid (5,5);
\draw [ultra thick, draw=red!70!black, fill=red!70!black, fill opacity=0.2] (-1,-1) rectangle  (0,2);
\end{scope}

\begin{scope}
\clip [reverseclip](-4.5,-2)--(0,0.1)--(-4.5,2);
\draw [ultra thick, draw=beaublue, fill=beaublue, fill opacity=0.2] (-1,-1) rectangle node[text=black, opacity=1]{$K_{i',j'}$}  (0,0);
\draw [ultra thick, draw=beaublue, fill=beaublue, fill opacity=0.2] (-1,0) rectangle node[text=black, opacity=1]{}  (0,1);
\end{scope}

  \end{tikzpicture}

%% file: 2-srd/neigh4.tex
  \begin{tikzpicture}
  \tikzset{
  reverseclip/.style={insert path={(-5,-2.1) rectangle (1.1,1.1)}}
}
\definecolor{beaublue}{rgb}{0.74, 0.83, 0.9};
\begin{scope}
\clip [draw](-4.5,-2)--(0,0.1)--(-4.5,2); 

\draw (-5,-5) grid (5,5);
\draw [ultra thick, draw=green!70!black, fill=green!70!black, fill opacity=0.2] (-2,-2) rectangle  (0,2);
\end{scope}

\begin{scope}
\clip [reverseclip](-4.5,-2)--(0,0.1)--(-4.5,2);
\draw [ultra thick, draw=beaublue, fill=beaublue, fill opacity=0.2] (-1,-1) rectangle node[text=black, opacity=1]{$K_{i',j'}$}  (0,0);
\draw [ultra thick, draw=beaublue, fill=beaublue, fill opacity=0.2] (-1,0) rectangle node[text=black, opacity=1]{}  (0,1);
\draw [ultra thick, draw=beaublue, fill=beaublue, fill opacity=0.2] (-2,0) rectangle node[text=black, opacity=1]{}  (-1,1);
\draw [ultra thick, draw=beaublue, fill=beaublue, fill opacity=0.2] (-2,-1) rectangle node[text=black, opacity=1]{}  (-1,0);
\draw [ultra thick, draw=beaublue, fill=beaublue, fill opacity=0.2] (-2,-2) rectangle node[text=black, opacity=1]{}  (-1,-1);
\draw [ultra thick, draw=beaublue, fill=beaublue, fill opacity=0.2] (-1,-2) rectangle node[text=black, opacity=1]{}  (0,-1);
\draw [ultra thick, draw=beaublue, fill=beaublue, fill opacity=0.2] (0,-2) rectangle node[text=black, opacity=1]{}  (1,-1);
\draw [ultra thick, draw=beaublue, fill=beaublue, fill opacity=0.2] (0,-1) rectangle node[text=black, opacity=1]{}  (1,0);
\draw [ultra thick, draw=beaublue, fill=beaublue, fill opacity=0.2] (0,0) rectangle node[text=black, opacity=1]{}  (1,1);
\end{scope}

  \end{tikzpicture}

%% file: 2-srd/bounding_box.tex
  \begin{tikzpicture}
  \tikzset{
  reverseclip/.style={insert path={(-4,-4) rectangle (4cm,4cm)}}
}

\definecolor{beaublue}{rgb}{0.74, 0.83, 0.9};

\begin{scope}
\clip (-2.2,0.9) rectangle (2.4,3.2);
\clip[draw,reverseclip] (.25,-2) circle (4);
\draw (-5,-5) grid (5,5);
\draw [thick, draw=beaublue, fill=beaublue, fill opacity=0.2] (-1,2.) rectangle (0,1.8) ;
\draw [thick, draw=beaublue, fill=beaublue, fill opacity=0.2] (1,2.) rectangle (2,1.6) ;
\end{scope}

\draw [thick, draw=beaublue] (-1,2.) rectangle (0,1.8) ;
\draw [thick, draw=beaublue] (1,2.) rectangle (2,1.6) ;

\draw [thick, red,fill opacity=1] (-1,1.8) circle (1pt);
\node[below=1pt of {(-1,1.8)}, outer sep=4pt, scale = 0.6] {$(x,y)_{i,j,\min}$};

\draw [thick,red,fill opacity=1] (1,1.6) circle (1pt);
\node[below=1pt of {(1,1.6)}, scale = 0.6] {$(x,y)_{i+2,j,\min}$};

    \coordinate (ij_top_left)  at (-1,2.1);
    \coordinate (ij_top_right) at ( 0,2.1);
    \coordinate (avex) at ($(ij_top_left)!0.5!(ij_top_right)$);
    \draw [thin,<->] (ij_top_left) -- (ij_top_right);
    \node[above=0.01pt of {avex}, scale = 0.6] {$\Delta x_{i,j}$};
    
    \coordinate (ij_right_top)  at (0.1,1.8);
    \coordinate (ij_right_bot) at ( 0.1,2);
    \coordinate (avey) at ($(ij_right_top)!0.5!(ij_right_bot)+(0,-0.04)$);
    \draw [thin,<->, scale = 0.5] (ij_right_top) -- (ij_right_bot);
    \node[right=1pt of {avey}, scale = 0.6] {$\Delta y_{i,j}$};
    
    \coordinate (ij2_right_top)  at (2.1,1.6);
    \coordinate (ij2_right_bot) at ( 2.1,2);
    \coordinate (avey) at ($(ij2_right_top)!0.5!(ij2_right_bot)$);
    \draw [thin,<->] (ij2_right_top) -- (ij2_right_bot);
    \node[right=0.01pt of {avey}, scale = 0.6] {$\Delta y_{i+2,j}$};

    \coordinate (ij2_top_left)  at (1,2.1);
    \coordinate (ij2_top_right) at (2,2.1);
    \coordinate (ave2x) at ($(ij2_top_left)!0.5!(ij2_top_right)$);
    \draw [thin,<->] (ij2_top_left) -- (ij2_top_right);
    \node[above=0.01pt of {ave2x}, scale = 0.6] {$\Delta x_{i+2,j}$};
    
  \end{tikzpicture}

%% file: 2-srd/mapping.tex
  \begin{tikzpicture}
  
\node [scale = 0.3] (x) at (-0.1,1.2) {$\eta$};
\node [scale = 0.3] (y) at (1.2,-0.1) {$\xi$};
\draw [->] (0,0) -- (1.2,0);
\draw [->] (0,0) -- (0,1.2);
\tikzset{mygridstyle/.style={dash pattern=on .8pt off .4pt on .6pt off .4pt on .4pt off .4pt on .2pt off .4pt}}
\draw [step=0.1,gray, gray, mygridstyle, line width=0.01pt] (0,0) grid (1,1);

\newcommand*{\TickSize}{0.025}%

  \foreach \x in {0,0.2,...,1.} {%
    \draw  ($(\x,0) + (0,\TickSize)$) -- ($(\x,0) + (0,-\TickSize)$)
        node [below, scale = 0.3] {\pgfmathprintnumber[fixed, precision=1]{\x}};
}

  \foreach \y in {0,0.2,...,1.} {%
    \draw ($(0,\y) + (\TickSize,0)$) -- ($(0,\y) + (-\TickSize,0)$)
        node [left, scale = 0.3] {\pgfmathprintnumber[fixed, precision=1]{\y}};
}

    \definecolor{beaublue}{rgb}{0.74, 0.83, 0.9};
    \coordinate (A)  at (0,0.9);
    \coordinate (BBCC) at (0.15,0.84);
    \coordinate (BC) at (0.3,0.765);
    \coordinate (BBCB) at (0.45,0.675);
    \coordinate (B)  at (0.57,0.57);
    \coordinate (AABB)  at (0.7,0.44);
    \coordinate (AB) at (0.805,0.31);
    \coordinate (ABCC) at (0.92,0.15);
    \coordinate (C)  at (1,0);
    
    \draw [fill = beaublue, fill opacity = 0.2, draw = beaublue] (A) -- (0,1) -- (1,1) -- (C) -- (ABCC) -- (AB) -- (AABB) -- (B) -- (BBCB) -- (BC) -- (BBCC) -- (A) ;
    \draw [draw = beaublue] plot [smooth, tension=1] coordinates {(A) (B) (C)};

\draw [ultra thin, fill=red,red,fill opacity=1] (0.800000000000000,0.800714469370758) circle (0.4pt);
\draw [ultra thin,fill=red,red, fill opacity=1] (0.066666666666667,0.933571489790253) circle (0.4pt);
\draw [ultra thin,fill=red,red, fill opacity=1] (0.133333333333333,0.867142979580505) circle (0.4pt);
\draw [ultra thin,fill=red,red, fill opacity=1] (0.200000000000000,0.933571489790253) circle (0.4pt);
\draw [ultra thin,fill=red,red, fill opacity=1] (0.466666666666667,0.667857448951263) circle (0.4pt);
\draw [ultra thin,fill=red,red, fill opacity=1] (0.466666666666667,0.933571489790253) circle (0.4pt);
\draw [ultra thin,fill=red,red, fill opacity=1] (0.533333333333333,0.800714469370758) circle (0.4pt);
\draw [ultra thin,fill=red,red, fill opacity=1] (0.733333333333333,0.601428938741515) circle (0.4pt);
\draw [ultra thin,fill=red,red, fill opacity=1] (0.800000000000000,0.335714897902525) circle (0.4pt);
\draw [ultra thin,fill=red,red, fill opacity=1] (0.800000000000000,0.933571489790253) circle (0.4pt);
\draw [ultra thin,fill=red,red, fill opacity=1] (0.933333333333333,0.136429367273283) circle (0.4pt);
\draw [ultra thin,fill=red,red, fill opacity=1] (0.933333333333333,0.402143408112273) circle (0.4pt);
\draw [ultra thin,fill=red,red, fill opacity=1] (0.933333333333333,0.601428938741515) circle (0.4pt);
\draw [ultra thin,fill=red,red, fill opacity=1] (0.933333333333333,0.800714469370758) circle (0.4pt);
\draw [ultra thin,fill=red,red, fill opacity=1] (0.933333333333333,0.933571489790253) circle (0.4pt);

    \end{tikzpicture}

%% file: 2-srd/bounding_box_hat.tex
  \begin{tikzpicture}
  \tikzset{
  reverseclip/.style={insert path={(-4,-4) rectangle (4cm,4cm)}}
}

\definecolor{beaublue}{rgb}{0.74, 0.83, 0.9};

\begin{scope}
\clip (-2.2,0.9) rectangle (2.4,3.2);
\clip[draw,reverseclip] (.25,-2) circle (4);
\draw (-5,-5) grid (5,5);
\draw [thick, draw=beaublue, fill=beaublue, fill opacity=0.2] (-1,3.) rectangle (0,1.8) ;
\draw [thick, draw=beaublue, fill=beaublue, fill opacity=0.2] (1,3.) rectangle (2,1.6) ;
\end{scope}

\draw [thick, draw=beaublue] (-1,3.) rectangle (0,1.8) ;
\draw [thick, draw=beaublue] (1,3.) rectangle (2,1.6) ;

\draw [thick, red,fill opacity=1] (-1,1.8) circle (1pt);
\node[below=1pt of {(-1,1.8)}, outer sep=4pt, scale = 0.6] {$(\hat x,\hat y)_{i,j,\min}$};

\draw [thick,red,fill opacity=1] (1,1.6) circle (1pt);
\node[below=1pt of {(1,1.6)}, scale = 0.6] {$(\hat x,\hat y)_{i+2,j,\min}$};

    \coordinate (ij_top_left)  at (-1,2.1);
    \coordinate (ij_top_right) at ( 0,2.1);
    \coordinate (avex) at ($(ij_top_left)!0.5!(ij_top_right)$);
    \draw [thin,<->] (ij_top_left) -- (ij_top_right);
    \node[above=0.01pt of {avex}, scale = 0.6] {$\Delta \hat x_{i,j}$};
    
    \coordinate (ij_right_top)  at (0.1,1.8);
    \coordinate (ij_right_bot) at ( 0.1,3);
    \coordinate (avey) at ($(ij_right_top)!0.5!(ij_right_bot)+(0,-0.04)$);
    \draw [thin,<->] (ij_right_top) -- (ij_right_bot);
    \node[right=1pt of {avey}, scale = 0.6] {$\Delta \hat y_{i,j}$};
    
    \coordinate (ij2_right_top)  at (2.1,1.6);
    \coordinate (ij2_right_bot) at ( 2.1,3.);
    \coordinate (avey) at ($(ij2_right_top)!0.5!(ij2_right_bot)$);
    \draw [thin,<->] (ij2_right_top) -- (ij2_right_bot);
    \node[right=0.01pt of {avey}, scale = 0.6] {$\Delta \hat y_{i+2,j}$};

    \coordinate (ij2_top_left)  at (1,2.1);
    \coordinate (ij2_top_right) at (2,2.1);
    \coordinate (ave2x) at ($(ij2_top_left)!0.5!(ij2_top_right)$);
    \draw [thin,<->] (ij2_top_left) -- (ij2_top_right);
    \node[above=0.01pt of {ave2x}, scale = 0.6] {$\Delta \hat x_{i+2,j}$};
    
  \end{tikzpicture}

%% file: 2-srd/mapping_hat.tex
  \begin{tikzpicture}
  
\node [scale = 0.3] (x) at (-0.1,1.2) {$\eta$};
\node [scale = 0.3] (y) at (1.2,-0.1) {$\xi$};
\draw [->] (0,0) -- (1.2,0);
\draw [->] (0,0) -- (0,1.2);
\tikzset{mygridstyle/.style={dash pattern=on .8pt off .4pt on .6pt off .4pt on .4pt off .4pt on .2pt off .4pt}}
\draw [step=0.1,gray, gray, mygridstyle, line width=0.01pt] (0,0) grid (1,1);

\newcommand*{\TickSize}{0.025}%

  \foreach \x in {0,0.2,...,1.} {%
    \draw  ($(\x,0) + (0,\TickSize)$) -- ($(\x,0) + (0,-\TickSize)$)
        node [below, scale = 0.3] {\pgfmathprintnumber[fixed, precision=1]{\x}};
}

  \foreach \y in {0,0.2,...,1.} {%
    \draw ($(0,\y) + (0,0)$) -- ($(0,\y) + (-\TickSize,0)$)
        node [left, scale = 0.3] {\pgfmathprintnumber[fixed, precision=1]{\y}};
}

    \definecolor{beaublue}{rgb}{0.74, 0.83, 0.9};

\draw [fill = beaublue, fill opacity = 0.2, draw = beaublue] (0,0.3) -- (1,0.3) -- (1,1) -- (0,1) -- (0,0.3);
\draw [ultra thin, fill=red,red,fill opacity=1] (0.684174025172418,0.337410015715938) circle (0.4pt);
\draw [ultra thin,fill=red,red, fill opacity=1] (0.315825974827582,0.962589984284062) circle (0.4pt);
\draw [ultra thin,fill=red,red, fill opacity=1] (0.944091898161729,0.516971046508595) circle (0.4pt);
\draw [ultra thin,fill=red,red, fill opacity=1] (0.055908101838271,0.783028953491405) circle (0.4pt);
\draw [ultra thin,fill=red,red, fill opacity=1] (0.157536078259683,0.411535369785387) circle (0.4pt);
\draw [ultra thin,fill=red,red, fill opacity=1] (0.842463921740317,0.888464630214613) circle (0.4pt);
\draw [ultra thin,fill=red,red, fill opacity=1] (0.500000000000000,0.650000000000000) circle (0.4pt);

    \begin{scope}[yscale=0.3]
    \coordinate (A)  at (0,0.9);
    \coordinate (BBCC) at (0.15,0.84);
    \coordinate (BC) at (0.3,0.765);
    \coordinate (BBCB) at (0.45,0.675);
    \coordinate (B)  at (0.57,0.57);
    \coordinate (AABB)  at (0.7,0.44);
    \coordinate (AB) at (0.805,0.31);
    \coordinate (ABCC) at (0.92,0.15);
    \coordinate (C)  at (1,0);
    
    \draw [fill = beaublue, fill opacity = 0.2, draw = beaublue] (A) -- (0,1) -- (1,1) -- (C) -- (ABCC) -- (AB) -- (AABB) -- (B) -- (BBCB) -- (BC) -- (BBCC) -- (A) ;
    \draw [draw = beaublue] plot [smooth, tension=1] coordinates {(A) (B) (C)};
\end{scope}
\draw [ultra thin, fill=red,red,fill opacity=1] (0.800000000000000,0.240214340811227) circle (0.4pt);
\draw [ultra thin,fill=red,red, fill opacity=1] (0.066666666666667,0.280071446937076) circle (0.4pt);
\draw [ultra thin,fill=red,red, fill opacity=1] (0.133333333333333,0.260142893874152) circle (0.4pt);
\draw [ultra thin,fill=red,red, fill opacity=1] (0.200000000000000,0.280071446937076) circle (0.4pt);
\draw [ultra thin,fill=red,red, fill opacity=1] (0.466666666666667,0.200357234685379) circle (0.4pt);
\draw [ultra thin,fill=red,red, fill opacity=1] (0.466666666666667,0.280071446937076) circle (0.4pt);
\draw [ultra thin,fill=red,red, fill opacity=1] (0.533333333333333,0.240214340811227) circle (0.4pt);
\draw [ultra thin,fill=red,red, fill opacity=1] (0.733333333333333,0.180428681622455) circle (0.4pt);
\draw [ultra thin,fill=red,red, fill opacity=1] (0.800000000000000,0.100714469370758) circle (0.4pt);
\draw [ultra thin,fill=red,red, fill opacity=1] (0.800000000000000,0.280071446937076) circle (0.4pt);
\draw [ultra thin,fill=red,red, fill opacity=1] (0.933333333333333,0.040928810181984) circle (0.4pt);
\draw [ultra thin,fill=red,red, fill opacity=1] (0.933333333333333,0.120643022433682) circle (0.4pt);
\draw [ultra thin,fill=red,red, fill opacity=1] (0.933333333333333,0.180428681622455) circle (0.4pt);
\draw [ultra thin,fill=red,red, fill opacity=1] (0.933333333333333,0.240214340811227) circle (0.4pt);
\draw [ultra thin,fill=red,red, fill opacity=1] (0.933333333333333,0.280071446937076) circle (0.4pt);

    \end{tikzpicture}

%% file: 3-rh/l1.tex
\begin{tikzpicture}
    \begin{loglogaxis}[
      width=11cm,
      height=6cm,
      xlabel={$N$},
      ylabel={$L_1$ error},
      legend cell align=left,
      legend style={at={(0.01, 0.01)}, anchor=south west},
      legend columns=1,
      legend style={font=\small},
      scale only axis,
      grid=both,
      xmin=20,
      xmax = 1000,
      ymin = 1e-13,
      ymax = 1e0
      ]
      
\pgfplotstableread[col sep = comma]{
N,p1
50,  1.1236699955761462e-02
75,  4.2553738745174279e-03
100, 2.1342949094355578e-03
125, 1.2436271290036120e-03
150, 8.0711684627351435e-04
175, 5.6748344250471541e-04
200, 4.2403575558762796e-04
225, 3.2799954845417735e-04
250, 2.6163367128636633e-04
275, 2.1368844732972677e-04
300, 1.7858315069209179e-04
325, 1.5083737605008996e-04
350, 1.2932127023105111e-04
375, 1.1240123519465861e-04
400, 9.8412522837477027e-05
425, 8.6894740290603109e-05
450, 7.7183913161998609e-05
475, 6.9244994831546246e-05
500, 6.2348707033433475e-05
525, 5.6464458033735770e-05
550, 5.1313101576719463e-05
575, 4.6960661787049441e-05
600, 4.3065150838818475e-05
625, 3.9627311091006768e-05
650, 3.6566399797607023e-05
675, 3.3922741731727421e-05
700, 3.1512821159091163e-05
725, 2.9355576601398656e-05
750, 2.7413685565460307e-05
775, 2.5650561726103329e-05
800, 2.4053049952659452e-05
}\first

\pgfplotstableread[col sep = comma]{
N,p2
50,  2.2003637439405263e-04
75,  4.2674352336617549e-05
100, 1.5044462023846899e-05
125, 6.7224800403316736e-06
150, 3.5447753661430644e-06
175, 2.1437944666672796e-06
200, 1.3961859750619936e-06
225, 9.6055582822442374e-07
250, 6.8281427100693428e-07
275, 5.0752471163384680e-07
300, 3.8601584302796248e-07
325, 3.0121514240214753e-07
350, 2.3840343777991707e-07
375, 1.9291138653231761e-07
400, 1.5800410104800230e-07
425, 1.3116590577491603e-07
450, 1.0974134068442082e-07
475, 9.3107444875217251e-08
500, 7.9637073836675336e-08
525, 6.8532397642942296e-08
550, 5.9450531133056574e-08
575, 5.1854517831430222e-08
600, 4.5553521912171536e-08
625, 4.0236041639224326e-08
650, 3.5718158618776150e-08
675, 3.1823678449123065e-08
700, 2.8498631098519227e-08
725, 2.5586129395192181e-08
750, 2.3084210253222525e-08
775, 2.0884450541397941e-08
800, 1.8975134240757726e-08
}\second

\pgfplotstableread[col sep = comma]{
N,p3
50 , 6.8935530515055029e-06
75 , 1.0972081217303980e-06
100, 3.4106051019171835e-07
125, 1.3169017597847995e-07
150, 5.7444350255736307e-08
175, 3.0019649618085498e-08
200, 1.7866347282788121e-08
225, 1.1068027557819398e-08
250, 7.2789623016040204e-09
275, 4.9271447393150928e-09
300, 3.5094615580245565e-09
325, 2.4483454280733705e-09
350, 1.8187296689065332e-09
375, 1.3963623283938204e-09
400, 1.0831090246666068e-09
425, 8.5301003353708082e-10
450, 6.6627521593115182e-10
475, 5.3449281760803048e-10
500, 4.3280116475387868e-10
525, 3.6051723360617720e-10
550, 2.9595620662346797e-10
575, 2.4975241383705447e-10
600, 2.0961713550176611e-10
625, 1.7703346664347626e-10
650, 1.5037685005368506e-10
675, 1.3005942990686467e-10
700, 1.1256634033983998e-10
725, 9.8141191638027069e-11
750, 8.5059874893237957e-11
775, 7.4374830212501892e-11
800, 6.5473226574491316e-11
}\third

\pgfplotstableread[col sep = comma]{
N,p4
50,  4.3593303757674417e-07
60,  2.1765189755469900e-07
70,  6.9332320896238360e-08
80,  4.3378608247120799e-08
90,  1.8732830522879445e-08
100, 1.2601801333402425e-08
110, 6.5703713980456158e-09
120, 4.5975681777297400e-09
130, 2.9073876055775915e-09
140, 2.0600322308413781e-09
150, 1.4583114464534556e-09
160, 1.0298067368254729e-09
170, 8.2869366830260509e-10
180, 5.7158148096296488e-10
190, 4.5936260556296706e-10
200, 3.3347102473674383e-10
210, 2.6864277226253605e-10
220, 2.0381333415266667e-10
230, 1.6850435358721586e-10
240, 1.3070773531616133e-10
250, 1.1052313173939298e-10
260, 8.7153548557020113e-11
270, 7.4463280132467172e-11
280, 6.0058653845721428e-11
290, 5.1422920014204812e-11
300, 4.4601637520401058e-11
310, 3.6354108830624465e-11
320, 3.1790926621407892e-11
330, 2.6401640647857101e-11
340, 2.3403932065336699e-11
350, 1.9631208761844980e-11
360, 1.7674513792502034e-11
370, 1.4929496764237545e-11
380, 1.3374609672082517e-11
390, 1.1463438282234097e-11
400, 1.0372508871290771e-11
}\fourth

\pgfplotstableread[col sep = comma]{
N,p5
50, 3.2820164610802226e-08
60, 1.4757455296527604e-08
70, 3.4149564689376120e-09
80, 2.0802480963159199e-09
90, 6.9109544809088392e-10
100,4.6106404027231144e-10
110,1.9690481185038818e-10
120,1.3630054088402375e-10
130,7.7871296936638318e-11
140,5.1513223440456348e-11
150,3.3346919506921468e-11
160,2.2314335964845012e-11
170,1.6985076644637069e-11
180,1.1205826787012601e-11
190,8.3234559320128185e-12
200,5.8742073103572831e-12
210,4.5059149361893438e-12
220,3.4879506423406955e-12
230,2.5809833336410042e-12
240,1.8876793321975593e-12
250,1.7606944578782725e-12
260,1.4033756362165658e-12
270,1.6785934697090656e-12
280,1.1130914889186397e-12
290,9.1666518662227940e-13
300,7.6510969910045137e-13
}\fifth

\xdef\IDX{10}
\pgfplotstablesave[skip rows between index={0}{\IDX}]{\first}{first.dat};
\pgfplotstablesave[skip rows between index={0}{\IDX}]{\second}{second.dat};
\pgfplotstablesave[skip rows between index={0}{\IDX}]{\third}{third.dat};
\pgfplotstablesave[skip rows between index={0}{10}]{\fourth}{fourth.dat};
\pgfplotstablesave[skip rows between index={0}{10},
                   skip rows between index={20}{40}]{\fifth}{fifth.dat};
\definecolor{custom}{rgb}{0.93, 0.53, 0.18}
      
\addplot[no markers, black] table[
     x=N,
     y={create col/linear regression={y=p1}},
     no marks]
  {first.dat};
\xdef\slopeA{\pgfplotstableregressiona}
\pgfmathsetmacro{\msA}{-\slopeA}
  
\addplot[no markers, red] table[
     x=N,
     y={create col/linear regression={y=p2}}, 
     no marks]
  {second.dat};
\xdef\slopeB{\pgfplotstableregressiona}
\pgfmathsetmacro{\msB}{-\slopeB}

\addplot[no markers, custom] table[
     x=N,
     y={create col/linear regression={y=p3}}]
  {third.dat};
\xdef\slopeC{\pgfplotstableregressiona}
\pgfmathsetmacro{\msC}{-\slopeC}

\addplot[no markers, blue] table[
     x=N,
     y={create col/linear regression={y=p4}}]
  {fourth.dat};
\xdef\slopeD{\pgfplotstableregressiona}
\pgfmathsetmacro{\msD}{-\slopeD}

\addplot[no markers, green] table[
     x=N,
     y={create col/linear regression={y=p5}}]
  {fifth.dat};
\xdef\slopeE{\pgfplotstableregressiona}
\pgfmathsetmacro{\msE}{-\slopeE}

\pgfkeys{
   /pgf/number format/textnumber/.style={
     fixed zerofill,
     precision=2,
     },
  }

\addplot[black, only marks,mark size=1.pt, forget plot] table[x=N, y=p1] {\first};
\addlegendentry{P1Q2 (\pgfmathprintnumber[textnumber]{\msA})};
\addplot[red, only marks,mark size=1.pt, forget plot] table[x=N, y=p2] {\second};
\addlegendentry{P2Q2 (\pgfmathprintnumber[textnumber]{\msB})};
\addplot[custom, only marks,mark size=1.pt, forget plot] table[x=N, y=p3] {\third};
\addlegendentry{P3Q3 (\pgfmathprintnumber[textnumber]{\msC})};
\addplot[blue, only marks,mark size=1.pt, forget plot] table[x=N, y=p4] {\fourth};
\addlegendentry{P4Q4 (\pgfmathprintnumber[textnumber]{\msD})};
\addplot[green, only marks,mark size=1.pt, forget plot] table[x=N, y=p5] {\fifth};
\addlegendentry{P5Q5 (\pgfmathprintnumber[textnumber]{\msE})};
    \end{loglogaxis}
    
  \end{tikzpicture}

%% file: 3-rh/linf.tex
\begin{tikzpicture}
    \begin{loglogaxis}[
      width=11cm,
      height=6cm,
      xlabel={$N$},
      ylabel={$L_\infty$ error},
      legend cell align=left,
      legend style={at={(0.01, 0.01)}, anchor=south west},
      legend columns=1,
      legend style={font=\small},
      scale only axis,
      grid=both,
      xmin=20,
      xmax = 1000,
      ymin = 1e-13,
      ymax = 1e0
      ]
      
\pgfplotstableread[col sep = comma]{
N,p1
50,  6.4145064704033139e-02
75,  3.4044677989997663e-02
100, 2.1487688185113013e-02
125, 1.5198288080667588e-02
150, 1.0379999543147722e-02
175, 7.9706275663743442e-03
200, 6.6294492407773520e-03
225, 5.4692638895835755e-03
250, 4.5416523443657097e-03
275, 3.8596431065148762e-03
300, 3.4219743016977722e-03
325, 2.8378130076487373e-03
350, 2.5800857571599556e-03
375, 2.3489913380252414e-03
400, 2.1101100814140406e-03
425, 1.8933429329862639e-03
450, 1.7035997446221973e-03
475, 1.6034259037827070e-03
500, 1.4452058924507849e-03
525, 1.3512672134093906e-03
550, 1.2591018096564932e-03
575, 1.1732932296734111e-03
600, 1.0655801021377531e-03
625, 1.0124086100370411e-03
650, 9.6397001717107855e-04
675, 9.0639101136580358e-04
700, 8.3628572686367253e-04
725, 8.1565589283261186e-04
750, 7.6599041686298275e-04
775, 7.1004592425616408e-04
800, 6.8425728306409361e-04
}\first

\pgfplotstableread[col sep = comma]{
N,p2
50, 5.5101439190578527e-03
75, 1.7882903256910831e-03
100,6.4758801794329024e-04
125,3.2930518572221734e-04
150,1.7979926737965446e-04
175,1.1443192998283891e-04
200,7.4556163776245965e-05
225,5.2234577513188274e-05
250,3.6789468124298619e-05
275,2.8073450827448276e-05
300,2.0922498321862548e-05
325,1.7785345272691711e-05
350,1.3548122770351556e-05
375,1.1046237583967056e-05
400,9.1747294400090773e-06
425,7.4586920272845525e-06
450,6.4734004990851801e-06
475,5.6485712240927910e-06
500,5.0414239464169341e-06
525,4.2147210563281412e-06
550,3.6524971056550726e-06
575,3.1999741430532858e-06
600,2.9713619547622194e-06
625,2.6101559479352332e-06
650,2.4466925834798481e-06
675,2.1130625357801414e-06
700,1.9289241357522968e-06
725,1.7172292657852140e-06
750,1.5428958710650420e-06
775,1.4114273667598098e-06
800,1.3171635934172699e-06
}\second

\pgfplotstableread[col sep = comma]{
N,p3
50 ,4.4773817393888027e-04
75 ,9.2911377663873473e-05
100,3.2151832081397935e-05
125,1.7791478513079495e-05
150,5.7360436042097618e-06
175,2.9057240493179926e-06
200,2.1115808654936785e-06
225,1.6305526167048612e-06
250,7.8868581659952142e-07
275,5.3402460334361912e-07
300,3.6892924282216910e-07
325,3.6200659839225047e-07
350,2.1519424853178037e-07
375,2.0408739982713531e-07
400,1.4131803999228865e-07
425,1.0797820951258252e-07
450,7.7952745342058094e-08
475,7.5907605867620020e-08
500,4.7147867623942297e-08
525,4.5143597027230697e-08
550,3.5374466106352642e-08
575,3.1587819221412516e-08
600,2.2917829051172234e-08
625,2.5920088339570491e-08
650,1.7537793939315094e-08
675, 1.6302734689999099e-08
700, 1.2531752774380678e-08
725, 1.2124069512231550e-08
750, 1.0028139152407078e-08
775, 1.1231172963865887e-08
800, 8.9750626652840282e-09
}\third

\pgfplotstableread[col sep = comma]{
N,p4
50, 3.6315713893286805e-05
60, 1.7430230632997645e-05
70, 1.2310759090916612e-05
80, 4.7111168118840752e-06
90, 4.1432952767039843e-06
100,1.6516127540611159e-06
110,1.5930866166247171e-06
120,8.1690780467935653e-07
130,5.7362567018603983e-07
140,3.9603127877096966e-07
150,3.0358595037860425e-07
160,2.4531422954110482e-07
170,1.5805114053635272e-07
180,1.3702497747569709e-07
190,9.3246397048218910e-08
200,8.2461297512548271e-08
210,5.9501094579195168e-08
220,5.1324687388998314e-08
230,4.1332313127639964e-08
240,2.9395074374605201e-08
250,2.5365612055239239e-08
260,1.9391794336964807e-08
270,1.7416910413903253e-08
280,1.3041347224351796e-08
290,1.2133884952358898e-08
300,9.7069255758364115e-09
310,8.6008724586861263e-09
320,7.9051111767647342e-09
330,6.2073984841148899e-09
340,5.5914632368647688e-09
350,4.5112509838851622e-09
360,4.6497692918201494e-09
370,3.7312525291710585e-09
380,3.5367336836955587e-09
390,2.8289317555163507e-09
400,2.7118324252839443e-09
}\fourth

\pgfplotstableread[col sep = comma]{
N,p5
50, 6.2415503405688355e-06
60, 2.8607158383681108e-06
70, 5.7236049288977853e-07
80, 6.0769783466474436e-07
90, 1.6153978033006666e-07
100,1.5888759849480572e-07
110,4.8168877775367491e-08
120,4.9381868705289378e-08
130,2.9345109009426551e-08
140,1.7301206606257225e-08
150,1.2476584654308454e-08
160,7.6848109542559939e-09
170,7.3705555096725828e-09
180,4.1694137031811351e-09
190,3.4568825024727801e-09
200,2.7733840779475827e-09
210,1.9794378081083863e-09
220,1.2541113703079532e-09
230,1.1818269834984463e-09
240,7.7915457419308609e-10
250,7.1012318247909434e-10
260,4.2319170390214822e-10
270,4.6412762522152207e-10
280,2.7112870282230972e-10
290,2.8235606908122435e-10
300,3.0840464210202256e-10
}\fifth

\definecolor{custom}{rgb}{0.93, 0.53, 0.18}

\xdef\IDX{10}
\pgfplotstablesave[skip rows between index={0}{\IDX}]{\first}{first.dat};
\pgfplotstablesave[skip rows between index={0}{\IDX}]{\second}{second.dat};
\pgfplotstablesave[skip rows between index={0}{\IDX}]{\third}{third.dat};
\pgfplotstablesave[skip rows between index={0}{10}]{\fourth}{fourth.dat};
\pgfplotstablesave[skip rows between index={0}{10},
                   skip rows between index={20}{40}]{\fifth}{fifth.dat};
    
\addplot[no markers, black] table[
     x=N,
     y={create col/linear regression={y=p1}},
     no marks]
  {first.dat};
\xdef\slopeA{\pgfplotstableregressiona}
\pgfmathsetmacro{\msA}{-\slopeA}
  
\addplot[no markers, red] table[
     x=N,
     y={create col/linear regression={y=p2}}, 
     no marks]
  {second.dat};
\xdef\slopeB{\pgfplotstableregressiona}
\pgfmathsetmacro{\msB}{-\slopeB}

\addplot[no markers, custom] table[
     x=N,
     y={create col/linear regression={y=p3}}]
  {third.dat};
\xdef\slopeC{\pgfplotstableregressiona}
\pgfmathsetmacro{\msC}{-\slopeC}

\addplot[no markers, blue] table[
     x=N,
     y={create col/linear regression={y=p4}}]
  {fourth.dat};
\xdef\slopeD{\pgfplotstableregressiona}
\pgfmathsetmacro{\msD}{-\slopeD}

\addplot[no markers, green] table[
     x=N,
     y={create col/linear regression={y=p5}}]
  {fifth.dat};
\xdef\slopeE{\pgfplotstableregressiona}
\pgfmathsetmacro{\msE}{-\slopeE}
      
\pgfkeys{
   /pgf/number format/textnumber/.style={
     fixed zerofill,
     precision=2,
     },
  }
      
\addplot[black, only marks,mark size=1.pt, forget plot] table[x=N, y=p1] {\first};
\addlegendentry{P1Q2 (\pgfmathprintnumber[textnumber]{\msA})};
\addplot[red, only marks,mark size=1.pt, forget plot] table[x=N, y=p2] {\second};
\addlegendentry{P2Q2 (\pgfmathprintnumber[textnumber]{\msB})};
\addplot[custom, only marks,mark size=1.pt, forget plot] table[x=N, y=p3] {\third};
\addlegendentry{P3Q3 (\pgfmathprintnumber[textnumber]{\msC})};
\addplot[blue, only marks,mark size=1.pt, forget plot] table[x=N, y=p4] {\fourth};
\addlegendentry{P4Q4 (\pgfmathprintnumber[textnumber]{\msD})};
\addplot[green, only marks,mark size=1.pt, forget plot] table[x=N, y=p5] {\fifth};
\addlegendentry{P5Q5 (\pgfmathprintnumber[textnumber]{\msE})};

    \end{loglogaxis}
  \end{tikzpicture}

%% file: 3-rh/boundary1.tex
\pgfplotstableread[col sep = space]{3-rh/boundary1.txt}\bdryA

\begin{tikzpicture}
    \begin{axis}[
      width=6.5cm,
      height=4.5cm,
      xlabel={$\theta$},
      ylabel near ticks,
      ylabel={},
      legend columns=1,
      legend style={font=\small},
      scale only axis,
            grid=both,
            cycle list/Set1-5,
      xmin = 0,
      xmax = 3.15,
      ymin = -0.1,
      ymax =  1.1,
    xtick={0,0.78539, 1.5708, 2.356194, 3.14159},
    xticklabels={0, $\pi/4$, $\pi/2$, $3\pi/4$, $\pi$},
    ytick distance={0.25}
      ]

      \foreach \i in {0,...,35}{
      \pgfmathsetmacro\resultA{0+\i*2}
      \pgfmathsetmacro\resultB{1+\i*2}
      \addplot+[line width=2pt] table[x index=\resultA,y index=\resultB] {\bdryA};
      }

    \end{axis}
  \end{tikzpicture}

%% file: 3-rh/boundary2.tex
\pgfplotstableread[col sep = space]{3-rh/boundary2.txt}\bdryA

\begin{tikzpicture}
    \begin{axis}[
      width=6.5cm,
      height=4.5cm,
      xlabel={$\theta$},
      ylabel near ticks,
      ylabel={},
      legend cell align=left,
      legend style={at={(0.01, 0.01)}, anchor=south west},
      legend columns=1,
      legend style={font=\small},
      scale only axis,
            grid=both,
            cycle list/Set1-5,
      xmin = 0,
      xmax = 3.15,
      ymin = -0.1,
      ymax =  1.1,
    xtick={0,0.78539, 1.5708, 2.356194, 3.14159},
    xticklabels={0, $\pi/4$, $\pi/2$, $3\pi/4$, $\pi$},
    ytick distance={0.25}
      ]

      \foreach \i in {0,...,19}{
      \pgfmathsetmacro\resultA{0+\i*2}
      \pgfmathsetmacro\resultB{1+\i*2}
      \addplot+[line width=2pt] table[x index=\resultA,y index=\resultB] {\bdryA};
      }

    \end{axis}
  \end{tikzpicture}

%% file: 6-ringleb/overlaps1_ringleb.tex
  \begin{tikzpicture}
  \tikzset{
  reverseclip/.style={insert path={(-4,-4) rectangle (4cm,4cm)}}
}

\definecolor{beaublue}{rgb}{0.74, 0.83, 0.9};

\def\firstcircle{(7.7,-4.75) circle (11)}
\def\secondcircle{(1.,-15) circle (18)}

\begin{scope}
\clip[] (-1.1,1.85) -- (-0.15,2.95) --(0.1,2.97) --(1.25,3) -- (1.25,4.1) -- (-3,4.1) -- (-3,0) -- (-2.3,0) -- (-1.1,1.85);
\draw [ultra thick, draw=beaublue, fill=beaublue, fill opacity=0.2] (-1,2) rectangle node[text=black, opacity=1]{}  (0,3);
\draw [ultra thick, draw=beaublue, fill=beaublue, fill opacity=0.2] (-1,2) rectangle node[text=black, opacity=1]{}  (-2,3);
\draw [ultra thick, draw=beaublue, fill=beaublue, fill opacity=0.2] (-1,1) rectangle node[text=black, opacity=1]{}  (-2,2);
\draw [ultra thick, draw=beaublue, fill=beaublue, fill opacity=0.2] (-1,3) rectangle node[text=black, opacity=1]{}  (-2,4);
\draw [ultra thick, draw=beaublue, fill=beaublue, fill opacity=0.2] (0,3) rectangle node[text=black, opacity=1]{}  (-1,4);
\draw [ultra thick, draw=beaublue, fill=beaublue, fill opacity=0.2] (1,3) rectangle node[text=black, opacity=1]{}  (0,4);
\node at (-0.5,2.5) {$K_{i,j}$};
\end{scope}

\begin{scope}
\clip (-2.3,-0.1) rectangle (2.1,3.2);
\clip[] \firstcircle;
\clip[] \secondcircle;
\draw [thick] \firstcircle;
\draw [thick] \secondcircle;
\draw (-5,-5) grid (5,5);

\draw [ultra thick, draw=green!70!black, fill=green!70!black, fill opacity=0.2] (-2,1) rectangle  (1,4);

\node at (-0.5,2.5) {$K_{i,j}$};

\end{scope}

  \end{tikzpicture}

%% file: 6-ringleb/overlaps2_ringleb.tex
  \begin{tikzpicture}
  \tikzset{
  reverseclip/.style={insert path={(-4,-4) rectangle (4cm,4cm)}}
}

\definecolor{beaublue}{rgb}{0.74, 0.83, 0.9};

\def\firstcircle{(7.7,-4.75) circle (11)}
\def\secondcircle{(1.,-15) circle (18)}

\begin{scope}
\clip[] (-1.1,1.85) -- (-0.15,2.95) --(0.1,2.97) --(1.25,3) -- (1.25,4.1) -- (-3,4.1) -- (-3,0) -- (-2.3,0) -- (-1.1,1.85);
\draw [ultra thick, draw=beaublue, fill=beaublue, fill opacity=0.2] (-2,1) rectangle node[text=black, opacity=1]{$K$\scalebox{0.65}{${}_{i-1,j-1}$} }  (-1,2);
\end{scope}

\begin{scope}
\clip (-2.3,-0.1) rectangle (2.1,3.2);
\clip[] \firstcircle;
\clip[] \secondcircle;
\draw [thick] \firstcircle;
\draw [thick] \secondcircle;
\draw (-5,-5) grid (5,5);

\draw [ultra thick, draw=green!70!black, fill=green!70!black, fill opacity=0.2] (-2,1) rectangle  (0,2);

\node at (-1.5,1.5) {$K$\scalebox{0.65}{${}_{i-1,j-1}$}};

\end{scope}

  \end{tikzpicture}

%% file: 6-ringleb/l1.tex
\begin{tikzpicture}
    \begin{loglogaxis}[
      width=11cm,
      height=6cm,
      xlabel={$N$},
      ylabel={$L_1$ error},
      legend cell align=left,
      legend style={at={(0.01, 0.01)}, anchor=south west},
      legend columns=1,
      legend style={font=\small},
      scale only axis,
      grid=both,
      xmin=1,
      xmax = 1000,
      ymin = 1e-10,
      ymax = 1e-1
      ]
      
\pgfplotstableread[col sep = comma]{
N,p1
 10,3.2805172682627381e-03
 12,2.2212376525927698e-03
 14,2.4188033851603946e-03
 16,1.0990329228792210e-03
 18,7.5388995375913655e-04
 20,6.1222672039309555e-04
 22,4.7048966779712130e-04
 24,3.7374946868527710e-04
 25,3.1353976893733855e-04
 27,2.7298322416276650e-04
 30,2.2278049985674316e-04
 32,1.9296646831678160e-04
 35,1.6769376833568716e-04
 38,1.4258199731665336e-04
 40,1.2968328287831724e-04
 43,1.1464753743052798e-04
 48,8.8816368028053560e-05
 51,8.0114594922438284e-05
 57,6.6529207916069127e-05
 60,5.9359367718543263e-05
 65,5.0913906805850121e-05
 70,4.4226259040320349e-05
 75,3.8523874502999804e-05
 80,3.4402075180498688e-05
 86,2.9468204959763645e-05
 91,2.6399197481568617e-05
 96,2.3701924384305831e-05
 99,2.2649543848092506e-05
102,2.0944389214146889e-05
107,1.9114083424763768e-05
115,1.6702582502821833e-05
120,1.5436319988848762e-05
125,1.4236844685566117e-05
130,1.3018214822580592e-05
135,1.2072396028662921e-05
140,1.1263823252332707e-05
145,1.0517015066749465e-05
150,9.8351557037362489e-06
155,9.2469560929437744e-06
158,8.8593234188056944e-06
166,8.0207776111162702e-06
171,7.5579995359283659e-06
175,7.2045921619359107e-06
180,6.8179703292910151e-06
185,6.4610954820267868e-06
190,6.1366974569310596e-06
195,5.8365422619745273e-06
201,5.4725596772911660e-06
206,5.2200712850594652e-06
209,5.0707904459131851e-06
214,4.8330688859930733e-06
222,4.4814949923088323e-06
230,4.1869328121337293e-06
235,4.0129452825175548e-06
240,3.8243188511679639e-06
245,3.6750243722965512e-06
250,3.5294160149448308e-06
255,3.3945320405634726e-06
260,3.2646497209569342e-06
265,3.1424533266822646e-06
270,3.0295084835597646e-06
285,2.7069790137617034e-06
290,2.6119165165420114e-06
295,2.5260423582695396e-06
300,2.4464613638745576e-06
305,2.3649136163594683e-06
310,2.2893248179910730e-06
314,2.2286728012614937e-06
319,2.1595726995027198e-06
324,2.0939233212803219e-06
329,2.0220759585809294e-06
345,1.8407049841046794e-06
350,1.7840658749849353e-06
355,1.7333831128884598e-06
360,1.6850558979021580e-06
365,1.6384975251063707e-06
370,1.5952554546552697e-06
375,1.5536833229056532e-06
380,1.5119869849709123e-06
383,1.4857231423718563e-06
388,1.4478923957068594e-06
400,1.3605886736075539e-06
405,1.3266893740669153e-06
410,1.2949434227356789e-06
415,1.2636648853717605e-06
420,1.2337723724519733e-06
425,1.2019841522891542e-06
429,1.1814996843124973e-06
434,1.1536122440865505e-06
439,1.1253688676918983e-06
444,1.0998878396696057e-06
456,1.0421625716611334e-06
460,1.0234220279746814e-06
465,1.0016881074705626e-06
470,9.8003511879462284e-07
475,9.5949840219227543e-07
480,9.3915897763730663e-07
485,9.2001879328847493e-07
490,9.0144366064280498e-07
495,8.8272759347720954e-07
}\first

\pgfplotstableread[col sep = comma]{
N,p2
 10,3.5583918200953726e-04
 12,1.4280609564686279e-04
 14,9.6261246869015619e-05
 16,4.5616307849762241e-05
 18,3.4470319848271411e-05
 20,2.3019636669290884e-05
 22,2.2897764315860129e-05
 24,1.7257544480502654e-05
 25,1.2193661957380674e-05
 27,9.4833541016063951e-06
 30,7.6943773176478999e-06
 32,5.6393677905924453e-06
 35,4.9711594012909206e-06
 38,3.8149976875174936e-06
 40,3.3355553218271572e-06
 43,2.7769348042274632e-06
 48,1.7555650703255187e-06
 51,1.5240324121049992e-06
 57,1.1620730238767239e-06
 60,9.5711367813618741e-07
 65,7.6102240178886121e-07
 70,6.1594407139261861e-07
 75,4.9821271768142755e-07
 80,4.2051785945255893e-07
 86,3.2821319044055822e-07
 91,2.7697761178210454e-07
 96,2.3521130445486351e-07
 99,2.2052754757423865e-07
102,1.9491397871757981e-07
107,1.6769120255539051e-07
115,1.3846687412806650e-07
120,1.2348982421237914e-07
125,1.0780039602754855e-07
130,9.3466680866205303e-08
135,8.3541968551295071e-08
140,7.5786471075131745e-08
145,6.8493047053730518e-08
150,6.1712103434194561e-08
155,5.6094302963036157e-08
158,5.2073521489128886e-08
166,4.4712752789941529e-08
171,4.0807331879283401e-08
175,3.8257833673178138e-08
180,3.5391498252046508e-08
185,3.2595217002169762e-08
190,3.0064114255348339e-08
195,2.7897307138432577e-08
201,2.4746793457784426e-08
206,2.3159197625725839e-08
209,2.2581882593646253e-08
214,2.0915087529924455e-08
222,1.8534201110578103e-08
230,1.6470577493343209e-08
235,1.5449212615904867e-08
240,1.4330261380137812e-08
245,1.3601822113484142e-08
250,1.2779891771481249e-08
255,1.2102868929587088e-08
260,1.1452644337663559e-08
265,1.0786633791219536e-08
270,1.0272072369254647e-08
285,8.6884351001401156e-09
290,8.2181720498563210e-09
295,7.8313572735170754e-09
300,7.4559546946650485e-09
305,7.0898736539229251e-09
310,6.7392997251678871e-09
314,6.4729602538811927e-09
319,6.1525417086734057e-09
324,5.8797006503713090e-09
329,5.5743360156034267e-09
345,4.8649787673038590e-09
350,4.6480401842476440e-09
355,4.4326777197811109e-09
360,4.2644646995263649e-09
365,4.0868996944014364e-09
370,3.9341055957769655e-09
375,3.7744120480348606e-09
380,3.6182402550943161e-09
383,3.5145401007314834e-09
388,3.3909027100770153e-09
400,3.1010425121785068e-09
405,2.9742885880635596e-09
410,2.8773721572826419e-09
415,2.7631297139043466e-09
420,2.6721789036262236e-09
425,2.5691772120734711e-09
429,2.5075347650970977e-09
434,2.4089968744845652e-09
439,2.3320261011436106e-09
444,2.2446701084894554e-09
456,2.0658885061149390e-09
460,1.9920899774918728e-09
465,1.9311223794898464e-09
470,1.8699147559528165e-09
475,1.8154442166854204e-09
480,1.7576899695990604e-09
485,1.7090742191313133e-09
490,1.6584706710189261e-09
495,1.6057296845307952e-09
}\second

\pgfplotstableread[col sep = comma]{
N,p3
 10,5.5851592994334569e-05
 12,1.5073544646451032e-05
 14,1.1231641461284709e-05
 16,4.1768620784299493e-06
 18,3.2366735657250445e-06
 20,1.8627184353563951e-06
 22,1.3842009620229028e-06
 24,1.1333040343504741e-06
 25,8.0119517818959650e-07
 27,5.1351938159014060e-07
 30,3.7578804452627873e-07
 32,2.8854696331875075e-07
 35,2.1094087163706030e-07
 38,1.4318517081222770e-07
 40,1.1632639832550120e-07
 43,8.8573238963638611e-08
 48,5.3691145745952347e-08
 51,4.0595455314727986e-08
 57,2.8293401745511217e-08
 60,2.3217518082659528e-08
65,1.6546902055410362e-08
70,1.1772503172814581e-08
75,9.0961708660222137e-09
80,7.1980977138575995e-09
86,5.2484566847735307e-09
91,4.1504277826441194e-09
96,3.4886859588415833e-09
99,3.1071671673565567e-09
102,2.6599183907258573e-09
107,2.1322788350574476e-09
115,1.7182244834549887e-09
120,1.4599311347686425e-09
125,1.1579167154085747e-09
130,1.0256723007244445e-09
135,9.0184683757809523e-10
140,7.3550700302109495e-10
}\third

\pgfplotstableread[col sep = comma]{
N,p4
 10,5.8172540219671763e-06
 12,1.9333130767431937e-06
 14,1.5019085611412127e-06
 16,4.3694435318850677e-07
 18,3.0310397194612835e-07
 20,1.9816451095332491e-07
 22,1.1240745809135710e-07
 24,8.5400620167334306e-08
 25,6.2375061192291855e-08
 27,3.2363079800546891e-08
 30,2.2272397313866937e-08
 32,1.6797851318085030e-08
 35,1.0803359519963620e-08
 38,7.1246981898417488e-09
 40,5.2255562670906441e-09
 43,3.9448141276239590e-09
 48,2.1686588497705205e-09
 51,1.5330507827011392e-09
 57,1.1035326630116758e-09
}\fourth

\pgfplotstableread[col sep = comma]{
N,p5
 10,9.0870219617537606e-07
 12,3.2805876188560290e-07
 14,2.5983109310093029e-07
 16,4.7986186432514587e-08
 18,3.5656878301039374e-08
 20,2.1420731202190825e-08
 22,1.1667361711474012e-08
 24,6.7185257467256095e-09
 25,5.0523332562062500e-09
 27,3.3172708427717640e-09
 30,1.3952867784602353e-09
 32,1.0118240286749149e-09
}\fifth
\xdef\IDX{4}
\pgfplotstablesave[skip rows between index={0}{5}]{\first}{first.dat};
\pgfplotstablesave[skip rows between index={0}{5}]{\second}{second.dat};
\pgfplotstablesave[skip rows between index={0}{5}]{\third}{third.dat};
\pgfplotstablesave[skip rows between index={0}{5}]{\fourth}{fourth.dat};
\pgfplotstablesave[skip rows between index={0}{5}]{\fifth}{fifth.dat};
\definecolor{custom}{rgb}{0.93, 0.53, 0.18}

\addplot[no markers, black] table[
     x=N,
     y={create col/linear regression={y=p1}},
     no marks]
  {first.dat};
\xdef\slopeA{\pgfplotstableregressiona}
\pgfmathsetmacro{\msA}{-\slopeA}
  
\addplot[no markers, red] table[
     x=N,
     y={create col/linear regression={y=p2}}, 
     no marks]
  {second.dat};
\xdef\slopeB{\pgfplotstableregressiona}
\pgfmathsetmacro{\msB}{-\slopeB}

\addplot[no markers, custom] table[
     x=N,
     y={create col/linear regression={y=p3}}]
  {third.dat};
\xdef\slopeC{\pgfplotstableregressiona}
\pgfmathsetmacro{\msC}{-\slopeC}

\addplot[no markers, blue] table[
     x=N,
     y={create col/linear regression={y=p4}}]
  {fourth.dat};
\xdef\slopeD{\pgfplotstableregressiona}
\pgfmathsetmacro{\msD}{-\slopeD}

\addplot[no markers, green] table[
     x=N,
     y={create col/linear regression={y=p5}}]
  {fifth.dat};
\xdef\slopeE{\pgfplotstableregressiona}
\pgfmathsetmacro{\msE}{-\slopeE}

\pgfkeys{
  /pgf/number format/textnumber/.style={
     fixed zerofill,
     precision=2,
     },
  }

\addplot[black, only marks,mark size=1.pt, forget plot] table[x=N, y=p1] {\first};
\addlegendentry{P1Q2 (\pgfmathprintnumber[textnumber]{\msA})};
\addplot[red, only marks,mark size=1.pt, forget plot] table[x=N, y=p2] {\second};
\addlegendentry{P2Q2 (\pgfmathprintnumber[textnumber]{\msB})};
\addplot[custom, only marks,mark size=1.pt, forget plot] table[x=N, y=p3] {\third};
\addlegendentry{P3Q3 (\pgfmathprintnumber[textnumber]{\msC})};
\addplot[blue, only marks,mark size=1.pt, forget plot] table[x=N, y=p4] {\fourth};
\addlegendentry{P4Q4 (\pgfmathprintnumber[textnumber]{\msD})};
\addplot[green, only marks,mark size=1.pt, forget plot] table[x=N, y=p5] {\fifth};
\addlegendentry{P5Q5 (\pgfmathprintnumber[textnumber]{\msE})};
    \end{loglogaxis}
  \end{tikzpicture}

%% file: 6-ringleb/linf.tex
\begin{tikzpicture}
    \begin{loglogaxis}[
      width=11cm,
      height=6cm,
      xlabel={$N$},
      ylabel={$L_{\infty}$ error},
      legend cell align=left,
      legend style={at={(0.01, 0.01)}, anchor=south west},
      legend columns=1,
      legend style={font=\small},
      scale only axis,
      grid=both,
      xmin=1,
      xmax = 1000,
      ymin = 1e-10,
      ymax = 1e-1
      ]
      
\pgfplotstableread[col sep = comma]{
N,p1
 10,1.4207859522917254e-02
 12,1.8787371931656915e-02
 14,1.5121395326878284e-02
 16,8.4495250954309542e-03
 18,6.3076218974477127e-03
 20,6.1286916250821166e-03
 22,4.6489576399864063e-03
 24,3.7767323447353007e-03
 25,3.1810319856222513e-03
 27,2.7927571457028177e-03
 30,2.3312665514216846e-03
 32,2.1906009516616809e-03
 35,2.2276890022234763e-03
 38,1.7063011306507336e-03
 40,1.7598092169057722e-03
 43,1.4192120582046108e-03
 48,1.1492811094181432e-03
 51,1.1235033732658328e-03
 57,9.0402850473980667e-04
 60,8.5724573251311220e-04
 65,7.1369734812720687e-04
 70,7.1013137502284884e-04
 75,6.3158971653365370e-04
 80,5.1399216359027111e-04
 86,4.9607168594978379e-04
 91,4.3737491750961421e-04
 96,4.1805057301502746e-04
 99,3.7950860032309741e-04
102,3.7603704098465585e-04
107,3.6583699827941185e-04
115,3.3053975549768300e-04
120,2.9930707331005113e-04
125,2.8237575210676713e-04
130,2.4941051152105320e-04
135,2.5718536409802262e-04
140,2.2539710736579899e-04
145,2.2396165575977456e-04
150,1.9521078996020957e-04
155,2.0319596723927802e-04
158,1.9847888695911919e-04
166,1.8588158874055516e-04
171,1.7551893811107444e-04
175,1.6712604555402244e-04
180,1.6392640139251125e-04
185,1.5696776410323299e-04
190,1.4620820772692777e-04
195,1.4590130973202164e-04
201,1.3788683478810349e-04
206,1.3116071294150267e-04
209,1.2813485078888576e-04
214,1.2143803458763980e-04
222,1.1768724789684715e-04
230,1.1443036739822166e-04
235,1.0737003506489096e-04
240,1.0662992208132493e-04
245,1.0229603444422697e-04
250,9.8343414077217695e-05
255,9.3374458093098411e-05
260,9.0287228146967990e-05
265,8.6225738877954150e-05
270,8.5824943577783586e-05
285,7.8758936854717021e-05
290,7.9045888635831396e-05
295,7.5099960987423664e-05
300,7.4743602496929462e-05
305,7.3095677703083339e-05
310,6.9211439997718394e-05
314,7.0116861475955083e-05
319,6.8864122841461040e-05
324,6.8068070440352457e-05
329,6.3638466351711998e-05
345,6.1307269115018492e-05
350,5.9053355626947912e-05
355,5.6753079152982444e-05
360,5.7065025415248272e-05
365,5.6039218622627729e-05
370,5.4967152404716835e-05
375,5.4169145878968372e-05
380,5.0935873939450538e-05
383,5.1847066564691957e-05
388,5.0944854875112355e-05
400,4.8933252958072160e-05
405,4.7506406191222261e-05
410,4.7543766338309368e-05
415,4.6687201061490491e-05
420,4.4676889040862200e-05
425,4.4255816754601440e-05
429,4.4424413484511760e-05
434,4.2946272639232674e-05
439,4.2823481031684096e-05
444,4.1388897923488877e-05
456,3.9513236738342172e-05
460,3.9509923057234886e-05
465,3.8913326421363692e-05
470,3.8334615652946269e-05
475,3.8002688625971004e-05
480,3.6662830810274549e-05
485,3.6462816796256625e-05
490,3.6389862346886126e-05
495,3.4910604257976097e-05
}\first

\pgfplotstableread[col sep = comma]{
N,p2
 10,3.0579542218168720e-03
 12,1.5631310426794665e-03
 14,1.2968604842921172e-03
 16,8.5764721203041017e-04
 18,7.9196561846461844e-04
 20,5.3982984384459698e-04
 22,4.6838352109657144e-04
 24,4.0962510672959418e-04
 25,3.3107736897575180e-04
 27,2.0769726990532700e-04
 30,1.9508162774517501e-04
 32,1.3638450223563137e-04
 35,1.3147057260043393e-04
 38,1.2615568370322183e-04
 40,1.0368888558998091e-04
 43,9.3943046878486314e-05
 48,6.8231628346127060e-05
 51,4.0040802477991910e-05
 57,4.0377735106345014e-05
 60,2.5576349776179619e-05
 65,2.4566259425040649e-05
 70,1.7588361211462455e-05
 75,1.6063472844329851e-05
 80,1.5183355240422358e-05
 86,1.2746438375677016e-05
 91,1.1069071673586173e-05
 96,9.6532988028696920e-06
 99,7.8255685929118357e-06
102,1.0379675471350325e-05
107,9.1323618588434741e-06
115,5.5492646163779469e-06
120,5.3745732352394882e-06
125,4.8003645697347608e-06
130,4.2818229505359184e-06
135,3.9369728689342409e-06
140,3.6587206202476352e-06
145,3.3837320487783984e-06
150,3.1148139933989683e-06
155,2.9037651645413831e-06
158,2.1943647211530859e-06
166,2.9006789590280135e-06
171,2.7121006501706901e-06
175,2.2678873424775148e-06
180,2.1394420598452868e-06
185,2.0040963847334226e-06
190,1.8938762242903806e-06
195,1.7948441557447126e-06
201,1.0681036682891332e-06
206,1.2421539847817797e-06
209,1.3428258286429795e-06
214,1.2722345780202815e-06
222,9.2259480966028207e-07
230,8.4612039585874044e-07
235,7.8360990662940111e-07
240,6.8388057439250360e-07
245,6.2712663795316814e-07
250,6.5282116534781665e-07
255,6.1687801422571908e-07
260,5.8272138991899425e-07
265,5.5496507322949640e-07
270,5.6125285397001790e-07
285,4.8758089077338695e-07
290,4.7307039519051131e-07
295,4.4249874231461206e-07
300,4.3184417020114552e-07
305,4.2175334868943537e-07
310,4.1200754363757142e-07
314,3.5849209767402357e-07
319,3.4923042591117337e-07
324,3.5115102592353509e-07
329,3.3478905858697061e-07
345,3.2627176482780129e-07
350,3.2034510244471903e-07
355,3.0175680043598163e-07
360,2.9729232575714803e-07
365,2.9237600251974527e-07
370,2.8792334338945125e-07
375,2.7445674544157583e-07
380,2.6912692019376294e-07
383,2.3511997904179793e-07
388,2.2216161432719872e-07
400,2.4660528863407194e-07
405,2.3706279994595292e-07
410,2.3417018424343183e-07
415,2.2603412941357703e-07
420,2.2342159822219543e-07
425,2.1804561922689203e-07
429,2.0066546135844021e-07
434,1.9152461638860530e-07
439,1.9012056495792962e-07
444,1.8243700339315438e-07
456,1.4123046443437204e-07
460,1.3170571422360666e-07
465,1.2839825380694947e-07
470,1.2617477973542179e-07
475,1.2546445482541912e-07
480,1.2622741418777395e-07
485,1.1848055792018641e-07
490,1.1795560916905146e-07
495,1.1711296210581423e-07
}\second

\pgfplotstableread[col sep = comma]{
N,p3
 10,5.0090824181758808e-04
 12,2.7837548014475466e-04
 14,2.5338808274777413e-04
 16,7.5196052143011727e-05
 18,9.6786684444682969e-05
 20,3.7422141563103750e-05
 22,5.1238350655791187e-05
 24,4.9965652888839251e-05
 25,2.0018283530043313e-05
 27,1.6078862531521132e-05
 30,1.0584199976904074e-05
 32,1.2314937095081468e-05
 35,7.9217511516072037e-06
 38,6.3244261377937150e-06
 40,4.4091540459811895e-06
 43,5.0526765265068718e-06
 48,3.4763852447827759e-06
 51,2.4294587933626488e-06
 57,1.4643596034646933e-06
 60,1.1799269983292504e-06
  65,7.3365338015829451e-07
 70,7.3408996703605567e-07
 75,5.4534071403100626e-07
 80,4.3966986496002392e-07
 86,3.1587763871243624e-07
 91,2.7114086342816535e-07
 96,2.3874233101572884e-07
 99,1.9909082482616469e-07
 102,2.3166615314984540e-07
 107,1.7353296333233459e-07
 115,1.1760096518820262e-07
 120,1.0666225203070923e-07
 125,8.2781106192086895e-08
 130,6.3459557941847322e-08
 135,6.2388720745154558e-08
 140,4.4623735595550329e-08
}\third

\pgfplotstableread[col sep = comma]{
N,p4
 10,1.2504599593510246e-04
 12,5.8682098434403862e-05
 14,3.4829995676388670e-05
 16,1.0862659241106520e-05
 18,1.2829320920504372e-05
 20,7.4674289305409758e-06
 22,4.9547062095367522e-06
 24,3.7328882876996872e-06
 25,2.4853114852785296e-06
 27,1.9924142147242918e-06
 30,1.2119395458309867e-06
 32,1.1305495081215255e-06
 35,7.3924460164409567e-07
 38,4.7911018641322300e-07
 40,3.3852277059320812e-07
 43,2.2428218193759619e-07
 48,1.4514228963147247e-07
 51,1.4085893118576109e-07
 57,8.1295700726791154e-08
}\fourth

\pgfplotstableread[col sep = comma]{
N,p5
 10,1.9612677539804224e-05
 12,1.0173777688082808e-05
 14,6.0590026298434907e-06
 16,1.9566881640509237e-06
 18,1.5425124916257360e-06
 20,9.4439728248740806e-07
 22,5.8758449217322095e-07
 24,4.5712815321952149e-07
 25,3.7733488111868496e-07
 27,1.8056509698549661e-07
 30,1.0036762299048263e-07
 32,9.4415027018790454e-08
}\fifth
\xdef\IDX{4}
\pgfplotstablesave[skip rows between index={0}{5}]{\first}{first.dat};
\pgfplotstablesave[skip rows between index={0}{5}]{\second}{second.dat};
\pgfplotstablesave[skip rows between index={0}{5}]{\third}{third.dat};
\pgfplotstablesave[skip rows between index={0}{5}]{\fourth}{fourth.dat};
\pgfplotstablesave[skip rows between index={0}{5}]{\fifth}{fifth.dat};
\definecolor{custom}{rgb}{0.93, 0.53, 0.18}

\addplot[no markers, black] table[
     x=N,
     y={create col/linear regression={y=p1}},
     no marks]
  {first.dat};
\xdef\slopeA{\pgfplotstableregressiona}
\pgfmathsetmacro{\msA}{-\slopeA}
  
\addplot[no markers, red] table[
     x=N,
     y={create col/linear regression={y=p2}}, 
     no marks]
  {second.dat};
\xdef\slopeB{\pgfplotstableregressiona}
\pgfmathsetmacro{\msB}{-\slopeB}

\addplot[no markers, custom] table[
     x=N,
     y={create col/linear regression={y=p3}}]
  {third.dat};
\xdef\slopeC{\pgfplotstableregressiona}
\pgfmathsetmacro{\msC}{-\slopeC}

\addplot[no markers, blue] table[
     x=N,
     y={create col/linear regression={y=p4}}]
  {fourth.dat};
\xdef\slopeD{\pgfplotstableregressiona}
\pgfmathsetmacro{\msD}{-\slopeD}

\addplot[no markers, green] table[
     x=N,
     y={create col/linear regression={y=p5}}]
  {fifth.dat};
\xdef\slopeE{\pgfplotstableregressiona}
\pgfmathsetmacro{\msE}{-\slopeE}

\pgfkeys{
  /pgf/number format/textnumber/.style={
     fixed zerofill,
     precision=2,
     },
  }

\addplot[black, only marks,mark size=1.pt, forget plot] table[x=N, y=p1] {\first};
\addlegendentry{P1Q2 (\pgfmathprintnumber[textnumber]{\msA})};
\addplot[red, only marks,mark size=1.pt, forget plot] table[x=N, y=p2] {\second};
\addlegendentry{P2Q2 (\pgfmathprintnumber[textnumber]{\msB})};
\addplot[custom, only marks,mark size=1.pt, forget plot] table[x=N, y=p3] {\third};
\addlegendentry{P3Q3 (\pgfmathprintnumber[textnumber]{\msC})};
\addplot[blue, only marks,mark size=1.pt, forget plot] table[x=N, y=p4] {\fourth};
\addlegendentry{P4Q4 (\pgfmathprintnumber[textnumber]{\msD})};
\addplot[green, only marks,mark size=1.pt, forget plot] table[x=N, y=p5] {\fifth};
\addlegendentry{P5Q5 (\pgfmathprintnumber[textnumber]{\msE})};
    \end{loglogaxis}
  \end{tikzpicture}

%% file: 9-ap/boundary1.tex
\pgfplotstableread[col sep = comma]{9-ap/boundary1.txt}\bdryA

\begin{tikzpicture}
    \begin{axis}[
      width=6.5cm,
      height=4.5cm,
      xlabel={$\theta$},
      ylabel near ticks,
      ylabel={$(P - 1/\gamma)/10^{-4}$},
      legend columns=1,
      legend style={font=\small},
      scale only axis,
            grid=both,
            cycle list/Set1-5,
      xmin = -3.3,
      xmax = 3.3,
      ymin = -0.09,
      ymax =  0.09,
    xtick={-3.14159,-1.5708,0, 1.5708,  3.14159},
    xticklabels={$-\pi$,$-\pi/2$,0,  $\pi/2$, $\pi$},
    xlabel shift = -6pt,
    ytick distance={0.04}
      ]

  \pgfplotstablegetcolsof{\bdryA}
  \pgfmathparse{\pgfplotsretval/2-1}
      \foreach \i in {0,...,41}{
      \pgfmathsetmacro\resultA{0+\i*2}
      \pgfmathsetmacro\resultB{1+\i*2}
      \addplot+[line width=2pt] table[x index=\resultA,y index=\resultB] {\bdryA};
      }

    \end{axis}
  \end{tikzpicture}